\RequirePackage{fix-cm}
\documentclass[11pt,a4paper]{article}

\usepackage[a4paper,top=2.20cm,bottom=2.20cm,left=2.35cm,right=2.35cm,headheight=15pt,headsep=0.38cm,footskip=0.78cm]{geometry}
\usepackage{fontspec}
\usepackage{amsmath,amssymb,amsthm,mathtools}
\usepackage{bm,mathrsfs}
\usepackage{enumitem}
\usepackage{array,booktabs,longtable,tabularx,multirow}
\usepackage{graphicx,tikz-cd}
\usepackage{microtype,setspace,fancyhdr,titlesec}
\usepackage[numbers,sort&compress]{natbib}
\usepackage{aliascnt}
\usepackage[colorlinks=false,pdfborder={0 0 0},unicode=true]{hyperref}
\usepackage[nameinlink,noabbrev]{cleveref}
\usepackage{bookmark}
\usepackage{emptypage}
\usepackage{xcolor}

\newcommand{\zihao}[1]{}

\setlist{nosep,leftmargin=2.0em,topsep=2pt,partopsep=0pt,parsep=0pt,itemsep=0pt}
\allowdisplaybreaks[4]

\titleformat{\section}{\bfseries\large}{\thesection}{0.6em}{}
\titleformat{\subsection}{\bfseries\normalsize}{\thesubsection}{0.6em}{}
\titleformat{\subsubsection}{\bfseries\normalsize}{\thesubsubsection}{0.6em}{}
\titlespacing*{\section}{0pt}{1.35ex plus .35ex minus .2ex}{0.55ex plus .15ex}
\titlespacing*{\subsection}{0pt}{1.0ex plus .25ex minus .15ex}{0.35ex plus .1ex}
\titlespacing*{\subsubsection}{0pt}{0.8ex plus .2ex minus .1ex}{0.25ex plus .1ex}

\newtheoremstyle{compactplain}{4pt}{4pt}{\itshape}{} {\bfseries}{.}{0.5em}{}
\newtheoremstyle{compactdefinition}{4pt}{4pt}{}{} {\bfseries}{.}{0.5em}{}
\newtheoremstyle{compactremark}{4pt}{4pt}{}{} {\itshape}{.}{0.5em}{}
\theoremstyle{compactplain}
\newtheorem{theorem}{Theorem}[section]
\newaliascnt{proposition}{theorem}
\newtheorem{proposition}[proposition]{Proposition}
\aliascntresetthe{proposition}
\newaliascnt{lemma}{theorem}
\newtheorem{lemma}[lemma]{Lemma}
\aliascntresetthe{lemma}
\newaliascnt{corollary}{theorem}
\newtheorem{corollary}[corollary]{Corollary}
\aliascntresetthe{corollary}
\newaliascnt{claim}{theorem}

\aliascntresetthe{claim}
\theoremstyle{compactdefinition}
\newaliascnt{definition}{theorem}
\newtheorem{definition}[definition]{Definition}
\aliascntresetthe{definition}
\newaliascnt{notation}{theorem}
\newtheorem{notation}[notation]{Notation}
\aliascntresetthe{notation}
\crefname{notation}{Notation}{Notation}
\Crefname{notation}{Notation}{Notation}
\newaliascnt{example}{theorem}

\aliascntresetthe{example}
\theoremstyle{compactremark}
\newaliascnt{remark}{theorem}
\newtheorem{remark}[remark]{Remark}
\aliascntresetthe{remark}
\numberwithin{equation}{section}

\crefname{theorem}{Theorem}{Theorems}
\Crefname{theorem}{Theorem}{Theorems}
\crefname{proposition}{Proposition}{Propositions}
\Crefname{proposition}{Proposition}{Propositions}
\crefname{lemma}{Lemma}{Lemmas}
\Crefname{lemma}{Lemma}{Lemmas}
\crefname{corollary}{Corollary}{Corollaries}
\Crefname{corollary}{Corollary}{Corollaries}
\crefname{claim}{Claim}{Claims}
\Crefname{claim}{Claim}{Claims}
\crefname{definition}{Definition}{Definitions}
\Crefname{definition}{Definition}{Definitions}
\crefname{example}{Example}{Examples}
\Crefname{example}{Example}{Examples}
\crefname{remark}{Remark}{Remarks}
\Crefname{remark}{Remark}{Remarks}

\DeclareMathOperator{\im}{im}
\DeclareMathOperator{\coker}{coker}
\DeclareMathOperator{\Map}{Map}
\DeclareMathOperator{\Deck}{Deck}
\DeclareMathOperator{\id}{id}
\newcommand{\RP}{\mathbb{RP}}
\newcommand{\ZZ}{\mathbb Z}
\newcommand{\RR}{\mathbb R}

\newcommand{\AuthorName}{Xinyu Liu and Xuezhi Zhao}

\fancypagestyle{paper}{%
  \fancyhf{}
  \fancyhead[L]{\small Self-maps of $\RP^3\#\RP^3$}
  \fancyhead[R]{\zihao{-5}\AuthorName}
  
  \fancyfoot[C]{\thepage}
}
\fancypagestyle{plain}{%
  \fancyhf{}
  
  \fancyfoot[C]{\thepage}
}
\hypersetup{pdftitle={Pointed and Free Homotopy Classification of Self-Maps of RP3\#RP3},pdfauthor={Xinyu Liu and Xuezhi Zhao},pdfsubject={Algebraic topology},pdfkeywords={connected sum, pointed homotopy, free homotopy, infinite dihedral group, equivariant homotopy, mapping spaces}}

\theoremstyle{compactplain}
\newtheorem{maintheorem}{Theorem}

\crefname{maintheorem}{Theorem}{Theorems}
\Crefname{maintheorem}{Theorem}{Theorems}
\title{Pointed and Free Homotopy Classification of Self-Maps of $\mathbb{RP}^3\#\mathbb{RP}^3$}
\author{
  Xinyu Liu\\
  {\small Capital Normal University}\\
  {\small\texttt{2240502058@cnu.edu.cn}}
  \and
  Xuezhi Zhao\\
  {\small Capital Normal University}\\
  {\small\texttt{zhaoxve@mail.cnu.edu.cn}}
}
\date{}

\begin{document}
\maketitle
\begin{abstract}
We classify self-maps of $M=\mathbb{RP}^3\#\mathbb{RP}^3$ up to pointed and free homotopy. Using the universal cover $S^2\times\mathbb R$, we reduce the problem to paths in sphere mapping spaces with even or odd endpoint conditions and compute secondary invariants via the evaluation fibration, Whitehead products, and double cosets. Maps inducing the trivial fundamental-group endomorphism form four families, each in bijection with $\mathbb Z$ after choosing a reference class. Writing $G=\langle a,b\mid a^2=b^2=1\rangle$ and $t=ab$, every other realizable endomorphism has the form $t\mapsto t^m$, $a\mapsto t^ka$. The realizable multipliers $d$ on $\pi_2(M)\cong\mathbb Z$ are precisely the odd integers, and for fixed $(m,k,d)$ the pointed classes are parametrized by $\mathbb Z/|d|\mathbb Z$. We compute the composition law on this submonoid and determine the pointed self-homotopy equivalences. Composition on either side preserves the trivial-endomorphism classes, and the composite of any two such classes is the constant class. The target fundamental group acts trivially on these classes; on the other classes, its action is $(m,k,d,c)\mapsto(\varepsilon m,\varepsilon k+2r,\varepsilon d,c)$ for $r\in\mathbb Z$ and $\varepsilon=\pm1$. Thus the latter free classes have unique normal forms indexed by $\operatorname{sgn}(d)m$, $k\bmod 2$, $D=|d|$, and $c\in\mathbb Z/D\mathbb Z$. The free self-homotopy equivalence group is $C_2\times C_2$. For fixed $(m,k,d)$, all $|d|$ classes remain distinct under free homotopy despite inducing identical homomorphisms on every homotopy group.
\end{abstract}

\section{Introduction}\label{sec:introduction}

The classification of maps up to homotopy is a basic problem in algebraic topology. For a non-simply connected space, the induced homomorphisms on homotopy groups must be considered together with the action of the fundamental group and the obstructions to extending homotopies. Classical treatments of these issues include the obstruction-theoretic approach of Olum to maps into projective and lens spaces\cite{olum-cocycle} and the general framework of Whitehead\cite{whitehead-elements}. An explicit classification must determine both which algebraic data are realizable and how many homotopy classes realize each admissible choice. These two questions need not have the same answer as the problem of determining mapping degrees.

We study the closed orientable three-manifold
\[
M=\RP^3\#\RP^3.
\]
Its fundamental group is the infinite dihedral group
\[
G=\langle a,b\mid a^2=b^2=1\rangle,
\qquad t=ab,
\]
and its universal cover is $S^2\times\RR$. This manifold belongs to the $S^2\times\RR$ geometry in the classification of three-dimensional geometries; see Scott\cite{scott}. Thus $M$ has an infinite fundamental group and a noncontractible universal cover. In particular, the classification of its self-maps involves information beyond an endomorphism of $G$.

\medskip
\noindent\textbf{Related work.}
Mapping degrees and fixed point invariants provide two established approaches to self-maps of three-manifolds. Sun, Wang, Wu, and Zheng\cite{sun-wang-wu-zheng} determine self-mapping degree sets in the setting of geometrizable closed orientable three-manifolds. Their Theorem~1.1(2) records, in particular, that every integer occurs as the degree of a self-map of $\RP^3\#\RP^3$. This describes the range of the degree invariant, whereas the present problem asks for the individual pointed homotopy classes above prescribed induced homomorphisms and the identifications among them under free homotopy.

The work most directly related to our setting is the 2015 paper of Gon\c{c}alves, Wong, and Zhao\cite{gwz-s2r}, which determines Nielsen numbers of self-maps of closed three-manifolds with $S^2\times\RR$ geometry using the induced homomorphisms on $\pi_1$ and $\pi_2$. It supplies the fixed point theoretic context for the present classification problem. The same authors also determine the possible Nielsen numbers for spherical three-manifolds\cite{gwz-spherical}. These contributions, together with Jiang's treatment of Nielsen fixed point theory\cite{jiang}, relate algebraic information about a map to its fixed point behavior. Our focus here is the finer homotopy information retained after fixing the induced maps on the first two homotopy groups, together with the effect of forgetting the basepoint.

There is also a substantial mapping-space background. Hansen\cite{hansen} studies spaces of maps from closed surfaces to $S^2$, and Gon\c{c}alves and Spreafico\cite{goncalves-spreafico} compute fundamental groups of components of spaces of maps from closed surfaces to $\RP^2$. Oshima\cite{oshima} computes homotopy groups of the space of self-maps of $\RP^3$. These works concern spaces closely related to the sphere and projective-space mapping spaces that occur below. Our calculation uses the evaluation fibration and its Whitehead-product description, as developed in the classical work of Whitehead and Lang\cite{whitehead-products,jhc-whitehead,lang}. For antipodal-equivariant sphere maps, we use the equivariant homotopy classification of Rubinsztein\cite{rubinsztein}, in the formulation recalled by Jezierski and Marzantowicz\cite{jezierski}. The general equivariant framework is described in\cite{may-equivariant}.

\medskip
\noindent\textbf{The classification problem and the results.}

\begin{maintheorem}[Homotopy classification of pointed self-maps]\label{thm:classification}
Let $M=\mathbb{RP}^3\#\mathbb{RP}^3$, let $x_0\in M$, and write
\[
G=\pi_1(M,x_0)=\langle a,b\mid a^2=b^2=1\rangle,\qquad t=ab.
\]
For a homomorphism $\varphi:G\to G$ and an integer $d$, put
\[
[M,M]_*^{\varphi,d}=\left\{[f]\in[M,M]_*\;\middle|\;f_{*,1}=\varphi,\quad f_{*,2}=d\,\mathrm{id}_{\pi_2(M,x_0)}\right\},
\]
where $[M,M]_*$ is the set of pointed homotopy classes of pointed self-maps, and $f_{*,i}$ is the homomorphism induced by $f$ on the $i$th homotopy group. The classification is as follows.
\begin{enumerate}[label=\textnormal{(\arabic*)}]
\item If $\varphi$ is trivial and $d=0$, then $[M,M]_*^{\varphi,0}$ splits into four families distinguished by $\alpha(f)\in H^2(M;\mathbb Z)\cong(\mathbb Z/2\mathbb Z)^2$. After choosing a reference homotopy class in each family, there is a bijection of sets
\[
[M,M]_*^{\varphi,0}\cong\bigsqcup_{\alpha\in H^2(M;\mathbb Z)}\mathbb Z.
\]
\item For any $m,k\in\mathbb Z$ and odd integer $d$, let $\varphi:G\to G$ be the homomorphism determined by
$\varphi(t)=t^m$ and $\varphi(a)=t^ka$.
After choosing a reference homotopy class, there is a bijection of sets
\[
[M,M]_*^{\varphi,d}\cong\mathbb Z/|d|\mathbb Z.
\]
In particular, there are exactly $|d|$ pointed homotopy classes of pointed self-maps with these two induced homomorphisms.
\end{enumerate}
In all other cases, $[M,M]_*^{\varphi,d}=\varnothing$.

The displayed bijections depend on the choices of reference homotopy classes.
\end{maintheorem}

\cref{thm:classification} classifies pointed self-maps by their induced homomorphisms on the first two homotopy groups. The construction of $\alpha(f)$ is given in \cref{lem:cohomology-endpoints}, and the proof of the theorem is in \cref{sec:classification}. The integer $d$ is the multiplier on the second homotopy group, whereas ordinary mapping degree is defined by the induced map on $H_3(M;\ZZ)$.

\cref{sec:classification} also gives the composition law and its consequences. In \cref{cor:double-reflection-composition}, matrices $\mathbf A=A(m,k)$ represent the fundamental-group endomorphisms, and compatible reference maps and generators determine the coordinates used for composition in the double-reflection part. Changing the reference homotopy classes may change this formula. The group of units is $[M,M]_*^\times\cong D_\infty\times C_2$ (\cref{cor:pointed-homotopy-equivalences}). The part $\mathcal T$ inducing the trivial fundamental-group homomorphism is preserved by composition on either side with arbitrary pointed homotopy classes and satisfies $\mathcal T^2=\{0\}$, where $0$ is the constant class (\cref{prop:trivial-composition}). For each realizable pair with nontrivial fundamental-group homomorphism, all $|d|$ classes induce the same homomorphisms on every homotopy group (\cref{cor:homotopy-groups-do-not-detect}).

Explicit coordinate multiplication is given on the double-reflection part; no explicit secondary-coordinate formulas are asserted for its two mixed products with $\mathcal T$. The contributions are the explicit pointed and free classifications, the computation of their secondary ambiguity, and the pointed composition law. The sphere and equivariant homotopy results used in the argument are classical inputs.

We now pass from pointed to free homotopy classes by computing the action of the target fundamental group. For a pointed self-map $f$, write $[f]_*$ for its pointed homotopy class and $[f]$ for its free homotopy class when both are under consideration.

\begin{maintheorem}[Free homotopy classification of self-maps]\label{thm:free-main}
Retain $M,x_0,G,a,t$ from Theorem~\ref{thm:classification}. For any homomorphism $\varphi:G\to G$ and integer $d$, write
\[
[M,M]^{\varphi,d}=\left\{[f]\in[M,M]\;\middle|\;
\begin{aligned}
&[f]\text{ has a pointed representative } f_0,\\
&(f_0)_{*,1}=\varphi,\\
&(f_0)_{*,2}=d\,\mathrm{id}_{\pi_2(M,x_0)}
\end{aligned}\right\},
\]
where $[M,M]$ is the set of free homotopy classes of self-maps of $M$. The classification is as follows.
\begin{enumerate}[label=\textnormal{(\arabic*)}]
\item If $\varphi$ is trivial and $d=0$, then, after choosing a reference homotopy class in each of the four families distinguished by $\alpha\in H^2(M;\ZZ)\cong(\ZZ/2\ZZ)^2$, there is a bijection of sets
\[
[M,M]^{\varphi,0}\cong\bigsqcup_{\alpha\in(\ZZ/2\ZZ)^2}\ZZ.
\]
\item For any integers $m,k$ and odd integer $d$, let $\varphi:G\to G$ be the homomorphism determined by
\[
\varphi(t)=t^m,\qquad \varphi(a)=t^ka.
\]
After choosing a reference homotopy class, there is a bijection of sets
\[
[M,M]^{\varphi,d}\cong\ZZ/|d|\ZZ.
\]
Furthermore, suppose that $\varphi':G\to G$ is determined by $\varphi'(t)=t^{m'}$ and $\varphi'(a)=t^{k'}a$, and that $d'$ is odd. Then
\[
[M,M]^{\varphi,d}=[M,M]^{\varphi',d'}
\]
if and only if
\[
\begin{cases}
\operatorname{sgn}(d)m=\operatorname{sgn}(d')m',\\
k\equiv k'\pmod 2,\\
|d|=|d'|.
\end{cases}
\]
Otherwise, these two sets are disjoint.
\end{enumerate}
In all other cases, $[M,M]^{\varphi,d}=\varnothing$.

Write $\ZZ_{\mathrm{odd}}^{>0}=\{1,3,5,\ldots\}$. The unique pointed normal form of a free homotopy class of double-reflection type is $(A(\mu,\kappa),D,c)$, where $c\in\ZZ/D\ZZ$ is the coordinate of the pointed class from Theorem~\ref{thm:classification}(2) with the compatible choices described in \cref{def:secondary-coordinates} and
\[
\mu=\operatorname{sgn}(d)m,\qquad
\kappa\equiv k\pmod 2,\quad \kappa\in\{0,1\},\qquad D=|d|,
\]
and hence we obtain a bijection of sets
\[
[M,M]\cong\left(\bigsqcup_{\alpha\in(\ZZ/2)^2}\ZZ\right)
\sqcup\left(\bigsqcup_{\substack{\mu\in\ZZ,\ \kappa\in\{0,1\}\\D\in\ZZ_{\mathrm{odd}}^{>0}}}\ZZ/D\ZZ\right).
\]
The displayed bijection depends on the choices of reference homotopy classes and generators.
\end{maintheorem}

The proof of \cref{thm:free-main} is given in \cref{sec:free-homotopy} (\cref{thm:free-classification}). The target fundamental group acts trivially on the classes inducing the trivial endomorphism; on the double-reflection part it changes $(m,k,d,c)$ to $(\varepsilon m,\varepsilon k+2r,\varepsilon d,c)$, where $r\in\ZZ$ and $\varepsilon\in\{1,-1\}$. In particular, fixing the conjugacy class of the induced fundamental-group homomorphism and a positive odd $D=|d|$ yields $D$ free homotopy classes when $m=0$ and $2D$ otherwise (\cref{cor:free-count}). The group of self-homotopy equivalences up to free homotopy is $[M,M]^\times\cong C_2\times C_2$ (\cref{cor:free-homotopy-equivalences}). For fixed admissible $(m,k,d)$, all $|d|$ pointed classes remain distinct under free homotopy, although their pointed representatives induce identical homomorphisms on all homotopy groups (\cref{cor:free-homotopy-groups-do-not-detect}).

\medskip
\noindent\textbf{Method and significance.}
For a fixed endomorphism of $G$, we pass to twisted equivariant maps of $S^2\times\RR$. A fundamental domain then converts the problem into a path-space problem in $\Map(S^2,S^2)$, with even or odd conditions at the two endpoints. Writing $B_d=\Map_d(S^2,S^2)$, the evaluation fibration and the identity $[\iota_2,\iota_2]=\pm2\eta_2$ give
\[
\pi_1(B_d)\cong
\begin{cases}
\ZZ,&d=0,\\
\ZZ/(2|d|)\ZZ,&d\ne0.
\end{cases}
\]
The remaining task is to compute the images of the endpoint mapping spaces in this fundamental group. The even endpoint spaces have trivial image in $\pi_1(B_0)$; for odd $d$, the odd endpoint space has the unique subgroup of order two as its image. The double-coset description of path components then gives exactly the classification sets in \cref{thm:classification}. To pass to free homotopy, we use the orbit description of the forgetful map~\cite[\S4.A]{hatcher}. The action is computed using the deck transformations and the pointed composition formula; on the trivial-endomorphism part, the vanishing of the first Chern class of a pulled-back tangent bundle gives a homotopy to the antipodal postcomposition~\cite{hatcher-vb}.

This calculation makes the residual homotopy information explicit in a three-manifold with nontrivial fundamental-group action. It also separates the numerical and algebraic invariants used in degree and fixed point problems from the secondary data needed for pointed and free classification. In particular, it gives a concrete setting in which equality of all induced homotopy-group homomorphisms does not imply even free homotopy, while still allowing an explicit calculation of composition on the double-reflection part.

\Cref{sec:geometry} fixes the geometric model and conventions. \Cref{sec:lifts} develops the equivariant path-space model. \Cref{sec:mapping} collects the sphere mapping-space calculations, and \cref{sec:endpoints} determines the endpoint images. \Cref{sec:classification} proves \cref{thm:classification} and gives the pointed composition law and its consequences. \Cref{sec:free-homotopy} proves \cref{thm:free-main} and gives its consequences. The double-coset lemma and its endpoint-space applications appear in \cref{sec:double-coset-body}. Selected proofs for the geometric model and equivariant lifts are given in Appendices~\ref{app:geometry} and~\ref{app:lifts}; the complete proof of the double-coset lemma is in Appendix~\ref{sec:double-coset-appendix}.

\section{Geometric model and conventions}\label{sec:geometry}

\subsection{Notation and the classification problem}

\begin{notation}[Pointed homotopy classes]\label{stmt:original-0}
Fix a basepoint $x_0$ of $M=\RP^3\#\RP^3$ and write $G=\pi_1(M,x_0)$.
The set $[M,M]_*$ consists of continuous maps $(M,x_0)\to(M,x_0)$ modulo homotopies that fix $x_0$ at every time.
For an endomorphism $\varphi:G\to G$, set
\[
[M,M]_*^\varphi=\{[f]\in[M,M]_*\mid f_{*,1}=\varphi\},
\]
where $f_{*,n}$ denotes the induced endomorphism of $\pi_n(M,x_0)$ for $n\ge1$.
\end{notation}

\begin{notation}[Group notation]\label{stmt:original-1}
Under the identification established in \cref{sec:group-model}, let $a,b$ be the nonidentity elements of the two free factors of $G\cong C_2*C_2$, and put $t=ab$. Thus
\[
G=\langle a,b\mid a^2=b^2=1\rangle
\cong\langle t,a\mid a^2=1,\ ata=t^{-1}\rangle=D_\infty.
\]
Write $\varphi_{\mathrm{triv}}:G\to G$ for the trivial endomorphism. For $m,k\in\ZZ$, let
\[
\mathbf A=A(m,k):=\begin{pmatrix}m&k\\0&1\end{pmatrix},\qquad \varphi_{\mathbf A}(t)=t^m,\quad \varphi_{\mathbf A}(a)=t^ka,\quad \varphi_{\mathbf A}(b)=t^{k-m}a.
\]
These assignments determine an endomorphism $\varphi_{\mathbf A}:G\to G$.
Set $\ZZ_{\mathrm{odd}}=2\ZZ+1$.
Write $\mathscr M=\{A(m,k)\mid m,k\in\ZZ\}$.
\end{notation}

\begin{remark}[Conventions]\label{stmt:original-2}
All maps and homotopies are continuous. Mapping spaces and products are taken in the category of compactly generated spaces, with the compactly generated compact-open topology on mapping spaces. Write $I=[0,1]$ and let $\pi_0$ denote the set of path components. Unless stated otherwise, equivariant homotopies use the trivial group action on $I$.
\end{remark}

\subsection{A model for the space and its fundamental group}\label{sec:group-model}

\begin{proposition}[Models for the punctured summand and connected sum]\label{prop:punctured-model}
Let $N=\RP^3\setminus\operatorname{int}(D^3)$, where $D^3$ is a standard embedded closed ball. Then $N\simeq\RP^2$, $\pi_1(N)\cong C_2$, and
\[
N\cong(S^2\times I)/((x,1)\sim(-x,1)),\qquad
M\cong(S^2\times I)/\bigl((x,i)\sim(-x,i),\ i=0,1\bigr).
\]
\end{proposition}

\noindent The proof is given in Appendix~\ref{app:proof-prop:punctured-model}.

\begin{proposition}[Fundamental group]\label{prop:fundamental-group}
For the generators determined by the two punctured summands,
\[
\pi_1(M,x_0)\cong C_2*C_2=\langle a,b\mid a^2=b^2=1\rangle.
\]
\end{proposition}

\noindent The proof is given in Appendix~\ref{app:proof-prop:fundamental-group}.

\begin{lemma}[Dihedral normal form]\label{lem:dihedral-normal-form}
Set $t=ab$. Then $t$ has infinite order, $at^n=t^{-n}a$ for all $n\in\ZZ$, and every element of $G$ has a unique expression $t^n$ or $t^na$. Consequently,
\[
G\cong\langle t,a\mid a^2=1,\ ata=t^{-1}\rangle=D_\infty.
\]
\end{lemma}

\begin{proof}The normal-form theorem for free products gives each nonidentity $g\in C_2*C_2$ a unique reduced word
\(g=x_1x_2\cdots x_k, \qquad x_i\in\{a,b\},\)
with distinct adjacent letters. Set $t=ab$. Then
\(t^{-1}=(ab)^{-1}=ba,\)
and
\(ata=a(ab)a=ba=t^{-1}.\)
Thus
\begin{equation}\label{eq:at}
at=t^{-1}a.
\end{equation}
Induction gives, for every $n\in\ZZ$,
\begin{equation}\label{eq:atn}
at^n=t^{-n}a.
\end{equation}
Any word in $t,a$ can therefore be rewritten with all occurrences of $a$ on the right; pairs of $a$ cancel. Every element consequently has the form
\(t^n\quad\text{or}\quad t^na.\)
This representation is unique. Indeed, $t^n=t^m$ implies $t^{n-m}=1$; if $n-m\ne0$, this is a nonempty reduced word in $a,b$, a contradiction. Thus $n=m$. Similarly, $t^na=t^ma$ implies $n=m$. Finally, $t^n$ cannot equal $t^ma$: their reduced lengths in $a,b$ have different parities. Hence
\begin{equation}\label{eq:Dinf-normal}
G=\{t^n\mid n\in\ZZ\}\sqcup\{t^na\mid n\in\ZZ\}.
\end{equation}
This is the translation--reflection normal form for the infinite dihedral group $D_\infty$.
\end{proof}

\begin{lemma}[Elements of order at most two]\label{lem:involutions}
The solutions to $g^2=1$ in $G$ are precisely $1$ and the reflections $t^na$, $n\in\ZZ$. Each reflection has order two.
\end{lemma}

\begin{proof}If $g=t^n$ and $g^2=1$, then $t^{2n}=1$. Since $t$ has infinite order, $n=0$. Thus only the identity among translations has order at most two. On the other hand, \eqref{eq:atn} gives
\((t^na)^2=t^nat^na=t^n(at^n)a=t^nt^{-n}a^2=1.\)
Every $t^na$ therefore has order two, and
\begin{equation}\label{eq:order2-elements}
\{g\in G\mid g^2=1\}
=
\{1\}\cup\{t^na\mid n\in\ZZ\}.
\end{equation}
We call the elements $t^na$ reflections.
\end{proof}

\begin{theorem}[Decomposition of the endomorphism semigroup]\label{prop:group-endomorphisms}
Let
\[
G=\pi_1(M,x_0)=\langle t,a\mid a^2=1,\ tata=1\rangle,
\qquad t=ab.
\]
Under composition, the endomorphism semigroup has the disjoint decomposition
\[
\operatorname{End}(G)=\operatorname{End}_1\sqcup\operatorname{End}_2
\sqcup\operatorname{End}_3\sqcup\operatorname{End}_4,
\]
where
\begin{itemize}
\item $\operatorname{End}_1=\{\varphi_{\mathrm{triv}}\}$, where $\varphi_{\mathrm{triv}}(t)=\varphi_{\mathrm{triv}}(a)=1$;
\item $\operatorname{End}_2=\{\varphi_{\mathbf A}\mid\mathbf A\in\mathscr M\}$, where
$\varphi_{\mathbf A}(t)=t^m$ and $\varphi_{\mathbf A}(a)=t^ka$;
\item $\operatorname{End}_3=\{\alpha_k\mid k\in\ZZ\}$, where
$\alpha_k(t)=t^ka$ and $\alpha_k(a)=1$;
\item $\operatorname{End}_4=\{\beta_k\mid k\in\ZZ\}$, where
$\beta_k(t)=\beta_k(a)=t^ka$.
\end{itemize}
The parameters in each family are unique. \end{theorem}

\begin{proof}Consider a group endomorphism
\(\varphi:G\to G.\)
Since $a^2=b^2=1$,
\(\varphi(a)^2=\varphi(b)^2=1.\)
By \cref{lem:involutions}, each of $\varphi(a),\varphi(b)$ is either $1$ or a reflection $t^na$. Conversely, the universal property of the free product shows that any such pair determines a unique endomorphism. Write $R_j=t^ja$. Since $t=ab$, we may compute $\varphi(t)$ in each family.
\begin{enumerate}[label=$\operatorname{End}_{\arabic*}$.,leftmargin=4em]
\item If $(\varphi(a),\varphi(b))=(1,1)$, then $\varphi(t)=1$ and $\varphi=\varphi_{\mathrm{triv}}$.
\item If $(\varphi(a),\varphi(b))=(R_r,R_s)$, the dihedral normal form gives
\begin{equation}\label{eq:phi-t}
\varphi(t)=R_rR_s=t^r(at^s)a=t^rt^{-s}a^2=t^{r-s}.
\end{equation}
Taking $m=r-s$ and $k=r$ gives $\varphi=\varphi_{\mathbf A}$ for $\mathbf A=A(m,k)$. Conversely, $m,k$ determine the endpoint exponents uniquely by $r=k$ and $s=k-m$.
\item If $(\varphi(a),\varphi(b))=(1,R_k)$, then $\varphi(t)=1\cdot R_k=R_k$, so $\varphi=\alpha_k$.
\item If $(\varphi(a),\varphi(b))=(R_k,1)$, then $\varphi(t)=R_k\cdot1=R_k=\varphi(a)$, so $\varphi=\beta_k$. Conversely, $b=t^{-1}a$ gives $\beta_k(b)=R_k^{-1}R_k=1$.
\end{enumerate}
These four families exhaust the possible ordered pairs and are pairwise disjoint. Uniqueness of the reflection normal form gives uniqueness of the parameters in each family.
\end{proof}

\begin{remark}[Affine action on reflection indices]
Put $R_j=t^ja$. Since $b=R_{-1}$ and $a=R_0$, we have
\[
\varphi_{\mathbf A}(R_j)=t^{mj+k}a=R_{mj+k},\qquad
\mathbf A\binom j1=\binom{mj+k}{1}.
\]
Thus, if the reflection exponents at the two endpoints are $r,s$, then
\[
(r,s)=(k,k-m),\qquad (m,k)=(r-s,r).
\]
The mixed families satisfy
\[
(\alpha_k(a),\alpha_k(b))=(1,t^ka),\qquad
(\beta_k(a),\beta_k(b))=(t^ka,1).
\]
\end{remark}

\begin{proposition}[Matrix representation and automorphism group]\label{prop:matrix-semigroup}
The set $\operatorname{End}_2$ is a unital subsemigroup of $\operatorname{End}(G)$, and
\[
\iota:\operatorname{End}_2\longrightarrow M_2(\ZZ),\qquad
\varphi_{\mathbf A}\longmapsto\mathbf A=\begin{pmatrix}m&k\\0&1\end{pmatrix}
\]
is a semigroup embedding. Explicitly,
\[
\varphi_{\mathbf A}\circ\varphi_{A(n,\ell)}=\varphi_{A(mn,k+m\ell)},
\qquad A(m,k)A(n,\ell)=A(mn,k+m\ell).
\]
Set $T=A(1,1)$ and $R=A(-1,0)$. The subgroup of $\operatorname{GL}_2(\ZZ)$ generated by $T,R$ is
\[
\langle T,R\rangle=\{A(\epsilon,k)\mid \epsilon\in\{1,-1\},\ k\in\ZZ\}.
\]
Its inverse image under $\iota$ is precisely $\operatorname{Aut}(G)$, and
$\operatorname{Aut}(G)\cong D_\infty\cong G$.
\end{proposition}
\begin{proof}
On the generators we compute
\[
(\varphi_{\mathbf A}\circ\varphi_{A(n,\ell)})(t)=t^{mn},\qquad
(\varphi_{\mathbf A}\circ\varphi_{A(n,\ell)})(a)
=\varphi_{\mathbf A}(t^\ell a)=t^{m\ell+k}a.
\]
Thus the composite belongs to $\operatorname{End}_2$, and $\varphi_{A(1,0)}=\id_G$.
Direct matrix multiplication gives the stated formula. By uniqueness of normal forms,
$\varphi_{\mathbf A}=\varphi_{A(n,\ell)}$ if and only if $m=n$ and $k=\ell$, so $\iota$ is an embedding.

Every endomorphism in $\operatorname{End}_1$, $\operatorname{End}_3$, or $\operatorname{End}_4$
has image of order at most two and hence is not an automorphism.
The translations in the image of $\varphi_{\mathbf A}$ are exactly $t^{mj}$, $j\in\ZZ$:
every source element is $t^j$ or $t^ja$, with image $t^{mj}$ or $t^{mj+k}a$, respectively.
Surjectivity therefore requires $t$ to be in the image, which implies $m=\pm1$.
Conversely, if $m=\epsilon=\pm1$, the composition formula shows that
$\varphi_{A(\epsilon,-\epsilon k)}$ is a two-sided inverse of $\varphi_{A(\epsilon,k)}$.
Consequently $\operatorname{Aut}(G)=\{\varphi_{A(\epsilon,k)}\mid\epsilon=\pm1,\ k\in\ZZ\}$.

For every integer $k$, $T^k=A(1,k)$ and $T^kR=A(-1,k)$.
Moreover, $R^2=I$ and $RTR=T^{-1}$, so every word in $T,R$ and their inverses
reduces to $T^k$ or $T^kR$. These matrices are pairwise distinct, and $T$ has infinite order;
there are no further relations. The normal form for the infinite dihedral group now shows that
$t\mapsto T$, $a\mapsto R$ defines $G\cong\langle T,R\rangle$.
Restricting $\iota$ gives $\operatorname{Aut}(G)\cong G$.
\end{proof}

\subsection{The universal cover and its deck action}

\begin{definition}[Deck-action model]\label{def:cover-action}
On $\widetilde M=S^2\times\RR$, define
\begin{equation}\label{eq:deck-ab}
a(x,u)=(-x,-u),\qquad b(x,u)=(-x,2-u).
\end{equation}
With the convention $(gh)z=g(hz)$, these involutions define an action of $G=\langle a,b\mid a^2=b^2=1\rangle$. For $t=ab$ and $n\in\ZZ$, the action is given by
\begin{equation}\label{eq:deck-tn}
t^n(x,u)=(x,u-2n),
\end{equation}
\begin{equation}\label{eq:deck-tna}
t^na(x,u)=(-x,-u-2n).
\end{equation}
\end{definition}

\begin{lemma}[Freeness and proper discontinuity]\label{lem:proper-action}
The action of $G$ on $\widetilde M$ in \cref{def:cover-action} is free. For each $z\in\widetilde M$, there exists an open neighborhood $U$ of $z$ such that
\[
U\cap gU=\varnothing\qquad(g\in G\setminus\{1\}).
\]
In particular, the action is properly discontinuous.
\end{lemma}

\noindent The proof is given in Appendix~\ref{app:proof-lem:proper-action}.

\begin{theorem}[Universal covering and deck action]\label{thm:universal-cover}
The action in \cref{def:cover-action} has the following properties.
\begin{enumerate}[label=\textnormal{(\roman*)}]
\item The quotient $(S^2\times\RR)/G$ is homeomorphic to $M$, and the quotient map
\begin{equation}\label{eq:universal-cover-connected-sum}
p:\widetilde M=S^2\times\RR\longrightarrow M
\end{equation}
is a universal covering with deck group $G$.
\item The covering and projection onto $S^2$ induce isomorphisms $\pi_k(M)\cong\pi_k(S^2)$ for $k\ge2$. In particular, $\pi_2(M)\cong\ZZ$ and $\pi_3(M)\cong\ZZ$.
\item The homomorphism
\begin{equation}\label{eq:rho-character}
\rho:G\to\{\pm1\},\qquad \rho(t^n)=1,\qquad \rho(t^na)=-1
\end{equation}
describes the action on the sphere coordinate and on $\pi_2(\widetilde M)\cong\ZZ$: these actions are $x\mapsto\rho(g)x$ and multiplication by $\rho(g)$, respectively. Thus $\rho(a)=\rho(b)=-1$ and $\rho(t)=1$.
\end{enumerate}
\end{theorem}

\noindent The proof is given in Appendix~\ref{app:proof-thm:universal-cover}.

\section{Equivariant lifts and the path-space model}\label{sec:lifts}

\subsection{Equivariant lifts for a fixed fundamental-group endomorphism}

Fix an endomorphism
\(\varphi:G\to G.\)
We establish the precise correspondence between pointed self-maps with $f_*=\varphi$ and $\varphi$-equivariant lifts.

\subsubsection{The strict correspondence between pointed maps and equivariant lifts}

Write
\([M,M]_*^\varphi=\{[f]\in[M,M]_*\mid f_*=\varphi\}.\)
Fix $\widetilde x_0\in p^{-1}(x_0)$. Throughout, $\alpha*\beta$ means first traversing $\alpha$ and then $\beta$;
this is also the order used for multiplication of loop classes. We use the following identification of the fundamental group with the deck group:
for $g=[\alpha]\in\pi_1(M,x_0)$, let $\widetilde\alpha$ be the unique lift of the pointed loop
$\alpha:I\to M$ starting at $\widetilde x_0$. The corresponding deck transformation,
also denoted $g$, satisfies
\begin{equation}\label{eq:path-end-deck}
\widetilde\alpha(1)=g\widetilde x_0.
\end{equation}

\begin{lemma}[Equivariant lifting and pointed homotopy]
\label{lem:equivariance-pi1}
Let $p:\widetilde M=S^2\times\RR\to M$ be the universal covering of $M=\RP^3\#\RP^3$.
Fix $p(\widetilde x_0)=x_0$ and an endomorphism $\varphi:G\to G$, identifying
$G=\pi_1(M,x_0)=\Deck(p)$ as above. There is a bijection between the following sets of maps:
\[
\left\{
\begin{array}{c}
f:(M,x_0)\longrightarrow(M,x_0),\\
f_*=\varphi
\end{array}
\right\}
\longleftrightarrow
\left\{
\begin{array}{c}
F:\widetilde M\longrightarrow\widetilde M,\\
F(gz)=\varphi(g)F(z),\quad F(\widetilde x_0)=\widetilde x_0
\end{array}
\right\}.
\]
The forward correspondence takes $f$ to the unique pointed lift $F$ of $fp$; its inverse takes $F$ to the unique descended map $f$ satisfying
$pF=fp$. The relation $F(gz)=\varphi(g)F(z)$ holds for every
$g\in G$ and $z\in\widetilde M$; such a map is called $\varphi$-equivariant.

This correspondence preserves homotopies. If $F_i$ corresponds to $f_i$ for $i=0,1$, then
\[
f_0\simeq_*f_1
\quad\Longleftrightarrow\quad
F_0\simeq_{\varphi,\widetilde x_0}F_1,
\]
where the right-hand side means that there is a continuous map
$\widetilde H:\widetilde M\times I\to\widetilde M$ satisfying
\[
\begin{gathered}
\widetilde H(z,0)=F_0(z),\qquad
\widetilde H(z,1)=F_1(z),\\
\widetilde H(\widetilde x_0,s)=\widetilde x_0,\qquad
\widetilde H(gz,s)=\varphi(g)\widetilde H(z,s)
\end{gathered}
\]
for all $z\in\widetilde M$, $g\in G$, and $s\in I$. Consequently there is a natural bijection
\[
[M,M]_*^\varphi
\cong
\frac{\{F:\widetilde M\to\widetilde M\mid
F(gz)=\varphi(g)F(z),\ F(\widetilde x_0)=\widetilde x_0\}}
{\simeq_{\varphi,\widetilde x_0}}.
\]
\end{lemma}

\noindent The proof is given in Appendix~\ref{app:proof-lem:equivariance-pi1}.

\subsubsection{The action of $\varphi(g)$ on the sphere coordinate}

\begin{lemma}[Sphere-coordinate equivariance]\label{lem:sphere-coordinate}
Let $\varphi:G\to G$ and $F=(P,U):\widetilde M\to\widetilde M$ be $\varphi$-equivariant. Then $P(gz)=(\rho\varphi)(g)P(z)$ for all $g\in G$ and $z\in\widetilde M$.
\end{lemma}

\begin{proof}Write a map $F:\widetilde M\to\widetilde M$ as
\[
F(z)=(P(z),U(z)),
\qquad
P:\widetilde M\to S^2,
\quad
U:\widetilde M\to\RR.
\]
The unique normal form of $\varphi(g)$ is $t^n$ or $t^na$. By \eqref{eq:deck-tn}--\eqref{eq:deck-tna}, its action on the first coordinate is
\(y\longmapsto \rho(\varphi(g))y.\)
Define
\begin{equation}\label{eq:lambda-phi}
\lambda_\varphi:=\rho\circ\varphi:G\to\{\pm1\}.
\end{equation}
Comparing first coordinates in the equivariance relation
\(F(gz)=\varphi(g)F(z)\)
gives
\begin{equation}\label{eq:P-equivariance}
P(gz)=\lambda_\varphi(g)P(z).
\end{equation}
This is the source of the even and odd endpoint conditions below.
\end{proof}

\subsubsection{Removing the condition on the lifted basepoint}

Fix $\varphi:G\to G$ and retain the basepoints and deck-group identification from \cref{lem:equivariance-pi1}.
Mapping spaces carry the compact-open topology;
when using exponential adjunction, we take mapping spaces and products in the category
of compactly generated spaces. Parameter spaces have trivial $G$-action.

\begin{definition}[Twisted equivariant mapping spaces]\label{stmt:original-7}
Let $X=\widetilde M$ with its deck action and let $Y=\widetilde M_\varphi$ be the same underlying space with action $g\cdot_\varphi y=\varphi(g)y$. Set
\[
\begin{aligned}
\mathcal E_\varphi&=\Map_G(X,Y)=\Map_\varphi^G(\widetilde M,\widetilde M)
=\{F:\widetilde M\to\widetilde M\mid F(gz)=\varphi(g)F(z)\},\\
\mathcal F_\varphi&=\Map_{\varphi,*}^G(\widetilde M,\widetilde M)
=\{F\in\mathcal E_\varphi\mid F(\widetilde x_0)=\widetilde x_0\},\\
A&=G\widetilde x_0=p^{-1}(x_0)\subset X.
\end{aligned}
\]
The homomorphism identity for $\varphi$ gives $1\cdot_\varphi y=y$ and $(gh)\cdot_\varphi y=g\cdot_\varphi(h\cdot_\varphi y)$, so the twisted action is well defined.
\end{definition}

\begin{lemma}[Restriction and evaluation fibrations]\label{lem:section-evaluation-fibration}
The inclusion $A\hookrightarrow X$ is a $G$-cofibration. Consequently the restriction
\(r:\Map_G(X,Y)\longrightarrow\Map_G(A,Y), \qquad r(F)=F|_A,\)
is a Hurewicz fibration. Moreover, evaluation
\[
e:\Map_G(A,Y)\longrightarrow\widetilde M,\qquad
e(h)=h(\widetilde x_0)
\]
is a homeomorphism with inverse
\(y\longmapsto h_y,\qquad h_y(g\widetilde x_0)=\varphi(g)y.\)
Thus $\operatorname{ev}_{\widetilde x_0}=e\circ r$ is a Hurewicz fibration, and its fiber over
$\widetilde x_0$ is precisely $\mathcal F_\varphi$:
\begin{equation}\label{eq:equivariant-evaluation-fibration}
\mathcal F_\varphi\lhook\joinrel\longrightarrow
\mathcal E_\varphi
\xrightarrow{\operatorname{ev}_{\widetilde x_0}}\widetilde M.
\end{equation}
\end{lemma}

\noindent The proof is given in Appendix~\ref{app:proof-lem:section-evaluation-fibration}.

\begin{theorem}[Removal of the lifted basepoint condition]\label{thm:remove-lift-basepoint}
For every endomorphism $\varphi:G\to G$, the inclusion $j:\mathcal F_\varphi\hookrightarrow\mathcal E_\varphi$ induces a bijection
\begin{equation}\label{eq:based-unbased-equiv-components}
j_*:\pi_0(\mathcal F_\varphi)\xrightarrow{\cong}\pi_0(\mathcal E_\varphi).
\end{equation}
Together with \cref{lem:equivariance-pi1}, this yields
\[
[M,M]_*^\varphi\cong\pi_0(\mathcal F_\varphi)\cong\pi_0(\mathcal E_\varphi).
\]
These statements include the case in which all three sets are empty.
\end{theorem}

\noindent The proof is given in Appendix~\ref{app:proof-thm:remove-lift-basepoint}.

\begin{remark}[Scope of the component bijection]\label{stmt:original-10}
The bijection in \cref{thm:remove-lift-basepoint} identifies equivariant homotopy classes with equivariant homotopy classes relative to the lifted basepoint. It asserts neither uniqueness of a representative nor a homotopy equivalence of the two mapping spaces.
\end{remark}

\subsection{From equivariant maps to paths in a mapping space}

\subsubsection{The real coordinate contributes no discrete invariant}

Fix an endomorphism $\varphi:G\to G$. Put
\[
X=\widetilde M=S^2\times\RR,
\qquad
F=(P,U):X\longrightarrow S^2\times\RR.
\]
For $h\in G$, define $\tau_h:\RR\to\RR$ by
\[
\tau_{t^n}(v)=v-2n,
\qquad
\tau_{t^na}(v)=-v-2n
\qquad(n\in\ZZ).
\]
Set $\lambda_\varphi=\rho\circ\varphi$ and define
\[
\begin{aligned}
\mathcal E_{\lambda_\varphi}
&:=\{P:X\to S^2\mid
P(gz)=\lambda_\varphi(g)P(z)\text{ for all }g\in G,z\in X\},\\
\mathcal U_\varphi
&:=\{U:X\to\RR\mid
U(gz)=\tau_{\varphi(g)}(U(z))\text{ for all }g\in G,z\in X\}.
\end{aligned}
\]

\begin{lemma}[Contractibility of the real-coordinate space]\label{lem:R-coordinate-no-pi0}
$\mathcal U_\varphi$ is nonempty and contractible, and the coordinate projection
\[
\operatorname{pr}_P:\mathcal E_\varphi\longrightarrow\mathcal E_{\lambda_\varphi},
\qquad
F=(P,U)\longmapsto P,
\]
induces a natural bijection on path components
\[
\operatorname{pr}_{P,*}:\pi_0(\mathcal E_\varphi)
\xrightarrow{\cong}
\pi_0(\mathcal E_{\lambda_\varphi}).
\]
\end{lemma}

\begin{proof}
We first describe both coordinates of a target deck transformation. By
\eqref{eq:deck-tn}--\eqref{eq:deck-tna}, if $h=t^n$, then
\(h(y,v)=(y,v-2n),\)
and if $h=t^na$, then
\(h(y,v)=(-y,-v-2n).\)
Thus every $h\in G$ satisfies
\(h(y,v)=\bigl(\rho(h)y,\tau_h(v)\bigr),\)
where $\rho(h)\in\{\pm1\}$ and $\tau_h$ is an affine transformation of the form
\[
\tau_h(v)=\varepsilon_hv+c_h,
\qquad
\varepsilon_h\in\{\pm1\},\ c_h\in\RR.
\]

Let $F=(P,U)\in\mathcal E_\varphi$. In the equivariance relation
\[
F(gz)=\varphi(g)F(z)
\]
the left side is
\(F(gz)=\bigl(P(gz),U(gz)\bigr),\)
and the right side is
\[
\varphi(g)F(z)
=\varphi(g)\bigl(P(z),U(z)\bigr)
=\bigl(\rho(\varphi(g))P(z),
\tau_{\varphi(g)}(U(z))\bigr).
\]
Comparing coordinates yields
\[
P(gz)=\rho(\varphi(g))P(z)=\lambda_\varphi(g)P(z),
\qquad
U(gz)=\tau_{\varphi(g)}(U(z)).
\]
Hence $F$ is $\varphi$-equivariant precisely when its first coordinate lies in
$\mathcal E_{\lambda_\varphi}$ and its second in $\mathcal U_\varphi$. In particular, the projection
$\operatorname{pr}_P$ continuously maps $\mathcal E_\varphi$ into
$\mathcal E_{\lambda_\varphi}$.

We prove $\mathcal U_\varphi$ is nonempty. Since $a^2=b^2=1$, every endomorphism
$\varphi:G\to G$ satisfies
\(\varphi(a)^2=\varphi(a^2)=1, \qquad \varphi(b)^2=\varphi(b^2)=1.\)
By \cref{lem:involutions}, the only elements with square equal to one are the identity and
the reflections $t^na$. Thus each of $\varphi(a),\varphi(b)$ has one of these forms, giving four cases.

\paragraph{$\operatorname{End}_1$. $\varphi(a)=1$ and $\varphi(b)=1$.}

Choose the constant function
\(U_*(x,u)=0.\)
Since $\tau_1(v)=v$,
\(U_*(az)=0=\tau_1(0)=\tau_{\varphi(a)}(U_*(z)),\)
and
\(U_*(bz)=0=\tau_1(0)=\tau_{\varphi(b)}(U_*(z)).\)

\paragraph{$\operatorname{End}_3$. $\varphi=\alpha_k$, $\varphi(a)=1$ and $\varphi(b)=t^ka$.}

Choose the constant function
\(U_*(x,u)=-k.\)
For the generator $a$, the identity $\tau_1=\id$ gives
\(U_*(az)=-k=\tau_1(-k)=\tau_{\varphi(a)}(U_*(z)).\)
For $b$, use
$\tau_{t^ka}(v)=-v-2k$ to obtain
\[
U_*(bz)=-k
=k-2k
=-(-k)-2k
=\tau_{t^ka}(-k)
=\tau_{\varphi(b)}(U_*(z)).
\]

\paragraph{$\operatorname{End}_4$. $\varphi=\beta_k$, $\varphi(a)=t^ka$ and $\varphi(b)=1$.}

Choose the constant function
\(U_*(x,u)=-k.\)
For $a$, using
$\tau_{t^ka}(v)=-v-2k$ gives
\[
U_*(az)=-k
=k-2k
=-(-k)-2k
=\tau_{t^ka}(-k)
=\tau_{\varphi(a)}(U_*(z)).
\]
For $b$, the identity $\tau_1=\id$ gives
\(U_*(bz)=-k=\tau_1(-k)=\tau_{\varphi(b)}(U_*(z)).\)

\paragraph{$\operatorname{End}_2$. $\varphi=\varphi_{\mathbf A}$, $\varphi(a)=t^ka$ and $\varphi(b)=t^{k-m}a$.}

Choose
\begin{equation}\label{eq:Urs}
U_{\mathbf A}(x,u)=mu-k.
\end{equation}
This is continuous. For $a(x,u)=(-x,-u)$, the left side is
\[
\begin{aligned}
U_{\mathbf A}(az)
&=U_{\mathbf A}(-x,-u)\\
&=m(-u)-k\\
&=-mu-k.
\end{aligned}
\]
On the other hand, $\varphi(a)=t^ka$ acts on the real coordinate by
$v\mapsto-v-2k$, so the right side is
\[
\begin{aligned}
\tau_{\varphi(a)}(U_{\mathbf A}(z))
&=-\bigl(mu-k\bigr)-2k\\
&=-mu+k-2k\\
&=-mu-k.
\end{aligned}
\]
The expressions agree. For $b(x,u)=(-x,2-u)$, the left side is
\[
\begin{aligned}
U_{\mathbf A}(bz)
&=U_{\mathbf A}(-x,2-u)\\
&=m(2-u)-k\\
&=-mu+2m-k.
\end{aligned}
\]
On the other hand, $\varphi(b)=t^{k-m}a$ acts on the real coordinate by
$v\mapsto-v-2(k-m)$, so the right side is
\[
\begin{aligned}
\tau_{\varphi(b)}(U_{\mathbf A}(z))
&=-\bigl(mu-k\bigr)-2(k-m)\\
&=-mu+2m-k.
\end{aligned}
\]
Again the two expressions agree.

Checking $a,b$ suffices for all $g\in G$. Indeed, if equivariance holds for
$g,h\in G$, so that
\[
U(gz)=\tau_{\varphi(g)}(U(z)),
\qquad
U(hz)=\tau_{\varphi(h)}(U(z)),
\]
then the convention $(gh)z=g(hz)$ and composition of deck transformations give
\[
\begin{aligned}
U(ghz)
&=\tau_{\varphi(g)}(U(hz))\\
&=\tau_{\varphi(g)}\bigl(\tau_{\varphi(h)}(U(z))\bigr)\\
&=\tau_{\varphi(g)\varphi(h)}(U(z))\\
&=\tau_{\varphi(gh)}(U(z)).
\end{aligned}
\]
Since $a,b$ generate $G$, each function constructed above belongs to
$\mathcal U_\varphi$. This proves nonemptiness.

To prove contractibility, take $U_0,U_1\in\mathcal U_\varphi$ and
$0\le\theta\le1$. Define
\(U_\theta=(1-\theta)U_0+\theta U_1.\)
For $g\in G$, write
\[
\tau_{\varphi(g)}(v)=\varepsilon_gv+c_g,
\qquad
\varepsilon_g\in\{\pm1\},\ c_g\in\RR.
\]
Since $U_0,U_1\in\mathcal U_\varphi$,
\(U_i(gz)=\varepsilon_gU_i(z)+c_g \qquad(i=0,1).\)
Therefore
\[
\begin{aligned}
U_\theta(gz)
&=(1-\theta)U_0(gz)+\theta U_1(gz)\\
&=(1-\theta)\bigl(\varepsilon_gU_0(z)+c_g\bigr)
  +\theta\bigl(\varepsilon_gU_1(z)+c_g\bigr)\\
&=\varepsilon_g\bigl((1-\theta)U_0(z)+\theta U_1(z)\bigr)
  +\bigl((1-\theta)+\theta\bigr)c_g\\
&=\varepsilon_gU_\theta(z)+c_g\\
&=\tau_{\varphi(g)}(U_\theta(z)).
\end{aligned}
\]
Thus $U_\theta\in\mathcal U_\varphi$, so $\mathcal U_\varphi$ is a nonempty convex subspace of the function space.

More explicitly, fix one of the functions $U_*$ constructed above and define
\[
H:\mathcal U_\varphi\times[0,1]\longrightarrow\mathcal U_\varphi,
\qquad
H(U,\theta)=(1-\theta)U+\theta U_*.
\]
The calculation shows that $H(U,\theta)$ remains in $\mathcal U_\varphi$, with
\(H(U,0)=U, \qquad H(U,1)=U_*.\)
Also $H(U_*,\theta)=U_*$ for every $\theta$. The adjoint of
$H$ is jointly continuous by continuity of evaluation, scalar multiplication, and addition. The exponential law
therefore gives continuity in the mapping space. Thus $H$ contracts $\mathcal U_\varphi$ to the point $U_*$,
proving contractibility.

We finally prove that projection induces a bijection on components. If
$F_\theta=(P_\theta,U_\theta)$ is a path from $F_0$ to $F_1$ in $\mathcal E_\varphi$,
its first coordinate is a path from $P_0$ to $P_1$ in $\mathcal E_{\lambda_\varphi}$.
Thus $\operatorname{pr}_{P,*}$ is well defined.

For surjectivity, take $P\in\mathcal E_{\lambda_\varphi}$. Since
$\mathcal U_\varphi$ is nonempty, choose $U_*\in\mathcal U_\varphi$ and define
\(F=(P,U_*):X\longrightarrow S^2\times\RR.\)
For $g\in G$ and $z\in X$,
\[
\begin{aligned}
F(gz)
&=\bigl(P(gz),U_*(gz)\bigr)\\
&=\bigl(\lambda_\varphi(g)P(z),
         \tau_{\varphi(g)}(U_*(z))\bigr)\\
&=\varphi(g)\bigl(P(z),U_*(z)\bigr)\\
&=\varphi(g)F(z).
\end{aligned}
\]
Thus $F\in\mathcal E_\varphi$ and $\operatorname{pr}_P(F)=P$, proving that
$\operatorname{pr}_{P,*}$ is surjective.

For injectivity, let $F_0=(P_0,U_0)$ and $F_1=(P_1,U_1)$ belong to
$\mathcal E_\varphi$, with $P_0,P_1$ in the same component of
$\mathcal E_{\lambda_\varphi}$. Choose a continuous path
\[
[0,1]\longrightarrow\mathcal E_{\lambda_\varphi},
\qquad
\theta\longmapsto P_\theta,
\]
with endpoints $P_0,P_1$, and put
\[
U_\theta=(1-\theta)U_0+\theta U_1,
\qquad
F_\theta=(P_\theta,U_\theta).
\]
The convexity calculation implies $U_\theta\in\mathcal U_\varphi$. Since
$P_\theta\in\mathcal E_{\lambda_\varphi}$, for all $g\in G$ and $z\in X$,
\[
\begin{aligned}
F_\theta(gz)
&=\bigl(P_\theta(gz),U_\theta(gz)\bigr)\\
&=\bigl(\lambda_\varphi(g)P_\theta(z),
         \tau_{\varphi(g)}(U_\theta(z))\bigr)\\
&=\varphi(g)F_\theta(z).
\end{aligned}
\]
Hence $F_\theta\in\mathcal E_\varphi$. The path $\theta\mapsto F_\theta$ is continuous
because both coordinate paths are continuous, and its endpoints are $F_0,F_1$. Thus these maps lie in the same component,
proving injectivity.

Consequently $\operatorname{pr}_{P,*}$ is a natural bijection. Its definition uses only the projection forgetting
the real coordinate and does not depend on the auxiliary $U_*$ chosen for surjectivity.
\qedhere
\end{proof}

\subsubsection{Endpoint conditions on the fundamental domain}

\begin{definition}[Even and odd maps]\label{stmt:original-12}
Let $A:S^2\to S^2$ be the antipodal map. Define subspaces of $B=\Map(S^2,S^2)$ by
\[
E=\{q\in B\mid qA=q\},\qquad O=\{q\in B\mid qA=Aq\}.
\]
Elements of $E$ and $O$ are called even and odd maps, respectively.
\end{definition}

\begin{lemma}[Endpoint conditions]\label{lem:endpoint-conditions}
Let $P\in\mathcal E_{\lambda_\varphi}$ and $P_u=P(-,u)$. Then $u\mapsto P_u$ is a path in $B$, with
\[
P_0(-x)=\lambda_\varphi(a)P_0(x),\qquad P_1(-x)=\lambda_\varphi(b)P_1(x).
\]
Thus each endpoint is even when the corresponding sign is $+1$, and odd when it is $-1$.
\end{lemma}
\begin{proof}Restrict $P$ to $S^2\times[0,1]$ and define
\(P_u:S^2\to S^2, \qquad P_u(x)=P(x,u).\)
Then $u\mapsto P_u$ is a path in the mapping space
\(B:=\Map(S^2,S^2)\)
by exponential adjunction.

At $u=0$, we have $a(x,0)=(-x,0)$. By \eqref{eq:P-equivariance},
\begin{equation}\label{eq:endpoint-a}
P_0(-x)=\lambda_\varphi(a)P_0(x).
\end{equation}
At $u=1$, we have $b(x,1)=(-x,1)$, so
\begin{equation}\label{eq:endpoint-b}
P_1(-x)=\lambda_\varphi(b)P_1(x).
\end{equation}
By the definitions of $E$ and $O$,
$P_0\in E$ if $\lambda_\varphi(a)=+1$, and $P_0\in O$ if $\lambda_\varphi(a)=-1$;
$P_1\in E$ if $\lambda_\varphi(b)=+1$, and $P_1\in O$ if $\lambda_\varphi(b)=-1$.
\end{proof}
\subsubsection{The reverse construction: equivariant extension of paths}

\begin{definition}[Endpoint path spaces]\label{def:relative-path-space}
For $\varphi:G\to G$, put $\lambda_\varphi=\rho\varphi$ and define
\[
A_0=\begin{cases}E,&\lambda_\varphi(a)=1,\\O,&\lambda_\varphi(a)=-1,\end{cases}
\qquad
A_1=\begin{cases}E,&\lambda_\varphi(b)=1,\\O,&\lambda_\varphi(b)=-1.\end{cases}
\]
For $B=\Map(S^2,S^2)$, let
\[
\mathcal P(A_0,A_1;B)=\{\gamma\in B^I\mid\gamma(0)\in A_0,\ \gamma(1)\in A_1\}.
\]
Retain the notation $\mathcal E_{\lambda_\varphi}$ for maps $P:S^2\times\RR\to S^2$ satisfying $P(gz)=\lambda_\varphi(g)P(z)$ for all $g,z$.
\end{definition}

\begin{lemma}[Restriction and equivariant extension]\label{lem:strip-path-correspondence}
For an endomorphism $\varphi:G\to G$, let $\lambda_\varphi,A_0,A_1$ be as in \cref{def:relative-path-space}. The assignment
\[
P\longmapsto\bigl(u\longmapsto P(-,u)\bigr),\qquad 0\le u\le1,
\]
is a bijection from $\mathcal E_{\lambda_\varphi}$ to $\mathcal P(A_0,A_1;B)$; its inverse is equivariant extension from $S^2\times I$.
Under this correspondence, equivariant homotopies are exactly path homotopies with the two endpoints remaining in $A_0$ and $A_1$, respectively. In particular,
\begin{equation}\label{eq:pathspace-equivalence}
\mathcal E_{\lambda_\varphi}/\simeq_G
\xrightarrow{\cong}\mathcal P(A_0,A_1;B)/\simeq_{A_0,A_1}.
\end{equation}
\end{lemma}

\begin{proof}
For $P\in\mathcal E_{\lambda_\varphi}$, define its associated path by
\(\gamma_P(u)(x)=P(x,u), \qquad 0\le u\le1.\)
Since $S^2$ is compact Hausdorff, the exponential law gives a natural homeomorphism
\[
\Map(S^2\times I,S^2)
\cong
\Map\bigl(I,\Map(S^2,S^2)\bigr),
\]
so $\gamma_P$ is continuous. Equivariance and $a(x,0)=(-x,0)$ imply
\(\gamma_P(0)(-x)=\lambda_\varphi(a)\gamma_P(0)(x).\)
Thus $\gamma_P(0)\in A_0$; likewise $b(x,1)=(-x,1)$ implies
$\gamma_P(1)\in A_1$. This defines the restriction map
\[
R:\mathcal E_{\lambda_\varphi}\to\mathcal P(A_0,A_1;B).
\]

Conversely, given $\gamma\in\mathcal P(A_0,A_1;B)$, write
$P_u=\gamma(u):S^2\to S^2$ and first define on $D$
\begin{equation}\label{eq:strip-Q}
Q(x,u)=P_u(x)=\gamma(u)(x).
\end{equation}
Continuity of evaluation $B\times S^2\to S^2$ makes $Q$ continuous. For any
$z\in S^2\times\RR$, choose $g\in G$ and $z_0\in D$ with $z=gz_0$, and define
\begin{equation}\label{eq:P-extend}
\operatorname{Ext}(\gamma)(gz_0)
=\lambda_\varphi(g)Q(z_0).
\end{equation}

We check well-definedness. Write $z_0=(x,u)$. If $0<u<1$ and $t^nz_0\in D$, then
$0\le u-2n\le1$, forcing $n=0$. If $t^naz_0\in D$, then
$0\le-u-2n\le1$. For $n\ge0$ the coordinate is negative; for $n\le-1$,
$-u-2n\ge2-u>1$. Both are impossible. Thus an orbit meeting the interior meets it only once.

At $u=0$, the only additional representation in the fundamental strip comes from
\(a(x,0)=(-x,0).\)
The two representations give the same value in \eqref{eq:P-extend} precisely when
\(Q(-x,0)=\lambda_\varphi(a)Q(x,0),\)
which is exactly $\gamma(0)\in A_0$. At $u=1$, the only additional representation comes from
\(b(x,1)=(-x,1),\)
and the compatibility condition is
\(Q(-x,1)=\lambda_\varphi(b)Q(x,1),\)
namely $\gamma(1)\in A_1$. These exhaust all ambiguities in the interior and at the boundary. Thus
\eqref{eq:P-extend} is well defined.

For global continuity, on each closed subset $gD$ we have
\[
\operatorname{Ext}(\gamma)|_{gD}(z)
=\lambda_\varphi(g)Q(g^{-1}z),
\]
which is continuous. The family $\{gD\}_{g\in G}$ is locally finite: each real-coordinate projection
is a closed interval of length one, and every point has a neighborhood bounded in the second coordinate that meets only finitely many
sets $gD$. The local definitions agree at their boundaries. The locally finite closed gluing lemma
therefore gives a continuous map
\(\operatorname{Ext}(\gamma):S^2\times\RR\to S^2.\)

For $h\in G$, write $z=gz_0$ with $z_0\in D$. Since $\lambda_\varphi$ is a homomorphism,
\[
\operatorname{Ext}(\gamma)(hz)
=\lambda_\varphi(hg)Q(z_0)
=\lambda_\varphi(h)\operatorname{Ext}(\gamma)(z).
\]
so the extension belongs to $\mathcal E_{\lambda_\varphi}$.

The constructions are inverse on maps. If $P\in\mathcal E_{\lambda_\varphi}$ and
$z=gz_0$, then
\[
\operatorname{Ext}(R(P))(z)
=\lambda_\varphi(g)P(z_0)
=P(gz_0)=P(z).
\]
Thus $\operatorname{Ext}\circ R=\id$. If
$\gamma\in\mathcal P(A_0,A_1;B)$, first extending and then restricting to $D$ amounts to taking
$g=1$ in \eqref{eq:P-extend}, giving
\(R(\operatorname{Ext}(\gamma))(u)=\gamma(u).\)
Hence $R\circ\operatorname{Ext}=\id$.

Finally, let $P^0,P^1\in\mathcal E_{\lambda_\varphi}$ be joined by an equivariant homotopy
\[
H:(S^2\times\RR)\times I\to S^2,
\qquad
H(gz,s)=\lambda_\varphi(g)H(z,s).
\]
Define
\(\Gamma(u,s)(x)=H((x,u),s).\)
The exponential law makes $\Gamma:I\times I\to B$ continuous. Applying the relations at $u=0,1$,
namely $a(x,0)=(-x,0)$ and $b(x,1)=(-x,1)$, gives
\(\Gamma(0,s)\in A_0, \qquad \Gamma(1,s)\in A_1.\)
Thus $\Gamma$ is a relative path homotopy from $R(P^0)$ to $R(P^1)$.

Conversely, suppose $\Gamma:I\times I\to B$ satisfies
\[
\Gamma(u,0)=\gamma^0(u),
\quad
\Gamma(u,1)=\gamma^1(u),
\quad
\Gamma(0,s)\in A_0,
\quad
\Gamma(1,s)\in A_1.
\]
First define on $D\times I$
\(Q(x,u,s)=\Gamma(u,s)(x),\)
and then put
\begin{equation}\label{eq:homotopy-extend}
H(gz_0,s)=\lambda_\varphi(g)Q(z_0,s).
\end{equation}
For each $s$, the endpoint conditions ensure that \eqref{eq:homotopy-extend} is well defined at both boundaries.
The family $\{gD\times I\}_{g\in G}$ remains locally finite and closed, so the same gluing argument proves continuity of $H$.
The homomorphism property gives
\(H(hz,s)=\lambda_\varphi(h)H(z,s).\)
Thus $H$ is an equivariant homotopy from $\operatorname{Ext}(\gamma^0)$ to $\operatorname{Ext}(\gamma^1)$.
The two homotopy relations therefore correspond exactly.
\end{proof}

\begin{remark}[Endpoint convention]\label{stmt:original-15}
In $\simeq_{A_0,A_1}$ the endpoints may vary within their prescribed subspaces. Accordingly, the classification involves the images of both endpoint fundamental groups, as in \cref{lem:double-coset-path}.
\end{remark}

\section{Sphere mapping spaces and classification tools}\label{sec:mapping}

\subsection{Path components of the even and odd mapping spaces}

\begin{lemma}[Degree of an even map]\label{lem:even-degree-zero}
If $q:S^2\to S^2$ satisfies $q(-x)=q(x)$ for every $x\in S^2$, then $\deg q=0$. Equivalently,
\begin{equation}\label{eq:even-degree-zero}
E\subset B_0,\qquad B_d:=\Map_d(S^2,S^2).
\end{equation}
\end{lemma}

\begin{proof}
Let $A:S^2\to S^2$ be the antipodal map. For $q\in E$, $q\circ A=q$. Taking degrees
and using $\deg A=-1$ yields
\(\deg q=\deg(q\circ A)=\deg(q)\deg(A)=-\deg q.\)
Hence $2\deg q=0$. Since the degree is an integer, $\deg q=0$.
\end{proof}

\begin{theorem}[Components of the even mapping space]\label{thm:even-components}
Factoring even maps through $p_2:S^2\to\RP^2$ induces bijections
\begin{equation}\label{eq:E-components}
\pi_0(E)\cong[\RP^2,S^2]\cong H^2(\RP^2;\ZZ)\cong\ZZ/2\ZZ.
\end{equation}
Write $E_0$ for the component containing the constant maps and $E_1$ for the other component. Both are contained in $B_0$.
\end{theorem}

\begin{proof}
Since $q(-x)=q(x)$, every even map factors through the quotient
\[
p_2:S^2\to\RP^2
\]
uniquely as
\(q=\bar q\circ p_2, \qquad \bar q:\RP^2\to S^2.\)
Conversely, each $\bar q$ gives an even map. With the compact-open topology there is therefore a natural homeomorphism
\(E\cong\Map(\RP^2,S^2),\)
and hence
\(\pi_0(E)\cong[\RP^2,S^2].\)

Since $S^2$ is simply connected, the forgetful map
$[\RP^2,S^2]_*\to[\RP^2,S^2]$ is bijective.

The space $\RP^2$ is a two-dimensional CW complex and $S^2$ is simply connected. By the Hopf--Whitney theorem
\cite[Theorem~4.24, p.~102]{cohen}, we have 
\([\RP^2,S^2]_* \cong H^2(\RP^2;\ZZ).\)
Since $H^2(\RP^2;\ZZ)\cong\ZZ/2$, these identifications prove \eqref{eq:E-components}. \cref{lem:even-degree-zero} shows that
all maps in both components have degree zero.
\end{proof}

\begin{theorem}[Degree classification of odd maps]\label{thm:odd-map-classification}
Degree induces a bijection $\pi_0(O)\xrightarrow{\cong}2\ZZ+1$. More explicitly,
\begin{equation}\label{eq:O-decomposition}
O=\bigsqcup_{d\in\ZZ_{\mathrm{odd}}}O_d,\qquad O_d=\{q\in O\mid\deg q=d\},
\end{equation}
where each $O_d$ is nonempty and path connected. Thus two odd maps are homotopic through odd maps if and only if they have the same degree.
\end{theorem}

\begin{proof}
The relation $q(-x)=-q(x)$ is precisely equivariance for the antipodal $C_2$-actions
on source and target. Rubinsztein's classification theorem\cite{rubinsztein} states that, for a finite group
$G$ acting freely on $S^n$ with $n>1$, the forgetful map
\([S^n,S^n]^G\longrightarrow[S^n,S^n]\)
is injective; if $|G|=m$, its image consists of the classes of degree in $m\ZZ+1$. See also
Jezierski--Marzantowicz\cite[Proposition~2.1 and Theorem~2.2, pp.~209--210]{jezierski}.
Apply this with $G=C_2$ acting antipodally to obtain the assertions and the decomposition
\eqref{eq:O-decomposition}.
\end{proof}

\begin{lemma}[Nonrealizability of mixed endomorphisms]\label{lem:mixed-not-realizable}
Let $\varphi:G\to G$ send exactly one of $a,b$ to $1$ and the other to a reflection. Then $[M,M]_*^\varphi=\varnothing$.
\end{lemma}

\begin{proof}
Suppose $\varphi(a)=1$ and $\varphi(b)=t^ka$. Then
$\lambda_\varphi(a)=+1$ and $\lambda_\varphi(b)=-1$, so the path
$u\mapsto P_u$ must have $P_0\in E$ and $P_1\in O$. By \cref{lem:even-degree-zero},
$\deg P_0=0$, while \cref{thm:odd-map-classification} says that $\deg P_1$ is odd. But
$u\mapsto P_u$ is an ordinary homotopy of sphere maps, along which degree is constant, a contradiction.
The case $\varphi(a)=t^ka$, $\varphi(b)=1$ is identical.
\end{proof}

\begin{remark}[Degree obstruction]\label{stmt:original-20}
The obstruction in \cref{lem:mixed-not-realizable} is the incompatibility of the endpoint degrees: an even endpoint has degree zero, whereas an odd endpoint has odd degree.
\end{remark}

\subsection{The fundamental group of the sphere mapping space $B_d$}

\begin{definition}[Degree components and evaluation fibers]\label{stmt:original-21}
Fix the north pole $N\in S^2$ and use $x_*=N$ as the evaluation point. For $d\in\ZZ$, set
\[
B_d=\{q\in\Map(S^2,S^2)\mid\deg q=d\},\qquad
F_d=\{q\in B_d\mid q(N)=N\}.
\]
Thus $F_d=\Map_{*,d}(S^2,S^2)$.
\end{definition}

\subsubsection{Evaluation fibration and its fiber}

\begin{lemma}[Evaluation fibration]\label{lem:evaluation-Bd}
For every $d\in\ZZ$, the space $B_d$ is path connected, and evaluation at $N=x_*$ is a fibration with fiber $F_d$:
\begin{equation}\label{eq:evaluation-Bd}
F_d\longrightarrow B_d\xrightarrow{\operatorname{ev}}S^2,\qquad \operatorname{ev}(q)=q(N).
\end{equation}
\end{lemma}

\begin{proof}
The inclusion $\{x_*\}\hookrightarrow S^2$ is an NDR-pair. By the function-space restriction
fibration theorem\cite[Lemma~1.1.4, p.~5]{may-ponto}, applied after adjoining disjoint basepoints, 
\(\operatorname{ev}:\Map(S^2,S^2)\longrightarrow S^2\)
is a fibration. Restriction to the degree-$d$ component gives \eqref{eq:evaluation-Bd}, whose fiber over
$x_*$ is $F_d$. Homotopy classes of sphere self-maps are classified by degree
\cite[Corollary~4.25, p.~361]{hatcher}, so $B_d$ is path connected.
\end{proof}

\begin{lemma}[Fundamental group of the evaluation fiber]\label{lem:pi1-fiber}
For every $d\in\ZZ$ and every $q_d\in F_d$, translation to the degree-zero component, followed by change of basepoint and pointed exponential adjunction, gives an isomorphism
\begin{equation}\label{eq:pi1-fiber}
H_{q_d}:\pi_1(F_d,q_d)\xrightarrow{\cong}\pi_3(S^2,N)\cong\ZZ\{\eta_2\},
\end{equation}
where $\eta_2$ is the Hopf class. In particular, $F_d\simeq F_0$.
\end{lemma}

\begin{proof}
We must first address the basepoint of each degree component. For $d\ne0$, the basepoint of $F_d$ is a degree-$d$ map,
so one cannot directly write $\pi_1(F_d)=[S^1\wedge S^2,S^2]_*$.

Put
\(F=\Map_*(S^2,S^2).\)
Pointed exponential adjunction gives $F\simeq\Omega^2S^2$, hence
\[
\pi_0(F)\cong[S^2,S^2]_*\cong\pi_2(S^2)\cong\ZZ,
\qquad
F=\bigsqcup_{d\in\ZZ}F_d.
\]
Use the model $S^2=I^2/\partial I^2$ to identify $F$ with a double loop space.
Let $c$ be the constant pointed map. Concatenation in the first coordinate gives the continuous multiplication
\[
m(f,g)([s,t])=
\begin{cases}
f([2s,t]),&0\le s\le\tfrac12,\\
g([2s-1,t]),&\tfrac12\le s\le1.
\end{cases}
\]
Both expressions equal $x_*$ at $s=\tfrac12$, and the whole boundary maps to $x_*$,
so this defines a pointed map. Closed gluing and the exponential law imply continuity of $m:F\times F\to F$.
This multiplication is homotopy associative with homotopy unit $c$. As a loop space,
$F$ is a grouplike H-space.

First verify that component indices add under multiplication. Concatenation factors as
\[
m(f,g)=(f\vee g)\circ p,
\qquad p:S^2\longrightarrow S^2\vee S^2,
\]
where $p$ also collapses the middle line of the square, and $f\vee g$ restricts to $f,g$ on the two spheres.
Let $u,u_1,u_2$ be the fundamental classes with orientations determined by these coordinates.
In
$H_2(S^2\vee S^2;\ZZ)=\ZZ u_1\oplus\ZZ u_2$, we have
\(p_*(u)=u_1+u_2.\)
Indeed, collapsing either sphere in the wedge makes the composite cover the remaining target sphere once
with positive orientation on the corresponding half. Both coefficients are therefore one, so
\[
\begin{aligned}
m(f,g)_*(u)
&=(f\vee g)_*p_*(u)\\
&=(f\vee g)_*(u_1+u_2)\\
&=\deg(f)u+\deg(g)u\\
&=(\deg(f)+\deg(g))u.
\end{aligned}
\]
By the definition of degree,
\(\deg m(f,g)=\deg f+\deg g.\)
Equivalently, the degree isomorphism
$\pi_2(S^2)\cong\ZZ$ preserves addition in the homotopy group
\cite[Corollary~4.25, p.~361]{hatcher}.

Fix $q_d\in F_d$ and $q_{-d}\in F_{-d}$, and define continuous maps
\[
\begin{aligned}
L_{q_{-d}}:F_d&\longrightarrow F_0,&
f&\longmapsto m(q_{-d},f),\\
L_{q_d}:F_0&\longrightarrow F_d,&
h&\longmapsto m(q_d,h).
\end{aligned}
\]
For $f\in F_d$, $\deg m(q_{-d},f)=-d+d=0$; for $h\in F_0$,
$\deg m(q_d,h)=d+0=d$. Thus the stated codomains are correct.

We explicitly homotope both composites to the identity. Let
\(A:F\times F\times F\times I\longrightarrow F\)
be a continuous associativity homotopy with
\(A(a,b,f,0)=m(a,m(b,f)),\qquad A(a,b,f,1)=m(m(a,b),f).\)
In the square model, use three concatenation intervals of lengths
\[
\ell_1(v)=\tfrac12-\tfrac v4,\qquad
\ell_2(v)=\tfrac14,\qquad
\ell_3(v)=\tfrac14+\tfrac v4
\quad(0\le v\le1)
\]
placing $a,b,f$ successively on them and rescaling each interval linearly to $[0,1]$.
The lengths stay positive and sum to one; all boundary values equal $x_*$. This gives a continuous homotopy.
At $v=0,1$ the lengths are, respectively,
$(\tfrac12,\tfrac14,\tfrac14)$ and $(\tfrac14,\tfrac14,\tfrac12)$,
which give the prescribed endpoints.

Let $J:F\times I\to F$ be a left-unit homotopy satisfying
\[
J(f,0)=m(c,f),\qquad J(f,1)=f.
\]
For example, put $a(v)=(1-v)/2$ and define in the square model
\[
J(f,v)([s,t])=
\begin{cases}
x_*,&0\le s\le a(v),\\
f([(s-a(v))/(1-a(v)),t]),&a(v)\le s\le1.
\end{cases}
\]
Since $1-a(v)\ge1/2$, reparametrization is continuous. At $s=a(v)$ both expressions equal
$x_*$. Gluing and the exponential law give continuity and the required endpoints.
For fixed $f$, $J(f,v)$ joins $m(c,f)$ to $f$, so its degree remains $\deg f$.

Write
\(r_d=m(q_d,q_{-d}),\qquad r_0=m(q_{-d},q_d).\)
Both have degree zero. The pointed degree classification makes $r_d,r_0$ pointed homotopic to
$c$, so there are continuous paths
\[
\alpha_d,\alpha_0:I\longrightarrow F_0,
\qquad
\alpha_d(0)=r_d,\quad\alpha_0(0)=r_0,\quad
\alpha_d(1)=\alpha_0(1)=c.
\]
The exponential law identifies homotopies with paths in the mapping space; degree invariance ensures that these paths
remain in $F_0$.

For $f\in F_d$ and $v\in I$, define
\[
K_d(f,v)=
\begin{cases}
A(q_d,q_{-d},f,3v),&0\le v\le\tfrac13,\\
m(\alpha_d(3v-1),f),&\tfrac13\le v\le\tfrac23,\\
J(f,3v-2),&\tfrac23\le v\le1.
\end{cases}
\]
At $v=1/3$ the first two expressions are $m(r_d,f)$; at $v=2/3$ the last two are
$m(c,f)$. Continuity on the pieces and agreement at their boundaries prove continuity of $K_d$.
The first piece is a path starting at a degree-$d$ map, so its degree stays $d$;
the second has degree $0+d=d$, and the third also has degree $d$.
Thus $K_d$ takes values in $F_d$, and its endpoints are
\[
\begin{aligned}
K_d(f,0)&=m(q_d,m(q_{-d},f))=(L_{q_d}\circ L_{q_{-d}})(f),\\
K_d(f,1)&=f.
\end{aligned}
\]
This proves $L_{q_d}\circ L_{q_{-d}}\simeq\id_{F_d}$.

For the other composite, define for $h\in F_0$
\[
K_0(h,v)=
\begin{cases}
A(q_{-d},q_d,h,3v),&0\le v\le\tfrac13,\\
m(\alpha_0(3v-1),h),&\tfrac13\le v\le\tfrac23,\\
J(h,3v-2),&\tfrac23\le v\le1.
\end{cases}
\]
The formulas agree at $v=1/3$ with value $m(r_0,h)$ and at $v=2/3$ with value $m(c,h)$,
so $K_0$ is continuous. The first piece starts at degree $-d+d+0=0$,
the second has degree $0+0=0$, and the third also has degree zero. Thus it stays in $F_0$.
Furthermore,
\[
\begin{aligned}
K_0(h,0)&=m(q_{-d},m(q_d,h))=(L_{q_{-d}}\circ L_{q_d})(h),\\
K_0(h,1)&=h.
\end{aligned}
\]
so $L_{q_{-d}}\circ L_{q_d}\simeq\id_{F_0}$.
The two translations are homotopy inverses. Therefore
\begin{equation}\label{eq:Fd-F0-equivalence}
F_d\simeq F_0.
\end{equation}

The fundamental-group basepoint still needs adjustment. Since $L_{q_{-d}}(q_d)=r_0$, the equivalence initially induces
$\pi_1(F_d,q_d)\cong\pi_1(F_0,r_0)$.
Reverse the path by setting $\beta(v)=\alpha_0(1-v)$, which runs from $c$ to $r_0$.
Sending a loop $\ell$ at $r_0$ to $\beta*\ell*\overline\beta$
gives the basepoint-change isomorphism $\pi_1(F_0,r_0)\cong\pi_1(F_0,c)$.
The first isomorphism below includes this basepoint change.

We now apply smash--function adjunction only at the constant pointed map $c$:
\begin{align*}
\pi_1(F_d,q_d)
&\cong\pi_1(F_0,c)\\
&=[S^1,F_0]_*\\
&\cong[S^1,\Map_*(S^2,S^2)]_*\\
&\cong[S^1\wedge S^2,S^2]_*\\
&\cong[S^3,S^2]_*\\
&=\pi_3(S^2)\cong\ZZ\{\eta_2\}.
\end{align*}
The third line holds because loops based at $c$ and their pointed homotopies remain in $F_0$.
\end{proof}

\subsubsection{Connecting homomorphism and Whitehead square}

\begin{lemma}[Whitehead--Lang formula]\label{lem:whitehead-lang}
Let $n\ge2$, let $(X,x_0)$ be a path-connected CW complex with $x_0$ a vertex, and let
$f:(S^n,*)\to(X,x_0)$ be pointed. Here $\Map_f$ denotes the unbased path component containing $f$,
and $\Map_f^*$ its evaluation fiber over $x_0$, based at $f$. For the evaluation fibration
\[
\Map_f^*(S^n,X)\longrightarrow\Map_f(S^n,X)
\xrightarrow{\operatorname{ev}}X
\]
and $i\ge2$, the connecting homomorphism is
\[
\partial_f:\pi_i(X,x_0)\longrightarrow
\pi_{i-1}(\Map_f^*(S^n,X),f).
\]
Under the standard identification
\[
H_f:\pi_{i-1}(\Map_f^*(S^n,X),f)
\xrightarrow{\cong}\pi_{i+n-1}(X,x_0)
\]
we have
\begin{equation}\label{eq:connecting-whitehead-general}
H_f(\partial_f(\alpha))=\pm[\alpha,f].
\end{equation}
\end{lemma}

\begin{proof}
This is the Whitehead--Lang formula for evaluation fibrations. The original Whitehead-product
formulation is in G.~W.~Whitehead\cite[Theorem~(3.2)]{whitehead-products}, with the sign correction in
J.~H.~C.~Whitehead\cite{jhc-whitehead}; see in particular Lang\cite[Theorem~2.6, p.~205]{lang} and also Hansen\cite{hansen}.
Lang's Lemma~2.1 translates the fiber component at $f$ to the constant component,
and Lemma~2.2 applies pointed adjunction there. Taking the sphere special case of
Theorem~2.6 then identifies the connecting map, after these identifications, with the
Whitehead product with $[f]$. This is exactly the displayed formula (up to the sign convention).
Both sides lie in $\pi_{i+n-1}(X)$ since $(i-1)+n=i+n-1$. We use only the image and
cokernel, so the overall sign does not affect the conclusions.
\end{proof}

\begin{remark}[Sphere specialization]
Apply \cref{lem:whitehead-lang} with $n=2$, $X=S^2$, $i=2$, and
$[q_d]=d\iota_2$ to obtain
\begin{equation}\label{eq:connecting-whitehead}
H_{q_d}(\partial_d(\iota_2))
=\pm[\iota_2,d\iota_2]
=\pm d[\iota_2,\iota_2].
\end{equation}
\end{remark}

\begin{theorem}[Whitehead square and connecting homomorphism]\label{thm:partial-d}
Let $\iota_2\in\pi_2(S^2)$ be the identity class and $\eta_2\in\pi_3(S^2)$ a Hopf generator. Then
\begin{equation}\label{eq:whitehead-square}
[\iota_2,\iota_2]=\pm2\eta_2.
\end{equation}
For the evaluation fibration of \cref{lem:evaluation-Bd}, the connecting homomorphism, under \cref{lem:pi1-fiber}, is
\begin{equation}\label{eq:partial-d}
\partial_d:\ZZ\longrightarrow\ZZ,\qquad n\longmapsto\pm2dn.
\end{equation}
The sign is determined by the conventions for generators and connecting homomorphisms.
\end{theorem}

\begin{proof}
G.~W.~Whitehead gives a direct calculation in \emph{Elements of Homotopy Theory}, Chapter~XII, \S2.
Corollary~(2.5) relates the $P$-homomorphism to the Whitehead product;
the proof of Corollary~(2.7) uses
\(P(\iota_5)=[\iota_2,\iota_2]\)
and states that the Hopf-invariant isomorphism
\(H_0:\pi_3(S^2)\xrightarrow{\cong}\ZZ\)
satisfies $H_0(P(\iota_5))=2$. Thus $P(\iota_5)$ is twice a generator. Orientation and
Whitehead-product conventions affect only the sign, yielding \eqref{eq:whitehead-square}
\cite[Chapter~XII, \S2, Corollaries~(2.5) and~(2.7)]{whitehead-elements}.

By \eqref{eq:connecting-whitehead} and bilinearity of the Whitehead product,
\[
H_{q_d}(\partial_d(\iota_2))
=\pm d[\iota_2,\iota_2]
=\pm2d\eta_2.
\]
Choose the generator of $\pi_1(F_d)$ mapping under $H_{q_d}$ to $\eta_2$. This gives
\eqref{eq:partial-d}.
\end{proof}

\begin{theorem}[Fundamental groups of degree components]\label{thm:pi1-Bd}
For any $q_d\in F_d$, the fiber inclusion induces a surjection on fundamental groups with kernel $2d\ZZ$ under \cref{lem:pi1-fiber}. Consequently,
\begin{equation}\label{eq:pi1-B0}
\pi_1(B_0,q_0)\cong\ZZ,
\end{equation}
and, for $d\ne0$,
\begin{equation}\label{eq:pi1-Bd}
\pi_1(B_d,q_d)\cong\ZZ/(2|d|)\ZZ.
\end{equation}
These descriptions hold at any basepoint in $B_d$ by change of basepoint.
\end{theorem}

\begin{proof}
The long exact sequence of \eqref{eq:evaluation-Bd} contains
\begin{equation}\label{eq:les-evaluation}
\pi_2(S^2)
\xrightarrow{\partial_d}
\pi_1(F_d)
\longrightarrow
\pi_1(B_d)
\longrightarrow
\pi_1(S^2)=0.
\end{equation}
Hence
\(\pi_1(B_d)\cong\coker(\partial_d).\)
By \cref{thm:partial-d}, $\partial_d$ is multiplication by $\pm2d$. For $d=0$ its image is zero,
giving $\pi_1(B_0)\cong\ZZ$; for $d\ne0$ its image is $2|d|\ZZ$, giving
$\pi_1(B_d)\cong\ZZ/(2|d|)\ZZ$.
\end{proof}

\begin{remark}[Basepoints in nonzero degree]\label{stmt:original-27}
When $d\ne0$, exponential adjunction applies after translation from $F_d$ to $F_0$ and transport to the constant map, as in \cref{lem:pi1-fiber}. A loop based at a degree-$d$ map does not directly adjoint to a map $S^1\wedge S^2\to S^2$.
\end{remark}

\subsection{Double-coset classification of paths}\label{sec:double-coset-body}

\begin{lemma}[Double-coset classification of paths]\label{lem:double-coset-path}
Let $i:A\to B$ and $j:C\to B$ be continuous maps between nonempty path-connected spaces, and give
\[
\mathcal H(A,C;B)=\{(a,\gamma,c)\in A\times B^I\times C\mid
\gamma(0)=i(a),\ \gamma(1)=j(c)\}
\]
the subspace topology. Fix $a_0\in A$, $c_0\in C$ and a path $\omega:i(a_0)\rightsquigarrow j(c_0)$.
Put $\Gamma=\pi_1(B,i(a_0))$, $H=i_*\pi_1(A,a_0)$, and
\[
K=\tau_\omega\bigl(j_*\pi_1(C,c_0)\bigr),\qquad
\tau_\omega([\ell])=[\omega*\ell*\bar\omega].
\]
Then the reference path determines a bijection
\[
\pi_0\mathcal H(A,C;B)\cong H\backslash\Gamma/K.
\]
If $\Gamma$ is abelian, this is a bijection with $\Gamma/(H+K)$, and that quotient group acts freely and transitively on the component set.
\end{lemma}

\noindent The proof is given in Appendix~\ref{app:proof-lem:double-coset-path}.

\begin{remark}[Application to endpoint spaces]
We apply this to the following two pairs of endpoint spaces, with inclusions $i,j$.
The subgroups on the right are transported to a common basepoint along the chosen reference path.

\noindent\textbf{(1) Even endpoints.}
Fix $\varepsilon_0,\varepsilon_1\in\ZZ/2$ and take
\(A=E_{\varepsilon_0},\qquad C=E_{\varepsilon_1},\qquad B=B_0.\)
Here $E_\varepsilon$ is a path component of the even mapping space.
Since $\pi_1(B_0)\cong\ZZ$, the lemma gives
\[
\pi_0\mathcal H(E_{\varepsilon_0},E_{\varepsilon_1};B_0)
\cong
\frac{\pi_1(B_0)}
{\im\pi_1(E_{\varepsilon_0})+\im\pi_1(E_{\varepsilon_1})}.
\]

\noindent\textbf{(2) Odd endpoints.}
Fix an odd integer $d$ and take
\(A=O_d,\qquad C=O_d,\qquad B=B_d.\)
Here $O_d$ is the path-connected space of odd maps of degree $d$.
Since $\pi_1(B_d)\cong\ZZ/(2|d|)\ZZ$, we obtain
\[
\pi_0\mathcal H(O_d,O_d;B_d)
\cong
\frac{\pi_1(B_d)}
{\im\pi_1(O_d)+\im\pi_1(O_d)}.
\]
The images of these endpoint spaces in the fundamental groups of the corresponding mapping spaces are computed in \cref{sec:endpoints}.
\end{remark}

\section{Images of the endpoint mapping spaces}\label{sec:endpoints}

\subsection{The image of the even endpoint space in $\pi_1(B_0)$}

\begin{definition}[Antipodal precomposition]
Let $A(x)=-x$ and define $C_A:B_0\to B_0$ by $C_A(q)=q\circ A$. By the definition of $E$,
\begin{equation}\label{eq:CA-even-id}
C_A|_E=\id_E.
\end{equation}
\end{definition}

\begin{lemma}[Antipodal precomposition]\label{lem:CA-minus-one}
For every $q\in E$, the self-map $C_A:B_0\to B_0$ fixes $q$ and induces inversion on $\pi_1(B_0,q)$. In additive notation,
\begin{equation}\label{eq:CA-pi1-minus}
(C_A)_*=-\id:\pi_1(B_0,q)\longrightarrow\pi_1(B_0,q)\cong\ZZ.
\end{equation}
This includes a constant map and a representative $q_\varepsilon\in E_\varepsilon$ for each $\varepsilon\in\{0,1\}$.
\end{lemma}

\begin{proof}
Fix the north pole $N\in S^2$ and put $c(x)=N$. Choose a degree-$(-1)$ homeomorphism fixing $N$,
\(J:S^2\to S^2, \qquad J(N)=N, \qquad\deg J=-1.\)
for example reflection in a plane through the polar axis. Since $\deg J=\deg A=-1$ and degree classifies $[S^2,S^2]$,
there is a homotopy $H_t:S^2\to S^2$ from $A$ to $J$. Define
\(C_t:B_0\to B_0, \qquad C_t(q)=q\circ H_t.\)
Every $H_t$ has degree $-1$, so $C_t$ preserves $B_0$; also $C_t(c)=c$. Thus
$C_t$ is a pointed homotopy from $C_A$ to $C_J$, giving
\begin{equation}\label{eq:CA-CJ-same}
(C_A)_{*,c}=(C_J)_{*,c}.
\end{equation}

Consider evaluation at $N$:
\[
F_0=\Map_{*,0}(S^2,S^2)
\xrightarrow{i}B_0
\xrightarrow{\operatorname{ev}_N}S^2.
\]
Since $J(N)=N$,
\(\operatorname{ev}_N(C_J(q))=q(J(N))=q(N),\)
and $C_J(F_0)\subset F_0$. Hence $C_J$ gives a strict map of fibrations covering $\id_{S^2}$:
\[
\begin{tikzcd}[column sep=4.5em,row sep=1.6em]
F_0 \arrow[r,"C_J|_{F_0}"] \arrow[d,hook,"i"'] & F_0 \arrow[d,hook,"i"]\\
B_0 \arrow[r,"C_J"] \arrow[d,"\operatorname{ev}_N"'] & B_0 \arrow[d,"\operatorname{ev}_N"]\\
S^2 \arrow[r,"\id_{S^2}"'] & S^2.
\end{tikzcd}
\]

To compute on the fiber, put $F=\Map_*(S^2,S^2)$. Inclusion $F_0\hookrightarrow F$ induces
an isomorphism on $\pi_1$: any loop at $c$ has path-connected image containing $c$, so it lies in
$F_0$, and the same holds for pointed homotopies between such loops. Adjunction therefore gives
\[
H_c:\pi_1(F_0,c)\xrightarrow{\cong}\pi_3(S^2).
\]
The adjoint of a loop $\gamma:S^1\to F_0$ is
\[
\widehat\gamma:S^1\wedge S^2\to S^2,
\qquad
\widehat\gamma([s,x])=\gamma(s)(x).
\]
The adjoint of $C_J\circ\gamma$ is
\(\widehat\gamma\circ(1_{S^1}\wedge J).\)
The reduced homology K\"unneth isomorphism for the smash product gives
\[
\widetilde H_3(S^1\wedge S^2)
\cong\widetilde H_1(S^1)\otimes\widetilde H_2(S^2).
\]
The actions of $1_{S^1}$ and $J$ are $+1$ and $-1$, respectively, so
\(\deg(1_{S^1}\wedge J)=-1.\)
Precomposition by this degree-$(-1)$ map multiplies $\pi_3(S^2)$ by $-1$. Hence
\((C_J|_{F_0})_*=-1 \quad\text{on }\pi_1(F_0,c).\)

\cref{thm:partial-d} gives $\partial_0=0$; also
$\pi_1(S^2)=0$. The evaluation long exact sequence therefore gives an isomorphism
\[
i_*:\pi_1(F_0,c)\xrightarrow{\cong}\pi_1(B_0,c).
\]
Naturality of the strict fibration diagram gives
\[
\begin{tikzcd}[column sep=5em,row sep=3em]
\pi_1(F_0,c) \arrow[r,"(C_J|_{F_0})_*"] \arrow[d,"i_*"',"\cong"]
 & \pi_1(F_0,c) \arrow[d,"i_*","\cong"']\\
\pi_1(B_0,c) \arrow[r,"(C_J)_*"'] & \pi_1(B_0,c).
\end{tikzcd}
\]
Equivalently,
\((C_J)_*\circ i_*=i_*\circ(C_J|_{F_0})_*.\)
Thus $(C_J)_{*,c}=-1$, and \eqref{eq:CA-CJ-same} gives
$(C_A)_{*,c}=-1$.

Now move the basepoint to $q_\varepsilon\in E_\varepsilon$. Choose a path
\(\epsilon:c\rightsquigarrow q_\varepsilon \quad\text{in }B_0\)
and denote basepoint change by $\beta_\epsilon$. Since $C_A(c)=c$ and
$C_A(q_\varepsilon)=q_\varepsilon$, the path $C_A\epsilon$ also runs from $c$ to $q_\varepsilon$.
Naturality gives
\[
\begin{tikzcd}[column sep=5em,row sep=3em]
\pi_1(B_0,c) \arrow[r,"(C_A)_{*,c}"] \arrow[d,"\beta_\epsilon"',"\cong"]
 & \pi_1(B_0,c) \arrow[d,"\beta_{C_A\epsilon}","\cong"']\\
\pi_1(B_0,q_\varepsilon) \arrow[r,"(C_A)_{*,q_\varepsilon}"']
 & \pi_1(B_0,q_\varepsilon).
\end{tikzcd}
\]
Equivalently,
\[
(C_A)_{*,q_\varepsilon}\circ\beta_\epsilon
=\beta_{C_A\epsilon}\circ(C_A)_{*,c}.
\]
Basepoint changes along two paths with the same endpoints differ by an inner automorphism of
$\pi_1(B_0,q_\varepsilon)$, which is isomorphic to $\ZZ$ by \cref{thm:pi1-Bd}.
It is abelian, so inner automorphisms are trivial. Thus $\beta_{C_A\epsilon}=\beta_\epsilon$, giving
$(C_A)_{*,q_\varepsilon}=-1$.
\end{proof}

\begin{theorem}[Vanishing of the even endpoint image]\label{thm:E-image-zero}
For $\varepsilon\in\{0,1\}$ and $q_\varepsilon\in E_\varepsilon$, inclusion induces the zero homomorphism
\begin{equation}\label{eq:E-image-zero}
\pi_1(E_\varepsilon,q_\varepsilon)\longrightarrow\pi_1(B_0,q_\varepsilon).
\end{equation}
\end{theorem}

\begin{proof}
Fix $q_\varepsilon\in E_\varepsilon$ and let
$i_\varepsilon:E_\varepsilon\hookrightarrow B_0$. Take
$[\gamma]\in\pi_1(E_\varepsilon,q_\varepsilon)$. Every $\gamma(t)$ is even, so pointwise
\begin{equation}\label{eq:CA-gamma-strict}
(C_A\circ\gamma)(t)=\gamma(t)\circ A=\gamma(t).
\end{equation}
Put
\[
n=(i_\varepsilon)_*[\gamma]
\in\pi_1(B_0,q_\varepsilon)\cong\ZZ.
\]
Since $C_A\circ i_\varepsilon=i_\varepsilon$, we have $(C_A)_*n=n$; \cref{lem:CA-minus-one}
also gives $(C_A)_*n=-n$. Thus $2n=0$. Since $\ZZ$ is torsion-free, $n=0$. As
$[\gamma]$ was arbitrary, the image is zero.
\end{proof}

\begin{remark}[Evaluation and the antipodal map]\label{stmt:original-32}
Since $A(N)\ne N$, precomposition by $A$ need not preserve the evaluation fiber $F_0$. The calculation uses a homotopic reflection $J$ fixing $N$, for which the evaluation-fibration diagram commutes strictly.
\end{remark}

\subsection{The image of the odd endpoint space in $\pi_1(B_d)$}

\begin{definition}[Antipodal conjugation]\label{def:theta-involution}
For odd $d$, define $\Theta:B_d\to B_d$ by $\Theta(q)=AqA$. The identities $A^2=\id$ and $\deg A=-1$ imply $\Theta^2=\id$ and $\deg(AqA)=d$. Moreover,
\begin{equation}\label{eq:Od-fix-theta}
O_d=\operatorname{Fix}(\Theta),
\end{equation}
since $qA=Aq$ is equivalent to $AqA=q$.
\end{definition}

\begin{lemma}[Odd representatives fixing the evaluation point]\label{lem:odd-based-representative}
For every odd integer $d$, the intersection $O_d\cap F_d$ is nonempty. Equivalently, there is a degree-$d$ odd map $q_d:S^2\to S^2$ satisfying $q_d(N)=N$.
\end{lemma}

\begin{proof}
By \cref{thm:odd-map-classification}, there is an odd degree-$d$ map
$q:S^2\to S^2$. Write $q(N)=y$. Transitivity of $SO(3)$ on the sphere allows us to choose
$R\in SO(3)$ with $R(y)=N$. Set $q_d=R\circ q$. Then
$q_d(N)=N$ and $\deg q_d=d$. Rotations satisfy $R(-z)=-R(z)$, so
\(q_d(-x)=R(q(-x))=R(-q(x))=-R(q(x))=-q_d(x).\)
Thus $q_d\in O_d\cap F_d$.
\end{proof}

\subsubsection{Precomposition, postcomposition, and naturality for a reflection}

\begin{definition}[Reflection operators]\label{stmt:original-35}
For $N=(0,0,1)\in S^2\subset\RR^3$, fix the reflection
\[
J(x_1,x_2,x_3)=(-x_1,x_2,x_3).
\]
Then $\deg J=-1$, $J(N)=N$ and $JA=AJ$. For $d\in\ZZ$, define
\[
R_J:B_d\to B_{-d},\quad q\mapsto qJ,\qquad
L_J:B_{-d}\to B_d,\quad q\mapsto Jq,
\]
and $\Theta_J=L_JR_J$, so that $\Theta_J(q)=JqJ$.
\end{definition}

\begin{lemma}[Naturality under precomposition]\label{lem:precomp-natural}
For each $d\in\ZZ$, define
\(R_J:B_d\longrightarrow B_{-d},\qquad R_J(q)=q\circ J,\)
and its restriction to the fibers
\(r_J:F_d\longrightarrow F_{-d},\qquad r_J(q)=q\circ J.\)
Then $R_J$ covers $\id_{S^2}$ and gives a strict map of evaluation fibrations
\[
\begin{tikzcd}[column sep=4.8em,row sep=1.7em]
F_d \arrow[r,"r_J"] \arrow[d,hook] & F_{-d} \arrow[d,hook] \\
B_d \arrow[r,"R_J"] \arrow[d,"\operatorname{ev}"'] & B_{-d} \arrow[d,"\operatorname{ev}"] \\
S^2 \arrow[r,"\id_{S^2}"'] & S^2 .
\end{tikzcd}
\]

Choose $q_d\in F_d$; its image $q_dJ=q_d\circ J$ belongs to $F_{-d}$.
Use the compatible identifications from \cref{lem:pi1-fiber}, obtained by component translation and pointed exponential adjunction,
\[
H_{q_d}:\pi_1(F_d,q_d)\xrightarrow{\cong}\pi_3(S^2),
\qquad
H_{q_dJ}:\pi_1(F_{-d},q_dJ)\xrightarrow{\cong}\pi_3(S^2).
\]
Under these identifications,
\[
H_{q_dJ}\circ(r_J)_*
=(1_{S^1}\wedge J)^*\circ H_{q_d}
=-H_{q_d}.
\]
In particular,
\[
\boxed{(r_J)_*=-1}
\]
where the equality is understood under the displayed identifications. The notation $(1_{S^1}\wedge J)^*$ denotes precomposition on the domain.
\end{lemma}

\begin{proof}
Since $\deg J=-1$, for $q\in B_d$,
\(\deg(q\circ J)=\deg(q)\deg(J)=-d,\)
so $R_J$ takes values in $B_{-d}$. As $J(N)=N$,
\[
\operatorname{ev}(R_J(q))
=(q\circ J)(N)
=q(N)
=\operatorname{ev}(q).
\]
Consequently,
\(\operatorname{ev}\circ R_J =\id_{S^2}\circ\operatorname{ev},\)
so $R_J$ covers $\id_{S^2}$.

If $q\in F_d$, then $q(N)=N$, giving
\((r_J(q))(N)=(q\circ J)(N)=q(N)=N.\)
Thus $r_J(q)\in F_{-d}$ and $r_J$ is the restriction of $R_J$ to the fiber. This proves strict commutativity.

We next compute the action on the fiber fundamental group, starting with $F_0$ based at the constant map $c$. If
$\gamma:S^1\to F_0$ is a loop at $c$, its adjoint is
\[
\widehat\gamma:S^1\wedge S^2\longrightarrow S^2,
\qquad
\widehat\gamma([t,x])=\gamma(t)(x).
\]
For $r_J\circ\gamma$,
\[
\widehat{r_J\gamma}([t,x])
=(\gamma(t)\circ J)(x)
=\widehat\gamma([t,J(x)]).
\]
Hence
\(\widehat{r_J\gamma} =\widehat\gamma\circ(1_{S^1}\wedge J).\)
Thus under the identification
\(H_c:\pi_1(F_0,c)\xrightarrow{\cong}\pi_3(S^2)\)
the homomorphism $(r_J)_*$ corresponds to precomposition by
$1_{S^1}\wedge J$.

The reduced homology isomorphism for smash products is
\[
\widetilde H_3(S^1\wedge S^2;\ZZ)
\cong
\widetilde H_1(S^1;\ZZ)\otimes
\widetilde H_2(S^2;\ZZ).
\]
Here $1_{S^1}$ acts by $+1$ on $\widetilde H_1(S^1)$, while $J$ acts on
$\widetilde H_2(S^2)$ by $\deg J=-1$. Therefore
\(\deg(1_{S^1}\wedge J)=-1.\)
For every $[\beta]\in\pi_3(S^2)$,
\([\beta\circ(1_{S^1}\wedge J)]=-[\beta].\)
and hence
\(H_c\circ(r_J)_*=-H_c \quad\text{on }\pi_1(F_0,c).\)

We now handle the basepoint in an arbitrary degree component. Choose a square model compatible with $J$,
$(S^2,N)\cong(I^2/\partial I^2,*)$, in which
\[
J([s,t])=[s,1-t].
\]
Such a model follows from stereographic projection from $N$: $J$ becomes a linear reflection of the plane.
Use coordinatewise homeomorphisms $(0,1)^2\cong\RR^2$ and collapse the entire boundary to the point at infinity.
Concatenate in the first coordinate, using the multiplication $m$ of \cref{lem:pi1-fiber}.
Since $J$ reverses only the second coordinate, calculation on the two pieces gives
\[
\begin{aligned}
r_Jm(f,g)([s,t])
&=m(f,g)([s,1-t])\\
&=\begin{cases}
f([2s,1-t]),&0\le s\le\frac12,\\
g([2s-1,1-t]),&\frac12\le s\le1,
\end{cases}\\
&=m(r_Jf,r_Jg)([s,t]).
\end{aligned}
\]
Thus $r_J$ strictly preserves this concatenation multiplication.

To avoid confusing basepoints when $d=0$, write $q=q_d$ and $q'=q\circ J$;
these may be different maps even when they have the same degree.
Define translations
\[
T:F_d\to F_0,\quad T(f)=m(q',f),
\qquad
T':F_{-d}\to F_0,\quad T'(g)=m(q,g).
\]
Degree additivity puts both maps in $F_0$, and \cref{lem:pi1-fiber}
makes them homotopy equivalences. Since $J^2=\id$, $r_Jq'=q$, so
\(T'r_J(f)=m(q,r_Jf)=r_Jm(q',f)=r_JT(f).\)

Put $v=T(q)=m(q',q)$. Its degree is zero, so choose a path in $F_0$
$\rho:v\rightsquigarrow c$. Pointwise precomposition gives
\(\rho'=r_J\rho:T'(q')=m(q,q')\rightsquigarrow c.\)
Let $\beta_\rho$ be the fundamental-group isomorphism changing the basepoint from $v$ to $c$ along $\rho$,
and define $\beta_{\rho'}$ similarly. Explicitly, $\beta_\rho$ takes a loop $\ell$ at $v$
to the class of $\bar\rho*\ell*\rho$ based at $c$.
Naturality of concatenation and $T'r_J=r_JT$ give
\[
\beta_{\rho'}\circ T'_*\circ(r_J)_*
=(r_J|_{F_0})_*\circ\beta_\rho\circ T_*.
\]
Choose compatibly
\[
H_q=H_c\circ\beta_\rho\circ T_*,\qquad
H_{q'}=H_c\circ\beta_{\rho'}\circ T'_*.
\]
The computation in the constant component then gives
\[
\begin{aligned}
H_{q'}\circ(r_J)_*
&=H_c\circ(r_J|_{F_0})_*\circ\beta_\rho\circ T_*\\
&=-H_c\circ\beta_\rho\circ T_*=-H_q.
\end{aligned}
\]
Changing the basepoint path changes transport only by an inner automorphism; since
$\pi_1(F_0,c)\cong\ZZ$ is abelian, the result is unaffected.
Changing a translation factor within its degree component gives a homotopy of translations by applying $m$ to a connecting path.
Using the associated basepoint path gives the same identification. Thus these choices are compatible with
\cref{lem:pi1-fiber}, and $(r_J)_*$ is multiplication by $-1$ under these identifications.
\end{proof}

\begin{lemma}[Naturality under postcomposition]\label{lem:postcomp-natural}
For each $d\in\ZZ$, define
\(L_J:B_{-d}\longrightarrow B_d, \qquad L_J(q)=J\circ q,\)
and its restriction to the fibers
\(\ell_J:F_{-d}\longrightarrow F_d, \qquad \ell_J(q)=J\circ q.\)
Then $L_J$ covers $J:S^2\to S^2$ and gives a strict map of evaluation fibrations
\[
\begin{tikzcd}[column sep=4.8em,row sep=1.7em]
F_{-d} \arrow[r,"\ell_J"] \arrow[d,hook] & F_d \arrow[d,hook] \\
B_{-d} \arrow[r,"L_J"] \arrow[d,"\operatorname{ev}"'] & B_d \arrow[d,"\operatorname{ev}"] \\
S^2 \arrow[r,"J"'] & S^2 .
\end{tikzcd}
\]

Choose $q_d\in F_d$ and put $q_{-d}=q_dJ\in F_{-d}$. Under the identifications of \cref{lem:pi1-fiber},
\[
H_{q_dJ}:\pi_1(F_{-d},q_dJ)\xrightarrow{\cong}\pi_3(S^2),
\qquad
H_{Jq_dJ}:\pi_1(F_d,Jq_dJ)\xrightarrow{\cong}\pi_3(S^2),
\]
we have
\(H_{Jq_dJ}\circ(\ell_J)_*=J_*\circ H_{q_dJ}.\)
Moreover,
\[
\boxed{J_*:\pi_3(S^2)\longrightarrow\pi_3(S^2)\text{ is the identity},}
\]
and therefore
\[
\boxed{(\ell_J)_*=+1}
\]
where equality is understood under the two displayed identifications.
\end{lemma}

\begin{proof}
For $q\in B_{-d}$,
\(\deg(J\circ q)=\deg(J)\deg(q)=(-1)(-d)=d,\)
so $L_J$ takes values in $B_d$. Furthermore,
\[
\operatorname{ev}(L_J(q))
=(J\circ q)(N)
=J(q(N))
=J(\operatorname{ev}(q)).
\]
Hence
\(\operatorname{ev}\circ L_J =J\circ\operatorname{ev},\)
so $L_J$ covers $J$. If $q\in F_{-d}$, then $q(N)=N$ and $J(N)=N$, giving
\((\ell_J(q))(N)=J(q(N))=J(N)=N.\)
Thus $\ell_J(q)\in F_d$ and the diagram commutes strictly.

First take a loop $\gamma$ at $c$ in $F_0$. Its adjoint after pointwise postcomposition satisfies
\[
\widehat{\ell_J\gamma}(t,x)
=J(\gamma(t)(x))
=J\circ\widehat\gamma(t,x).
\]
Thus under smash--function adjunction, postcomposition corresponds to the target-induced map
\[
J_*:\pi_3(S^2)\longrightarrow\pi_3(S^2).
\]
Although $\deg J=-1$, the induced map on $\pi_3(S^2)$ is not multiplication by $-1$. Indeed, let
$\iota_2\in\pi_2(S^2)$ be the identity class. Since $\deg J=-1$,
\(J_*\iota_2=-\iota_2.\)
Naturality of the Whitehead product gives
\[
J_*[\iota_2,\iota_2]
=[J_*\iota_2,J_*\iota_2]
=[-\iota_2,-\iota_2]
=[\iota_2,\iota_2].
\]
The standard Whitehead-square calculation is
\([\iota_2,\iota_2]=\pm2\eta_2, \qquad \pi_3(S^2)=\ZZ\{\eta_2\}.\)
Therefore
\(2J_*(\eta_2)=2\eta_2.\)
Since $\pi_3(S^2)\cong\ZZ$ is torsion-free, multiplication by two is injective, so
\(J_*(\eta_2)=\eta_2.\)
Thus $J_*$ is the identity.

Finally, choose the identifications compatibly with the translations in \cref{lem:pi1-fiber}. Postcomposition satisfies
\(\ell_J(m(f,g))=m(\ell_J(f),\ell_J(g)),\)
because both sides apply $J$ to the value on each concatenated piece. Component translation, adjunction, and $\ell_J$ therefore form a commutative diagram, giving
\(H_{Jq_dJ}\circ(\ell_J)_*=J_*\circ H_{q_dJ}.\)
Since $J_*$ is the identity,
\(H_{Jq_dJ}\circ(\ell_J)_*=H_{q_dJ},\)
so $(\ell_J)_*$ is the identity under the standard identifications.
\end{proof}

\begin{theorem}[Conjugation and change of basepoint]\label{thm:reflection-conjugation}
Let $d$ be odd, $q\in O_d\cap F_d$, and $q'=JqJ$. For any homotopy $H_t:S^2\to S^2$ from $A$ to $J$, set
\[
\Theta_t(f)=H_tfH_t,\qquad \delta(t)=\Theta_t(q).
\]
Then $\delta$ is a path in $B_d$ from $q$ to $q'$. If
$\beta_\delta:\pi_1(B_d,q)\to\pi_1(B_d,q')$ denotes transport along $\delta$, then
\[
\beta_\delta^{-1}\circ(\Theta_J)_*=\Theta_*=-\id_{\pi_1(B_d,q)}.
\]
In the cyclic-group coordinates of \cref{thm:pi1-Bd}, this is $\Theta_*([n])=[-n]$.
\end{theorem}

\begin{proof}
We first prove that \(\Theta_J\) preserves \(B_d\).

Let
\(f\in B_d.\)
By definition,
\(\deg f=d.\)
Since
\(\deg J=-1,\)
the multiplicativity of degree gives
\[
\begin{aligned}
\deg(\Theta_J(f))
&=
\deg(J\circ f\circ J)\\
&=
\deg J\cdot\deg f\cdot\deg J\\
&=
(-1)d(-1)\\
&=
d.
\end{aligned}
\]
Hence
\(\Theta_J(f)\in B_d,\)
and therefore
\(\Theta_J(B_d)\subseteq B_d.\)

We next verify that \(q'\in F_d\). Since
\(q(N)=N, \qquad J(N)=N,\)
we have
\[
\begin{aligned}
q'(N)
&=
(J\circ q\circ J)(N)\\
&=
J(q(J(N)))\\
&=
J(q(N))\\
&=
J(N)\\
&=
N.
\end{aligned}
\]
Moreover,
\(\deg q' = (-1)d(-1) = d.\)
Thus
\(q'\in F_d.\)

Let
\(i:F_d\hookrightarrow B_d\)
denote the fiber inclusion.

By \cref{lem:precomp-natural}, precomposition with \(J\),
\(r_J:F_d\longrightarrow F_{-d}, \qquad r_J(f)=f\circ J,\)
satisfies
\(r_J(q)=q\circ J.\)
Under the standard identifications from \cref{lem:pi1-fiber},
\[
H_q:
\pi_1(F_d,q)\xrightarrow{\cong}\pi_3(S^2),
\]
and
\[
H_{qJ}:
\pi_1(F_{-d},qJ)\xrightarrow{\cong}\pi_3(S^2),
\]
\cref{lem:precomp-natural} gives
\(H_{qJ}\circ(r_J)_* = -\,H_q.\)

By \cref{lem:postcomp-natural}, postcomposition with \(J\),
\(\ell_J:F_{-d}\longrightarrow F_d, \qquad \ell_J(g)=J\circ g,\)
satisfies
\(\ell_J(qJ) = J\circ q\circ J = q'.\)
Moreover, under the standard identification
\[
H_{q'}:
\pi_1(F_d,q')
\xrightarrow{\cong}
\pi_3(S^2),
\]
\cref{lem:postcomp-natural} gives
\(H_{q'}\circ(\ell_J)_* = H_{qJ}.\)

Thus we obtain the commutative diagram
\[
\begin{array}{ccccc}
\pi_1(F_d,q)
&
\xrightarrow{\ (r_J)_*\ }
&
\pi_1(F_{-d},qJ)
&
\xrightarrow{\ (\ell_J)_*\ }
&
\pi_1(F_d,q')
\\[2mm]
{\scriptstyle H_q}\downarrow
&&
{\scriptstyle H_{qJ}}\downarrow
&&
{\scriptstyle H_{q'}}\downarrow
\\[2mm]
\pi_3(S^2)
&
\xrightarrow{\ -\operatorname{id}\ }
&
\pi_3(S^2)
&
\xrightarrow{\ \operatorname{id}\ }
&
\pi_3(S^2).
\end{array}
\]

Since
\(\Theta_J|_{F_d} = \ell_J\circ r_J,\)
we obtain
\[
\begin{aligned}
H_{q'}
\circ
(\Theta_J|_{F_d})_*
&=
H_{q'}
\circ
(\ell_J)_*
\circ
(r_J)_*\\
&=
H_{qJ}\circ(r_J)_*\\
&=
-\,H_q.
\end{aligned}
\]
Hence
\begin{equation}
H_{q'}
\circ
(\Theta_J|_{F_d})_*
=
-\,H_q.
\end{equation}

We now treat carefully the change of basepoint from \(q\) to \(q'\).

Since
\(q,q'\in F_d\)
and
\(\deg q=\deg q'=d,\)
and since \(F_d\) is path connected, there exists a path
\(\kappa:I\longrightarrow F_d\)
such that
\(\kappa(0)=q, \qquad \kappa(1)=q'.\)

For each \(v\in I\), define
\(b_v=\kappa(v)\circ J.\)
Since
\(\deg\kappa(v)=d\)
and
\(\deg J=-1,\)
we have
\(\deg b_v=-d.\)

Define
\(T_v:F_d\longrightarrow F_0, \qquad T_v(f)=m(b_v,f).\)
By \cref{lem:pi1-fiber},
\(\deg m(f,g) = \deg f+\deg g.\)
Hence, for every \(f\in F_d\),
\[
\begin{aligned}
\deg T_v(f)
&=
\deg m(b_v,f)\\
&=
\deg b_v+\deg f\\
&=
-d+d\\
&=
0.
\end{aligned}
\]
Therefore
\(T_v(F_d)\subseteq F_0.\)

By the component-translation result of \cref{lem:pi1-fiber},
\(T_v:F_d\longrightarrow F_0\)
is a homotopy equivalence for every \(v\).

Define
\(w:I\longrightarrow F_0, \qquad w(v)=T_v(\kappa(v)).\)
Then
\(w(0)=T_0(q), \qquad w(1)=T_1(q').\)

Let
\[
\beta_\kappa^F:
\pi_1(F_d,q)
\xrightarrow{\cong}
\pi_1(F_d,q')
\]
denote change of basepoint along \(\kappa\), and let
\[
\beta_w^0:
\pi_1(F_0,w(0))
\xrightarrow{\cong}
\pi_1(F_0,w(1))
\]
denote change of basepoint along \(w\).

Consider the continuous map
\[
\Phi:I\times I\longrightarrow F_0,
\qquad
\Phi(v,s)=T_v(\kappa(s)).
\]
Its four boundary restrictions are
\(\Phi(0,s)=T_0(\kappa(s)),\)
\[
\Phi(1,s)=T_1(\kappa(s)),
\]
\(\Phi(v,0)=T_v(q),\)
and
\(\Phi(v,1)=T_v(q').\)
Its diagonal is
\(\Phi(v,v) = T_v(\kappa(v)) = w(v).\)

Naturality of change of basepoint with respect to the family \(T_v\) therefore gives the commutative diagram
\[
\begin{array}{ccc}
\pi_1(F_d,q)
&
\xrightarrow{\ \beta_\kappa^F\ }
&
\pi_1(F_d,q')
\\[2mm]
{\scriptstyle (T_0)_*}\downarrow
&&
{\scriptstyle (T_1)_*}\downarrow
\\[2mm]
\pi_1(F_0,w(0))
&
\xrightarrow{\ \beta_w^0\ }
&
\pi_1(F_0,w(1)).
\end{array}
\]
Hence
\((T_1)_*\circ\beta_\kappa^F = \beta_w^0\circ(T_0)_*.\)

Let
\(c:S^2\longrightarrow S^2\)
be the constant based map.

Choose a path in \(F_0\),
\(\lambda_0:w(0)\rightsquigarrow c.\)
Define
\(\lambda_1=\overline w*\lambda_0: w(1)\rightsquigarrow c.\)
Then
\[
w*\lambda_1
=
w*\overline w*\lambda_0
\simeq_{\partial I}
\lambda_0.
\]
Consequently,
\(\beta_{\lambda_1}^0 \circ \beta_w^0 = \beta_{\lambda_0}^0.\)

Let
\[
H_c:
\pi_1(F_0,c)
\xrightarrow{\cong}
\pi_3(S^2)
\]
be the standard isomorphism obtained from the smash--function adjunction.

Compatibly with the component-translation construction of \cref{lem:pi1-fiber}, define
\(H_q = H_c \circ \beta_{\lambda_0}^0 \circ (T_0)_*\)
and
\(H_{q'} = H_c \circ \beta_{\lambda_1}^0 \circ (T_1)_*.\)

Then
\[
\begin{aligned}
H_{q'}\circ\beta_\kappa^F
&=
H_c
\circ
\beta_{\lambda_1}^0
\circ
(T_1)_*
\circ
\beta_\kappa^F\\
&=
H_c
\circ
\beta_{\lambda_1}^0
\circ
\beta_w^0
\circ
(T_0)_*\\
&=
H_c
\circ
\beta_{\lambda_0}^0
\circ
(T_0)_*\\
&=
H_q.
\end{aligned}
\]

Thus the following diagram commutes:
\[
\begin{array}{ccc}
\pi_1(F_d,q)
&
\xrightarrow{\ \beta_\kappa^F\ }
&
\pi_1(F_d,q')
\\[2mm]
{\scriptstyle H_q}\downarrow
&&
{\scriptstyle H_{q'}}\downarrow
\\[2mm]
\pi_3(S^2)
&
\xrightarrow{\ \operatorname{id}\ }
&
\pi_3(S^2).
\end{array}
\]

We now combine this with
\(H_{q'} \circ (\Theta_J|_{F_d})_* = -\,H_q.\)
For any
\[
u\in\pi_1(F_d,q),
\]
we have
\[
\begin{aligned}
H_{q'}
\left(
(\Theta_J|_{F_d})_*(u)
\right)
&=
-H_q(u)\\
&=
H_q(-u)\\
&=
H_{q'}
\left(
\beta_\kappa^F(-u)
\right).
\end{aligned}
\]
Since
\[
H_{q'}:
\pi_1(F_d,q')
\xrightarrow{\cong}
\pi_3(S^2)
\]
is an isomorphism, it follows that
\((\Theta_J|_{F_d})_*(u) = \beta_\kappa^F(-u).\)
Thus
\[
(\Theta_J|_{F_d})_*
=
\beta_\kappa^F\circ
\bigl(-\operatorname{id}\bigr).
\]

We next pass from the fiber \(F_d\) to the total space \(B_d\).

The long exact sequence of the evaluation fibration
\(F_d \longrightarrow B_d \xrightarrow{\operatorname{ev}} S^2\)
at the basepoint \(q\) contains
\[
\pi_1(F_d,q)
\xrightarrow{i_{*,q}}
\pi_1(B_d,q)
\longrightarrow
\pi_1(S^2,N).
\]
Since
\(\pi_1(S^2,N)=0,\)
the map
\(i_{*,q}: \pi_1(F_d,q) \twoheadrightarrow \pi_1(B_d,q)\)
is surjective.

Similarly,
\(i_{*,q'}: \pi_1(F_d,q') \twoheadrightarrow \pi_1(B_d,q')\)
is surjective.

Let
\[
\beta_\kappa^B:
\pi_1(B_d,q)
\xrightarrow{\cong}
\pi_1(B_d,q')
\]
denote change of basepoint along the same path
\(\kappa:I\longrightarrow B_d\)
after applying the inclusion \(F_d\hookrightarrow B_d\).

Naturality of change of basepoint with respect to
\(i:F_d\hookrightarrow B_d\)
gives the commutative diagram
\[
\begin{array}{ccc}
\pi_1(F_d,q)
&
\xrightarrow{\ \beta_\kappa^F\ }
&
\pi_1(F_d,q')
\\[2mm]
{\scriptstyle i_{*,q}}\downarrow
&&
{\scriptstyle i_{*,q'}}\downarrow
\\[2mm]
\pi_1(B_d,q)
&
\xrightarrow{\ \beta_\kappa^B\ }
&
\pi_1(B_d,q').
\end{array}
\]
Thus
\(\beta_\kappa^B\circ i_{*,q} = i_{*,q'}\circ\beta_\kappa^F.\)

On the other hand,
\(\Theta_J(F_d)\subseteq F_d\)
and
\(i\circ(\Theta_J|_{F_d}) = \Theta_J\circ i.\)
Hence we have the commutative diagram
\[
\begin{array}{ccc}
\pi_1(F_d,q)
&
\xrightarrow{\ (\Theta_J|_{F_d})_*\ }
&
\pi_1(F_d,q')
\\[2mm]
{\scriptstyle i_{*,q}}\downarrow
&&
{\scriptstyle i_{*,q'}}\downarrow
\\[2mm]
\pi_1(B_d,q)
&
\xrightarrow{\ (\Theta_J)_*\ }
&
\pi_1(B_d,q').
\end{array}
\]
Therefore
\[
(\Theta_J)_*\circ i_{*,q}
=
i_{*,q'}
\circ
(\Theta_J|_{F_d})_*.
\]

Let
\(u\in\pi_1(F_d,q).\)
Using the preceding two commutative diagrams together with
\((\Theta_J|_{F_d})_*(u) = \beta_\kappa^F(-u),\)
we obtain
\[
\begin{aligned}
(\Theta_J)_*
\bigl(i_{*,q}(u)\bigr)
&=
i_{*,q'}
\left(
(\Theta_J|_{F_d})_*(u)
\right)\\
&=
i_{*,q'}
\left(
\beta_\kappa^F(-u)
\right)\\
&=
\beta_\kappa^B
\left(
i_{*,q}(-u)
\right)\\
&=
\beta_\kappa^B
\left(
-\,i_{*,q}(u)
\right).
\end{aligned}
\]

Since
\(i_{*,q}: \pi_1(F_d,q) \twoheadrightarrow \pi_1(B_d,q)\)
is surjective, for every
\(x\in\pi_1(B_d,q)\)
there exists
\(u\in\pi_1(F_d,q)\)
such that
\(x=i_{*,q}(u).\)
Hence
\((\Theta_J)_*(x) = \beta_\kappa^B(-x).\)
Therefore
\[
(\beta_\kappa^B)^{-1}
\circ
(\Theta_J)_*
=
-\operatorname{id}_{\pi_1(B_d,q)}.
\]

It remains to compare \(\Theta_J\) with
\(\Theta(f)=A\circ f\circ A.\)

Since
\(\deg A=\deg J=-1\)
and
\([S^2,S^2]\cong\mathbb Z\)
is classified by degree, there exists a homotopy
\(H_t:S^2\longrightarrow S^2\)
such that
\(H_0=A, \qquad H_1=J.\)

Degree is invariant under homotopy, so
\(\deg H_t=-1 \qquad (t\in I).\)
Define
\[
\Theta_t:B_d\longrightarrow B_d,
\qquad
\Theta_t(f)=H_t\circ f\circ H_t.
\]
For every \(f\in B_d\),
\[
\begin{aligned}
\deg(\Theta_t(f))
&=
\deg H_t\cdot\deg f\cdot\deg H_t\\
&=
(-1)d(-1)\\
&=
d.
\end{aligned}
\]
Thus each
\(\Theta_t:B_d\longrightarrow B_d\)
is well defined.

Moreover,
\(\Theta_0=\Theta, \qquad \Theta_1=\Theta_J.\)

Since
\(q\in O_d,\)
we have
\(q(-x)=-q(x).\)
Writing
\[
A(x)=-x,
\]
this is equivalent to
\(q\circ A=A\circ q.\)
Therefore
\(\Theta(q) = A\circ q\circ A = q\circ A\circ A = q.\)

Define
\(\delta(t)=\Theta_t(q).\)
Then
\(\delta(0)=\Theta(q)=q\)
and
\(\delta(1)=\Theta_J(q)=q'.\)
Hence
\(\delta:q\rightsquigarrow q'.\)

We use the following convention for change of basepoint. If
\(\zeta:x\rightsquigarrow y\)
is a path in a space \(X\), define
\(\beta_\zeta: \pi_1(X,x) \longrightarrow \pi_1(X,y)\)
by
\(\beta_\zeta([\ell]) = [\overline{\zeta}*\ell*\zeta].\)

The homotopy
\(\Theta_t:B_d\longrightarrow B_d\)
has basepoint track
\[
t\longmapsto\Theta_t(q)=\delta(t).
\]
Naturality of fundamental groups under a homotopy with moving basepoint therefore gives the commutative diagram
\[
\begin{array}{ccc}
\pi_1(B_d,q)
&
\xrightarrow{\ \Theta_*\ }
&
\pi_1(B_d,q)
\\[2mm]
\operatorname{id}\downarrow
&&
{\scriptstyle\beta_\delta}\downarrow
\\[2mm]
\pi_1(B_d,q)
&
\xrightarrow{\ (\Theta_J)_*\ }
&
\pi_1(B_d,q').
\end{array}
\]
Thus
\((\Theta_J)_* = \beta_\delta\circ\Theta_*.\)

We now compare the two paths
\[
\delta,\kappa:q\rightsquigarrow q'.
\]
They have the same initial and terminal points.

For
\([\ell]\in\pi_1(B_d,q'),\)
we have
\((\beta_\kappa^B)^{-1}([\ell]) = [\kappa*\ell*\overline\kappa].\)
Applying \(\beta_\delta\), we obtain
\[
\begin{aligned}
\beta_\delta
\circ
(\beta_\kappa^B)^{-1}([\ell])
&=
[
\overline\delta*
\kappa*
\ell*
\overline\kappa*
\delta
]\\
&=
[
(\overline\delta*\kappa)
*
\ell
*
(\overline\delta*\kappa)^{-1}
].
\end{aligned}
\]
Thus
\(\beta_\delta \circ (\beta_\kappa^B)^{-1}\)
is an inner automorphism of
\(\pi_1(B_d,q').\)

By \cref{thm:pi1-Bd},
\[
\pi_1(B_d,q')
\cong
\mathbb Z/(2|d|)\mathbb Z.
\]
Hence
\(\pi_1(B_d,q')\)
is abelian, so every inner automorphism is the identity. Therefore
\(\beta_\delta \circ (\beta_\kappa^B)^{-1} = \operatorname{id},\)
and hence
\(\beta_\delta = \beta_\kappa^B.\)

Consequently,
\[
\begin{aligned}
\beta_\delta^{-1}
\circ
(\Theta_J)_*
&=
(\beta_\kappa^B)^{-1}
\circ
(\Theta_J)_*\\
&=
-\operatorname{id}_{\pi_1(B_d,q)}.
\end{aligned}
\]

Since
\((\Theta_J)_* = \beta_\delta\circ\Theta_*,\)
we also have
\(\Theta_* = \beta_\delta^{-1} \circ (\Theta_J)_*.\)
Therefore
\[
\Theta_*
=
\beta_\delta^{-1}
\circ
(\Theta_J)_*
=
-\operatorname{id}_{\pi_1(B_d,q)}.
\]

Finally, by \cref{thm:pi1-Bd},
\(\pi_1(B_d,q) \cong \mathbb Z/(2|d|)\mathbb Z,\)
and hence, under this identification,
\(\Theta_*([n]) = [-n].\)
This completes the proof.
\end{proof}

\begin{remark}[Precomposition versus postcomposition]\label{stmt:original-39}
Under the fiber identifications with $\pi_3(S^2)$, precomposition and postcomposition by $J$ act by $-1$ and $+1$, respectively. Their composite therefore induces inversion on $\pi_1(B_d)$ after basepoint transport.
\end{remark}

\subsubsection{An upper bound by the unique subgroup of order two}

\begin{lemma}[Upper bound for the odd endpoint image]\label{lem:Od-image-upper}
For odd $d$ and $q_d\in O_d\cap F_d$, the image of inclusion is contained in the subgroup of elements annihilated by $2$:
\begin{equation}\label{eq:Od-image-upper}
\im\bigl[\pi_1(O_d,q_d)\to\pi_1(B_d,q_d)\bigr]
\subseteq\{0,[|d|]\}\subset\ZZ/(2|d|)\ZZ.
\end{equation}
\end{lemma}

\begin{proof}
Represent an element $x$ of the image by a loop $\gamma$ in $O_d$. Since
$O_d=\operatorname{Fix}(\Theta)$, the equality $\Theta\circ\gamma=\gamma$ holds pointwise, giving
$\Theta_*x=x$. The conjugation theorem also gives $\Theta_*x=-x$, hence $2x=0$.
The only such elements in $\ZZ/(2|d|)\ZZ$ are
$0$ and $|d|\pmod{2|d|}$, proving \eqref{eq:Od-image-upper}.
\end{proof}

\subsubsection{A rotation loop realizes the nonzero element of order two}

\begin{definition}[Rotation loop]\label{def:rotation-loop}
Let $d$ be odd and choose $q_d\in O_d\cap F_d$. For rotation $R_\theta\in SO(3)$ through angle $\theta$ about the polar axis, define
\begin{equation}\label{eq:rotation-loop}
\gamma_d:[0,2\pi]\longrightarrow O_d\cap F_d,\qquad
\gamma_d(\theta)=R_\theta\circ q_d.
\end{equation}
This is a loop at $q_d$: $R_0=R_{2\pi}=\id$, $R_\theta(N)=N$, and $R_\theta(-x)=-R_\theta(x)$.
\end{definition}

\begin{lemma}[Connecting homomorphism for rotations]\label{lem:SO-connecting}
In the fibration
\begin{equation}\label{eq:SO-fibration}
SO(2)\longrightarrow SO(3)\xrightarrow{p}S^2,\qquad p(R)=R(N),
\end{equation}
where $SO(2)$ is the stabilizer of $N$, the connecting homomorphism
\[
\partial_{SO}:\pi_2(S^2,N)\longrightarrow\pi_1(SO(2),\id)
\]
is multiplication by $\pm2$ after identifying both groups with $\ZZ$.
\end{lemma}

\begin{proof}
Since $SO(3)\cong\RP^3$ has universal cover $S^3$,
\[
\pi_2(SO(3))=0,
\qquad
\pi_1(SO(3))\cong\ZZ/2,
\qquad
\pi_1(SO(2))\cong\ZZ.
\]
The long exact sequence gives
\[
0\longrightarrow\ZZ
\xrightarrow{\partial_{SO}}\ZZ
\longrightarrow\ZZ/2\longrightarrow0.
\]
Every injective homomorphism $\ZZ\to\ZZ$ is multiplication by a nonzero integer $m$, with cokernel of order $|m|$.
Here the cokernel is $\ZZ/2$, so $|m|=2$ and $\partial_{SO}=\times(\pm2)$.
\end{proof}

\begin{lemma}[Rotation-induced map of fibrations]\label{lem:rotation-fibration-map}
Let $d$ be odd and $q_d\in O_d\cap F_d$. The maps
\[
\alpha_d:SO(3)\to B_d,\quad R\mapsto Rq_d,\qquad
h_d:SO(2)\to F_d,\quad R\mapsto Rq_d,
\]
form the following strictly commutative map of fibrations:
\[
\begin{tikzcd}[column sep=4.2em]
SO(2) \arrow[r] \arrow[d,"h_d"'] & SO(3) \arrow[r,"p"] \arrow[d,"\alpha_d"] & S^2 \arrow[d,equal]\\
F_d \arrow[r] & B_d \arrow[r,"\operatorname{ev}"] & S^2.
\end{tikzcd}
\]
\end{lemma}

\begin{proof}
Since $q_d(N)=N$, for $R\in SO(3)$,
\(\operatorname{ev}(\alpha_d(R)) =(R\circ q_d)(N)=R(N)=p(R).\)
Hence $\operatorname{ev}\circ\alpha_d=p$. For $R\in SO(2)$, $R(N)=N$, so
$R\circ q_d\in F_d$, giving the fiber map $h_d$.
\end{proof}

\begin{theorem}[Rotation loop and odd endpoint image]\label{thm:rotation-loop-class}
Let $d$ be odd and $q_d\in O_d\cap F_d$. For the rotation loop $\gamma_d$ of \cref{def:rotation-loop}, its image in $\pi_1(B_d,q_d)\cong\ZZ/(2|d|)\ZZ$ is the unique nonzero element of order two:
\[
[\gamma_d]=[|d|].
\]
Consequently,
\begin{equation}\label{eq:Od-image}
\im\bigl[\pi_1(O_d,q_d)\to\pi_1(B_d,q_d)\bigr]
=\{0,[|d|]\}\cong\ZZ/2\ZZ.
\end{equation}
\end{theorem}

\begin{proof}
By \cref{thm:partial-d}, the connecting homomorphism of the lower evaluation fibration is
$\partial_d=\times(\pm2d)$; by \cref{lem:SO-connecting}, the upper one is
$\partial_{SO}=\times(\pm2)$. \cref{lem:rotation-fibration-map} and naturality of the long exact sequence
give the commutative square
\[
\begin{tikzcd}[column sep=5.5em]
\pi_2(S^2) \arrow[r,"\partial_{SO}"] \arrow[d,equal]
  & \pi_1(SO(2)) \arrow[d,"(h_d)_*"]\\
\pi_2(S^2) \arrow[r,"\partial_d"']
  & \pi_1(F_d,q_d).
\end{tikzcd}
\]
Write $(h_d)_*:\ZZ\to\ZZ$ as multiplication by $m$. Commutativity gives
$m(\pm2)=\pm2d$, hence $m=\pm d$. A generator of $\pi_1(SO(2))$ is a full rotation
$\theta\mapsto R_\theta$, whose image under $h_d$ is $\gamma_d$. Therefore
$\gamma_d$ represents $\pm d$ in $\pi_1(F_d)$.

The fiber inclusion induces the quotient map
\(\ZZ\longrightarrow\ZZ/(2|d|)\ZZ,\)
so
\([\gamma_d]=\pm d\pmod{2|d|}=|d|\pmod{2|d|}.\)
This is the unique nonzero element of order two. Since $\gamma_d$ lies in $O_d$, combining this with
\cref{lem:Od-image-upper} proves \eqref{eq:Od-image}.
\end{proof}

\begin{remark}[Independence of the order-two class]\label{stmt:original-45}
In $\ZZ/(2|d|)\ZZ$, the classes $[d]$ and $[-d]$ coincide and equal the unique nonzero element of order two. Thus sign choices for intermediate generators do not affect the rotation-loop image or the subgroup in \eqref{eq:Od-image}.
\end{remark}

\section{Classification, composition, and consequences}\label{sec:classification}

\begin{lemma}[The cohomology class and the endpoint indices]\label{lem:cohomology-endpoints}\label{rem:cohomology-endpoints}
Let $M=\mathbb{RP}^3\#\mathbb{RP}^3$ and $x_0\in M$, and consider the universal covering
\[
p:\widetilde M=S^2\times\mathbb R\longrightarrow M,\qquad p(\widetilde x_0)=x_0.
\]
Let $q:S^2\times\mathbb R\to S^2$ be the projection and let $u$ generate $H^2(S^2;\mathbb Z)$. If a pointed self-map $f:(M,x_0)\to(M,x_0)$ induces the trivial homomorphism on the fundamental group, then it has a unique pointed lift $\widehat f:M\to\widetilde M$. Put
\[
\alpha(f)=(q\circ\widehat f)^*u\in H^2(M;\mathbb Z).
\]
Then $\alpha(f)$ depends only on the pointed homotopy class of $f$. Under the isomorphism induced by restriction to the two endpoints,
\[
H^2(M;\mathbb Z)\cong(\mathbb Z/2\mathbb Z)^2,
\]
it corresponds to the component indices $(\varepsilon_0,\varepsilon_1)$ of the two even endpoint maps.
\end{lemma}

\begin{proof}
Let $p$ be the universal covering and let $\widetilde x_0\in\widetilde M$ satisfy $p:\widetilde M=S^2\times\mathbb R\longrightarrow M,\  p(\widetilde x_0)=x_0$. Suppose $f_{*,1}$ is trivial. Since $f_{*,1}(\pi_1(M,x_0))=\{1\}=p_*(\pi_1(\widetilde M,\widetilde x_0))$, the covering-space lifting criterion gives a unique pointed lift $\widehat f:M\longrightarrow\widetilde M,\  p\circ\widehat f=f,\ \widehat f(x_0)=\widetilde x_0$. 

By contrast, the lift used earlier is the unique pointed lift of $f\circ p$, $F_f:\widetilde M\longrightarrow\widetilde M,\  p\circ F_f=f\circ p$. Write $F_f(x,s)=(P_f(x,s),U_f(x,s))$. Since $\widehat f\circ p$ is also a pointed lift of $f\circ p$, uniqueness gives $F_f=\widehat f\circ p$. Let $q:S^2\times\mathbb R\to S^2$ be the projection, let $u$ generate $H^2(S^2;\mathbb Z)$, and put $h_f=q\circ\widehat f$. Define $\alpha([f])=h_f^*u$, also written $\alpha(f)$. Unique lifting of pointed homotopies makes this well defined, and $P_f=h_f\circ p,\ \alpha(f)=h_f^*u$. 

Use the quotient model
\[
M\cong\frac{S^2\times[0,1]}{(x,0)\sim(-x,0),\ (x,1)\sim(-x,1)}.
\]
Let $j_i:\mathbb{RP}^2\to M$, $j_i([x])=p(x,i)$, $i=0,1$, be the two endpoint inclusions, and let $p_2:S^2\to\mathbb{RP}^2$ be the antipodal quotient map. The identity $P_f=h_f\circ p$ gives
\[
P_f^i(x)=P_f(x,i)=h_f(p(x,i))=h_f(j_i([x])).
\]
Thus $P_f^i=(h_f\circ j_i)\circ p_2$. In particular, $P_f^i(x)=P_f^i(-x)$, so there is a unique continuous map $\overline P_f^{\,i}:\mathbb{RP}^2\longrightarrow S^2,\ \overline P_f^{\,i}([x])=P_f^i(x)$, such that $P_f^i=\overline P_f^{\,i}\circ p_2$. The same identities show that $\overline P_f^{\,i}=h_f\circ j_i$.

By \cref{thm:even-components},
\[
\pi_0(E)\cong[\mathbb{RP}^2,S^2]\cong H^2(\mathbb{RP}^2;\mathbb Z)\cong\mathbb Z/2\mathbb Z,
\]
where the cohomological classification sends $[v]$ to $v^*u$. The component $E_0$ contains the constant maps, and $E_1$ is the other component. Under $H^2(\mathbb{RP}^2;\mathbb Z)\cong\mathbb Z/2\mathbb Z$, therefore, $P_f^i\in E_{\varepsilon_i}\ \Longleftrightarrow\ (\overline P_f^{\,i})^*u=\varepsilon_i$. Consequently,
\[
j_i^*\alpha(f)=j_i^*(h_f^*u)=(h_f\circ j_i)^*u=(\overline P_f^{\,i})^*u=\varepsilon_i.
\]
Finally, restriction gives an isomorphism
\[
(j_0^*,j_1^*):H^2(M;\mathbb Z)\longrightarrow H^2(\mathbb{RP}^2;\mathbb Z)\oplus H^2(\mathbb{RP}^2;\mathbb Z).
\]
Indeed, use the open cover $M=U\cup V$ associated to the connected-sum decomposition in Section~2, where $U,V$ deformation retract onto the two endpoint copies of $\mathbb{RP}^2$ and $U\cap V\simeq S^2$. The cohomological Mayer--Vietoris sequence contains
\[
\begin{aligned}
H^1(U\cap V;\mathbb Z)&\longrightarrow H^2(M;\mathbb Z)\\
&\longrightarrow H^2(U;\mathbb Z)\oplus H^2(V;\mathbb Z)\longrightarrow H^2(U\cap V;\mathbb Z).
\end{aligned}
\]
Under the deformation-retraction isomorphisms this becomes
\[
0\longrightarrow H^2(M;\mathbb Z)\xrightarrow{(j_0^*,j_1^*)}(\mathbb Z/2\mathbb Z)^2\longrightarrow\mathbb Z.
\]
Every homomorphism from a finite group to $\mathbb Z$ is zero. Exactness therefore gives $(j_0^*,j_1^*):H^2(M;\mathbb Z)\xrightarrow{\cong}(\mathbb Z/2\mathbb Z)^2$. Hence $(j_0^*,j_1^*)\alpha(f)=(\varepsilon_0,\varepsilon_1)$.
If $f\simeq_*g$, lifting the pointed homotopy gives $\widehat f\simeq_*\widehat g$, hence $h_f\simeq_*h_g$. Thus $\alpha(f)=h_f^*u=h_g^*u=\alpha(g)$, verifying pointed homotopy invariance.

The four endpoint-component pairs correspond to the four values of $\alpha(f)$ as follows:
\[
\begin{aligned}
(E_0,E_0)&\longleftrightarrow(0,0),& (E_0,E_1)&\longleftrightarrow(0,1),\\
(E_1,E_0)&\longleftrightarrow(1,0),& (E_1,E_1)&\longleftrightarrow(1,1).
\end{aligned}
\]
Thus specifying the value of $\alpha(f)$ is equivalent to specifying the path components of the two endpoint maps.
\end{proof}

\subsection{Proof of the classification theorem}

\begin{proof}
We first clarify the path-space notation. When the endpoint spaces
$A_0,A_1$ are subspaces of $B$, a point of the space in
\cref{lem:double-coset-path},
$\mathcal H(A_0,A_1;B)$, is initially written as a triple
$(a_0,\gamma,a_1)$ with
\(\gamma(0)=a_0\in A_0, \qquad \gamma(1)=a_1\in A_1.\)
Since the path determines its endpoints, forgetting their coordinates
gives a natural homeomorphism
\[
\mathcal H(A_0,A_1;B)
\cong
\mathcal P(A_0,A_1;B)
:=
\{\gamma:I\to B\mid
\gamma(0)\in A_0,\ \gamma(1)\in A_1\}.
\]
We therefore identify the two path spaces below.

Fix $\varphi$. The correspondence between pointed maps and equivariant
lifts fixing the lifted basepoint,
\cref{thm:remove-lift-basepoint},
\cref{lem:R-coordinate-no-pi0},
and \cref{lem:strip-path-correspondence}
give successive natural bijections
\[
\begin{aligned}
[M,M]_*^\varphi
&\cong\pi_0(\mathcal F_\varphi)
\cong\pi_0(\mathcal E_\varphi)\\
&\cong\pi_0(\mathcal E_{\lambda_\varphi})
\cong\pi_0\mathcal P(A_0,A_1;B)\\
&\cong\pi_0\mathcal H(A_0,A_1;B).
\end{aligned}
\]
Here $B=\Map(S^2,S^2)$, and $A_0,A_1$ are determined by
$\lambda_\varphi(a),\lambda_\varphi(b)$.
Paths in this mapping space correspond exactly to homotopies of paths
with endpoints remaining in $A_0,A_1$.
Thus computing $\pi_0\mathcal H$, after separating degrees and endpoint
components, yields the final classification.
We compute
\(\pi_0\mathcal H(E_{\varepsilon_0},E_{\varepsilon_1};B_0)\)
and
\(\pi_0\mathcal H(O_d,O_d;B_d)\)
using the double-coset formula of
\cref{lem:double-coset-path}.

Before computing induced maps on $\pi_2$, we fix the basepoints
throughout.
By \cref{thm:remove-lift-basepoint}, choose in every nonempty
equivariant component a representative
\[
F=(P,U)
\]
satisfying
\[
F(\widetilde x_0)=\widetilde x_0.
\]
Its descent $f$ is pointed and satisfies
\(f_{*,1}=\varphi.\)
An equivariant homotopy used to make this choice restricts, at every
fixed $u$, to an ordinary homotopy of the sphere-coordinate maps
\(P_u:S^2\to S^2.\)
Thus this replacement changes neither their degree nor the path-space
component.
All subsequent based homotopy-group diagrams use this pointed
representative.

Here is a basepoint-explicit justification for computing on an interior
slice.
Write
\(\widetilde x_0=(x_b,u_b)\)
and fix
\(u_0\in(0,1).\)
The path
\(\tau(v) = \bigl(x_b,(1-v)u_b+vu_0\bigr)\)
joins $\widetilde x_0$ to
\[
z_0=(x_b,u_0).
\]
The slice inclusion
\(i_{u_0}:(S^2,x_b)\longrightarrow(\widetilde M,z_0)\)
is a pointed homotopy equivalence.
Naturality of basepoint transport for $F$ compares $F_*$ at
$\widetilde x_0$ with the map
\[
\pi_2(\widetilde M,z_0)
\longrightarrow
\pi_2(\widetilde M,F(z_0)),
\]
using $\tau$ and $F\tau$.
After projection to $S^2$, this latter map is induced by
\(P_{u_0}: (S^2,x_b) \longrightarrow (S^2,P_{u_0}(x_b)).\)
Both $\widetilde M$ and $S^2$ are simply connected, so any two paths
with the same endpoints give the same basepoint-change isomorphism on
$\pi_2$.
Use a fixed orientation of $S^2$ to identify these groups with $\ZZ$
via the Hurewicz isomorphism; projection induces the corresponding
identification for $\widetilde M$.
Basepoint transport preserves these homology identifications.
Consequently the multiplier of
\[
F_*:
\pi_2(\widetilde M,\widetilde x_0)
\longrightarrow
\pi_2(\widetilde M,\widetilde x_0)
\]
is exactly
\(\deg P_{u_0}.\)
Since
\(pF=fp\)
is an equality of pointed maps, the covering-induced isomorphism
\[
p_*:
\pi_2(\widetilde M,\widetilde x_0)
\longrightarrow
\pi_2(M,x_0)
\]
identifies this multiplier with that of $f_{*,2}$.
This justifies the slice calculations in both cases below.

\medskip
\noindent
\textbf{(1) The trivial fundamental-group endomorphism.}
\par\medskip

Suppose
\(\varphi(g)=1\quad(g\in G),\)
that is,
\(\varphi(a)=\varphi(b)=1.\)
Then
\(\lambda_\varphi(a)=\lambda_\varphi(b)=+1,\)
so both endpoints lie in $E$.
Since
\(E=E_0\sqcup E_1,\)
there are four possible endpoint-component pairs
\((\varepsilon_0,\varepsilon_1)\in(\ZZ/2)^2.\)
For a fixed pair, the double-coset lemma gives
\[
\pi_0\mathcal H(E_{\varepsilon_0},E_{\varepsilon_1};B_0)
\cong
\frac{\pi_1(B_0)}
{\im\pi_1(E_{\varepsilon_0})
+\im\pi_1(E_{\varepsilon_1})}.
\]
By \cref{thm:pi1-Bd} and \cref{thm:E-image-zero},
\[
\pi_1(B_0)\cong\ZZ,
\qquad
\im\pi_1(E_{\varepsilon_i})=0.
\]
Hence
\begin{equation}\label{eq:zero-phi-path-components}
\pi_0\mathcal H(E_{\varepsilon_0},E_{\varepsilon_1};B_0)
\cong\ZZ.
\end{equation}
Therefore
\begin{equation}\label{eq:zero-phi-classification}
[M,M]_*^{\varphi,0}
=
\left\{
[f]\in[M,M]_*
\mid
f_{*,1}=\varphi
\right\}
\cong
\bigsqcup_{(\varepsilon_0,\varepsilon_1)\in(\ZZ/2)^2}\ZZ.
\end{equation}
Each summand carries a free and transitive $\ZZ$-action; the displayed bijection depends on the choice of a reference path for each endpoint pair.

We also record the induced higher-dimensional information.
Since
\(P_0,P_1\in E,\)
\cref{lem:even-degree-zero} gives
\(\deg P_0=\deg P_1=0.\)
The path
\(u\longmapsto P_u\)
is an ordinary homotopy, so
\(\deg P_u=0 \qquad (u\in I).\)

Fix
\(u_0\in(0,1).\)
The inclusion
\[
i_{u_0}:S^2\longrightarrow S^2\times\RR,
\qquad
i_{u_0}(x)=(x,u_0)
\]
and projection
\(\operatorname{pr}_{S^2}:S^2\times\RR\longrightarrow S^2\)
are homotopy equivalences.
For the lift
\[
F=(P,U),
\]
we have
\(\operatorname{pr}_{S^2}\circ F\circ i_{u_0} = P_{u_0}.\)
Thus $F_*$ on
\(\pi_2(\widetilde M)\cong\ZZ\)
corresponds to
\[
(P_{u_0})_*
=
\times\deg P_{u_0}
=
0.
\]
Using the covering-induced isomorphism
\(\pi_2(\widetilde M)\cong\pi_2(M)\)
yields
\begin{equation}\label{eq:zero-phi-pi2}
f_{*,2}:
\pi_2(M)\longrightarrow\pi_2(M)
\quad\text{is the zero map}.
\end{equation}

The ordinary degree is also zero.
Since
\(f_{*,1}=\varphi,\)
the lifting criterion gives
\[
\bar f:M\longrightarrow\widetilde M
\qquad\text{such that}\qquad
f=p\circ\bar f.
\]
Since
\(\widetilde M\simeq S^2,\)
we have
\(H_3(\widetilde M;\ZZ)=0.\)
Thus
\(\bar f_*: H_3(M;\ZZ) \longrightarrow H_3(\widetilde M;\ZZ)\)
is zero, and hence
\(f_*=p_*\bar f_*\)
vanishes on $H_3$.
Therefore
\begin{equation}\label{eq:zero-phi-degree}
\deg f=0.
\end{equation}
The integer secondary coordinate in the trivial fundamental-group case is
consequently not the ordinary degree.

\medskip
\noindent
\textbf{(2) Endomorphisms sending both generators to reflections.}
\par\medskip

Suppose
\(\varphi(t)=t^m,\qquad \varphi(a)=t^ka,\qquad \varphi(b)=t^{k-m}a.\)
Then
\(\lambda_\varphi(a)=\lambda_\varphi(b)=-1,\)
so both endpoints lie in $O$.
By \cref{thm:odd-map-classification}, they lie in some
\[
O_d,\qquad O_{d'}.
\]
Since
\(u\longmapsto P_u\)
is an ordinary homotopy, degree is constant, so
\(d=d',\)
and this integer is odd.

This $d$ is precisely the multiplier of $f_{*,2}$ on
\(\pi_2(M)\cong\ZZ.\)
Indeed, choose an interior slice
\(S^2\times\{u_0\}\subset\widetilde M.\)
It represents a generator of $\pi_2(\widetilde M)$ under
\(\widetilde M\simeq S^2.\)
On this slice the lift
\(F=(P,U)\)
has real coordinate contributing nothing to $\pi_2$, and sphere
coordinate $P_{u_0}$ of degree $d$.
Hence
\[
F_*:
\pi_2(\widetilde M)
\longrightarrow
\pi_2(\widetilde M)
\quad\text{is multiplication by }d.
\]

The equality
\(p\circ F=f\circ p\)
gives the commutative diagram
\[
\begin{tikzcd}[column sep=5em,row sep=3em]
\pi_2(\widetilde M)
\arrow[r,"F_*"]
\arrow[d,"p_*"',"\cong"]
&
\pi_2(\widetilde M)
\arrow[d,"p_*","\cong"']
\\
\pi_2(M)
\arrow[r,"f_{*,2}"']
&
\pi_2(M).
\end{tikzcd}
\]
Thus
\(p_*\circ F_* = f_{*,2}\circ p_*,\)
so
\begin{equation}\label{eq:pi2-degree-d}
f_{*,2}
=
d\,\mathrm{id}_{\ZZ}.
\end{equation}

Fix $\varphi$ as above and odd $d$.
The double-coset lemma gives
\[
\pi_0\mathcal H(O_d,O_d;B_d)
\cong
\frac{\pi_1(B_d)}
{\im\pi_1(O_d)+\im\pi_1(O_d)}.
\]
By \cref{thm:pi1-Bd} and \cref{thm:rotation-loop-class},
\[
\pi_1(B_d)
=
\ZZ/(2|d|),
\qquad
\im\pi_1(O_d)
=
\langle|d|\rangle
=
\{0,|d|\}.
\]
Therefore
\[
\pi_0\mathcal H(O_d,O_d;B_d)
\cong
\frac{\ZZ/(2|d|)}
{\langle|d|\rangle}
\cong
\ZZ/|d|.
\]
That is,
\begin{equation}\label{eq:double-reflection-classes}
[M,M]_*^{\varphi,d}
=
\left\{
[f]\in[M,M]_*
\ \middle|\
f_{*,1}=\varphi,
\quad
f_{*,2}=d\,\mathrm{id}_{\ZZ}
\right\}
\end{equation}
carries a free and transitive $\ZZ/|d|$-action.
Choosing a reference map gives a set bijection
\[
[M,M]_*^{\varphi,d}
=
\left\{
[f]\in[M,M]_*
\ \middle|\
f_{*,1}=\varphi,
\quad
f_{*,2}=d\,\mathrm{id}_{\ZZ}
\right\}
\cong
\ZZ/|d|.
\]

This also explains why the answer is not
\(\ZZ/(2|d|).\)
If the two endpoints are fixed, the loop parameter lies in
\(\pi_1(B_d)=\ZZ/(2|d|).\)
In the actual equivariant homotopy problem the endpoints may move in
$O_d$, whose rotation loop contributes the subgroup
\[
\{0,|d|\}.
\]
We must therefore divide out this subgroup of order two.

\medskip
\noindent
\textbf{(3) Every listed parameter is realizable.}
\par\medskip

We show that the preceding conditions are sufficient as well as
necessary.

First consider the trivial fundamental-group endomorphism.
Fix an endpoint-component pair
\((\varepsilon_0,\varepsilon_1)\in(\ZZ/2)^2\)
and choose representatives
\[
q_{\varepsilon_0}\in E_{\varepsilon_0},
\qquad
q_{\varepsilon_1}\in E_{\varepsilon_1}.
\]
Since $B_0$ is path connected, choose a reference path
\[
\eta:I\longrightarrow B_0,
\qquad
\eta(0)=q_{\varepsilon_0},
\quad
\eta(1)=q_{\varepsilon_1}.
\]
Take any
\(n\in \pi_1(B_0,q_{\varepsilon_0}) \cong\ZZ,\)
and choose a loop at $q_{\varepsilon_0}$ representing $n$,
\[
\ell_n:I\longrightarrow B_0,
\qquad
\ell_n(0)=\ell_n(1)=q_{\varepsilon_0}.
\]
Concatenate this loop with the reference path to define
\(\gamma_n=\ell_n*\eta,\)
explicitly,
\[
\gamma_n(u)=
\begin{cases}
\ell_n(2u),
&0\le u\le\frac12,\\[2mm]
\eta(2u-1),
&\frac12\le u\le1.
\end{cases}
\]
Then
\[
\gamma_n(0)
=
q_{\varepsilon_0}
\in E_{\varepsilon_0},
\qquad
\gamma_n(1)
=
q_{\varepsilon_1}
\in E_{\varepsilon_1},
\]
so
\[
\gamma_n
\in
\mathcal H(E_{\varepsilon_0},E_{\varepsilon_1};B_0).
\]

Distinct integers give distinct path components.
Indeed, if
\([\gamma_n]=[\gamma_m],\)
the double-coset lemma implies that
\(n-m\)
belongs to
\[
\im\bigl[
\pi_1(E_{\varepsilon_0})
\to
\pi_1(B_0)
\bigr]
+
\im\bigl[
\pi_1(E_{\varepsilon_1})
\to
\pi_1(B_0)
\bigr].
\]
Both images are zero, as proved above, so
\(n=m.\)
Thus
\(n\longmapsto[\gamma_n]\)
gives distinct classes for this endpoint pair.
Conversely, move the endpoints of any admissible path within the
endpoint spaces to the chosen representatives, then concatenate with
$\bar\eta$ to obtain a loop at $q_{\varepsilon_0}$.
The double-coset lemma therefore ensures that every class is represented
by some $\gamma_n$.
This is a bijection.
Choosing $\eta$ chooses the zero of the integer parameter; without that
choice the set still carries a free and transitive $\ZZ$-action, but no origin for the integer coordinate has been specified.

Apply the fundamental-strip extension lemma to $\gamma_n$ to obtain a
$\lambda_\varphi$-equivariant map
\[
P_n:S^2\times\RR\longrightarrow S^2.
\]
For the trivial fundamental-group endomorphism,
\(\lambda_\varphi(a)=\lambda_\varphi(b)=+1,\)
and we may choose
\(U=0.\)
Then
\[
F_n=(P_n,0):
\widetilde M=S^2\times\RR
\longrightarrow
\widetilde M
\]
satisfies the equivariance relation corresponding to the trivial
fundamental-group endomorphism.
By \cref{thm:remove-lift-basepoint}, its equivariant homotopy
component contains a representative fixing $\widetilde x_0$, which
descends to a pointed self-map of $M$.
Thus in the trivial fundamental-group case every endpoint pair and every
integer parameter is realized.

Now consider the reflection case.
Given
\[
m,k\in\ZZ,
\qquad
d\in\ZZ\text{ odd},
\qquad
c\in\ZZ/|d|\ZZ,
\]
Rubinsztein's theorem supplies an odd degree-$d$ map
\[
q_d\in O_d.
\]
We have already proved
\(\pi_1(B_d,q_d) \cong \ZZ/(2|d|)\ZZ\)
and the path-space classification is given by the quotient
\(\rho_d: \ZZ/(2|d|)\ZZ \twoheadrightarrow \ZZ/|d|\ZZ.\)
Since $\rho_d$ is surjective, choose
\(\widetilde c \in \ZZ/(2|d|)\ZZ, \qquad \rho_d(\widetilde c)=c.\)
Take a loop at $q_d$ representing $\widetilde c$,
\[
\ell_{\widetilde c}:I\longrightarrow B_d,
\qquad
\ell_{\widetilde c}(0)
=
\ell_{\widetilde c}(1)
=
q_d.
\]
Since
\(q_d\in O_d,\)
both endpoints belong to $O_d$, so this loop defines a class in
\(\mathcal H(O_d,O_d;B_d)\)
corresponding exactly to $c$ in the quotient.

Put
\(P_u = \ell_{\widetilde c}(u): S^2\longrightarrow S^2.\)
Only the endpoints are required to lie in $O_d$; intermediate maps need
only lie in $B_d$.
Following \cref{lem:strip-path-correspondence}, set
\(Q(x,u)=P_u(x)\)
on the fundamental strip and define
\(P(gz_0) = \lambda_\varphi(g)Q(z_0), \qquad z_0\in S^2\times I.\)
Oddness of the endpoints ensures compatibility on common boundaries.
This equivariantly extends the path to
\(S^2\times\RR,\)
giving a $\lambda_\varphi$-equivariant map
\[
P:S^2\times\RR\longrightarrow S^2.
\]

Choose the affine-equivariant function from \cref{lem:R-coordinate-no-pi0},
\(U(x,u) = mu-k.\)
Then
\(F=(P,U): \widetilde M \longrightarrow \widetilde M\)
satisfies
\(F(gz) = \varphi(g)F(z) \qquad (g\in G,\ z\in\widetilde M).\)
Thus $F$ descends to a self-map of $M$.
Finally, by \cref{thm:remove-lift-basepoint}, choose in its
equivariant homotopy component a representative fixing the lifted
basepoint.
Its descent is pointed and satisfies
\(f_{*,1}=\varphi, \qquad f_{*,2}=d\,\mathrm{id}_{\ZZ},\)
with path-space parameter $c$.
Thus every tuple
\((m,k,d,c)\)
is realizable.

\cref{lem:mixed-not-realizable}
excludes the two mixed types.
The trivial fundamental-group and reflection cases have both been
computed and realized.
These exhaust all possibilities, completing the proof.
\end{proof}
\begin{remark}[Reference homotopy classes and secondary coordinates]\label{stmt:original-47}
Choosing a reference homotopy class in each family in the trivial case and in each nonempty classification set in the double-reflection case gives a noncanonical bijection
\[
[M,M]_*\cong
\left(\bigsqcup_{(\varepsilon_0,\varepsilon_1)\in(\ZZ/2\ZZ)^2}\ZZ\right)
\sqcup
\left(\bigsqcup_{\substack{m,k\in\ZZ\\d\in2\ZZ+1}}\ZZ/|d|\ZZ\right).
\]

Here $\varphi$ is determined by $(m,k)$ as in part (2) of the theorem. Each family in the trivial case carries a free and transitive action of $A=\mathbb Z$. In the double-reflection case, each nonempty classification set for prescribed $(\varphi,d)$ carries a free and transitive action of $A=\mathbb Z/|d|\mathbb Z$. These actions are induced by concatenation of mapping-space loops, with the endpoint images divided out. After fixing the identifications of the acting groups, choosing a reference homotopy class gives a bijection between the classification set and $A$; changing that class translates the coordinate. In particular, $|[M,M]_*^{\varphi,d}|=|d|$, giving one, three and five classes for $d=\pm1,\pm3,\pm5$, respectively. The secondary coordinate is not determined by $f_{*,1}$ and $f_{*,2}$.
\end{remark}

\subsection{Compatible reference maps and secondary coordinates}

\begin{lemma}[Multiplicative odd reference family]\label{lem:compatible-odd-reference-family}
There exists a family $(q_m)_{m\in2\ZZ+1}$ of maps $S^2\to S^2$ such that
\[
\deg q_m=m,\qquad q_m(-x)=-q_m(x),\qquad q_m(N)=N,
\]
for every odd $m$, and
\[
q_1=\id_{S^2},\qquad q_d\circ q_e=q_{de}\quad(d,e\in2\ZZ+1).
\]
In particular, $q_m\in O_m\cap F_m$.
\end{lemma}

\begin{proof}
Use the \emph{unreduced} suspension

\[
S(S^1)=\bigl(S^1\times[-1,1]\bigr)/\sim,
\]
where $S^1\times\{1\}$ is collapsed to one point and $S^1\times\{-1\}$ to another; there are no other identifications. The map

\[
\Psi:S(S^1)\longrightarrow S^2\subset\mathbb C\times\mathbb R,
\qquad \Psi[z,s]=(\sqrt{1-s^2}\,z,s)
\]
is well defined and continuous, because its defining map on the cylinder is continuous and constant on each collapsed end. It is bijective: for $|s|<1$ the first coordinate determines $z$ uniquely, and at $s=\pm1$ there is exactly one suspension vertex. Its domain is compact and its target is Hausdorff, so it is a homeomorphism. The upper vertex corresponds to $N=(0,1)$.

In these coordinates the antipodal involution is

\[
\widehat A[z,s]=[-z,-s],\qquad
\Psi\widehat A[z,s]=(-\sqrt{1-s^2}\,z,-s)=A\Psi[z,s].
\]
In particular, it exchanges the two vertices. It is not the involution $[z,s]\mapsto[-z,s]$.
For every odd integer $m$, the continuous map $(z,s)\mapsto(z^m,s)$ preserves each collapsed end, so it descends to a continuous map

\[
\widehat q_m:S(S^1)\to S(S^1),\qquad
\widehat q_m[z,s]=[z^m,s],\qquad
q_m=\Psi\widehat q_m\Psi^{-1}.
\]
It fixes both vertices, so $q_m(N)=N$. For all $z\in S^1$ and $s\in[-1,1]$, oddness of $m$ gives

\[
\widehat q_m\widehat A[z,s]=[(-z)^m,-s]
=[-z^m,-s]=\widehat A\widehat q_m[z,s].
\]
\[
\begin{tikzcd}[column sep=large,row sep=large]
S^2\arrow[r,"q_m"]\arrow[d,"A"'] & S^2\arrow[d,"A"]\\
S^2\arrow[r,"q_m"'] & S^2.
\end{tikzcd}
\]
Thus $q_m(-x)=-q_m(x)$. To compute its degree, put $\mu_m(z)=z^m$. The loop $t\mapsto e^{2\pi it}$ is sent to $t\mapsto e^{2\pi imt}$, so $(\mu_m)_*$ is multiplication by $m$ on $H_1(S^1;\mathbb Z)$. The suspension isomorphism $\sigma$, obtained from the reduced Mayer--Vietoris sequence for the two cones, is natural. Hence

\[
\begin{tikzcd}[column sep=large,row sep=large]
\widetilde H_1(S^1;\mathbb Z)\arrow[r,"(\mu_m)_*=\times m"]\arrow[d,"\sigma"',"\cong"] & \widetilde H_1(S^1;\mathbb Z)\arrow[d,"\sigma","\cong"']\\
\widetilde H_2(S(S^1);\mathbb Z)\arrow[r,"(\widehat q_m)_*"'] & \widetilde H_2(S(S^1);\mathbb Z).
\end{tikzcd}
\]
Since both vertical maps are isomorphisms, $(\widehat q_m)_*$ is multiplication by $m$. Conjugating by $\Psi_*$ shows that $(q_m)_*$ is also multiplication by $m$, so $\deg q_m=m$. This calculation includes negative odd integers.

Finally, for every pair of odd integers $d,e$,

\[
\widehat q_d\widehat q_e[z,s]=[(z^e)^d,s]
=[z^{de},s]=\widehat q_{de}[z,s].
\]
The equality also holds at the vertices, which all these maps fix. Conjugating by the same homeomorphism $\Psi$ gives the strictly commutative triangle

\[
\begin{tikzcd}[column sep=large,row sep=large]
S^2\arrow[r,"q_e"]\arrow[dr,"q_{de}"'] & S^2\arrow[d,"q_d"]\\
& S^2.
\end{tikzcd}
\]
Thus $q_dq_e=q_{de}$ as actual maps, not merely up to homotopy. Since $z^1=z$, the same construction gives $q_1=\mathrm{id}_{S^2}$. All the required properties follow.
\end{proof}

\begin{notation}[Notation and coordinates for composition]\label{def:secondary-coordinates}\label{not:section11-common}
Retain $M,G,t$ and $[M,M]_*^{\varphi,d}$ from \cref{thm:classification}. Put
\[
\mathscr M=\left\{A(m,k)=\begin{pmatrix}m&k\\0&1\end{pmatrix}\;\middle|\;m,k\in\mathbb Z\right\},
\qquad \varphi_{A(m,k)}(t)=t^m,\quad\varphi_{A(m,k)}(a)=t^ka.
\]
Let $\varphi_{\mathrm{triv}}$ denote the trivial homomorphism, and put
\[
\mathcal R=\{[f]\in[M,M]_*\mid f_{*,1}\ne\varphi_{\mathrm{triv}}\},
\qquad\mathcal T=[M,M]_*^{\varphi_{\mathrm{triv}},0}.
\]

Fix a family $(q_d)$ as in \cref{lem:compatible-odd-reference-family}. Let $\iota_2\in\pi_2(S^2,N)$ be the class of the identity map and choose generators $\xi_d\in\pi_1(F_d,q_d)$ so that $\partial_d(\iota_2)=2d\xi_d$. The image of $\xi_d$ identifies $\pi_1(B_d,q_d)$ with $\ZZ/(2|d|)\ZZ$.
For every $(\mathbf A,d)$, take the class of the constant path at $q_d$ as origin in the path-space model for $[M,M]_*^{\varphi_{\mathbf A},d}$. The quotient
\[
\ZZ/(2|d|)\ZZ\longrightarrow\ZZ/|d|\ZZ
\]
then defines its secondary coordinate $c$. Write the resulting class as $(\mathbf A,d,c)$.
Write $[f]\star[g]=[f\circ g]$ and let $0$ denote the constant-map class. Write $[n]_r$ for reduction modulo $r$ and $I_2$ for the identity matrix.
\end{notation}

\subsection{Composition law and coordinate monoid}

\begin{theorem}[Composition in the double-reflection part]\label{cor:double-reflection-composition}
If $[f]=(\mathbf A,d,c)$ and $[g]=(\mathbf B,e,c')$, then
\[
[f]\star[g]=\bigl(\mathbf A\mathbf B,de,
[e\det(\mathbf B)\widetilde c+d^2\widetilde c']_{|de|}\bigr),
\]
where $\widetilde c,\widetilde c'$ are arbitrary integer representatives; the result is independent of their choice. Thus $\mathcal R$ is a submonoid of $[M,M]_*$ with identity $(I_2,1,0)$.
\end{theorem}

\begin{proof}
By \cref{thm:pi1-Bd},
\(\pi_1(B_m,q_m)\cong\ZZ/(2|m|)\ZZ\)
for every nonzero odd integer $m$, and by \cref{thm:rotation-loop-class},
\[
\operatorname{im}\!\left(
\pi_1(O_m,q_m)\longrightarrow\pi_1(B_m,q_m)
\right)
=
\{0,|m|\}.
\]
Hence the path-space double-coset lemma gives
\[
\pi_0\mathcal H(O_m,O_m;B_m)
\cong
\frac{\ZZ/(2|m|)\ZZ}{\langle|m|\rangle}
\cong
\ZZ/|m|\ZZ.
\]
From now on we fix the reference maps, generators, and quotient isomorphisms above; with these choices, the elements of $[M,M]_*^{\varphi_{\mathbf A},m}$ are labelled by
\(c\in\ZZ/|m|\ZZ.\)

Put $\mathbf A=A(m,k)$ and $\mathbf B=A(n,\ell)$. The relation $ata=t^{-1}$ gives $at^ja=t^{-j}$. On the generators we obtain
\[
(\varphi_{\mathbf A}\circ\varphi_{\mathbf B})(t)=t^{mn},\qquad
(\varphi_{\mathbf A}\circ\varphi_{\mathbf B})(a)=\varphi_{\mathbf A}(t^\ell a)=t^{m\ell+k}a .
\]
Since $b=t^{-1}a$, the second endpoint has image
\[
(\varphi_{\mathbf A}\circ\varphi_{\mathbf B})(b)
=\varphi_{\mathbf A}(t^{\ell-n}a)=t^{m(\ell-n)+k}a .
\]
\[
\mathbf A\mathbf B=A(mn,k+m\ell),\qquad
\varphi_{\mathbf A}\circ\varphi_{\mathbf B}=\varphi_{\mathbf A\mathbf B}.
\]
By functoriality of homotopy groups,
\((f\circ g)_{*,2} =f_{*,2}\circ g_{*,2} =de\,\id_{\ZZ}.\)
Hence $[f\circ g]\in[M,M]_*^{\varphi_{\mathbf A\mathbf B},de}$, and it remains to compute its secondary coordinate.

Choose a loop
\(\alpha:I\longrightarrow B_d, \qquad \alpha(0)=\alpha(1)=q_d,\)
representing $c$, and a loop
\(\beta:I\longrightarrow B_e, \qquad \beta(0)=\beta(1)=q_e,\)
representing $c'$. By \cref{lem:strip-path-correspondence}, extend them equivariantly to
\(P,Q:S^2\times\RR\longrightarrow S^2.\)
In the double-reflection case,
\(\lambda_{\varphi_{\mathbf A}}=\rho, \qquad \lambda_{\varphi_{\mathbf B}}=\rho,\)
because the images of both generators are reflections and $a,b$ generate $G$. By \cref{lem:R-coordinate-no-pi0}, choose the standard affine second coordinates
\(U_{\mathbf A}(u)=mu-k, \qquad U_{\mathbf B}(u)=nu-\ell.\)
Set
\(F=(P,U_{\mathbf A}), \qquad G=(Q,U_{\mathbf B}).\)
Strictly speaking, $F$ and $G$ are representatives in the corresponding unbased equivariant components. Using the previously established bijections
$\pi_0(\mathcal F_\varphi)\cong\pi_0(\mathcal E_\varphi)$, choose based representatives
$F^{\mathrm b},G^{\mathrm b}$ and equivariant paths
$F_t:F^{\mathrm b}\rightsquigarrow F$ and
$G_t:G^{\mathrm b}\rightsquigarrow G$. Then
\(K_t:=F_t\circ G_t\)
is $(\varphi_{\mathbf A}\circ\varphi_{\mathbf B})$-equivariant, because for every $h\in G$,
\[
\begin{aligned}
K_t(hz)
&=F_t(G_t(hz))
=F_t(\varphi_{\mathbf B}(h)G_t(z))\\
&=\varphi_{\mathbf A}(\varphi_{\mathbf B}(h))F_t(G_t(z))
=(\varphi_{\mathbf A}\circ\varphi_{\mathbf B})(h)K_t(z).
\end{aligned}
\]
Thus computing the equivariant component of $F\circ G$ computes the pointed class $[f\circ g]$.

Set
\(V(u)=U_{\mathbf B}(u)=nu-\ell.\)
Then
\[
(F\circ G)(x,u)
=
\left(
P(Q(x,u),V(u)),
U_{\mathbf A}(V(u))
\right).
\]
The second coordinate is
\[
U_{\mathbf A}(V(u))=m(nu-\ell)-k=mnu-(k+m\ell)=U_{\mathbf A\mathbf B}(u).
\]
Write
\(P_v(x)=P(x,v).\)
On the fundamental strip,
\(P(x,u)=\alpha(u)(x) \qquad(0\le u\le1).\)
If $u\in[2j,2j+1]$, then using $t^h(x,v)=(x,v-2h)$ with $h=-j$ gives
\(t^{-j}(x,u-2j)=(x,u).\)
Since $\rho(t^{-j})=1$, equivariance yields
\(P_u=\alpha(u-2j).\)
If $u\in[2j-1,2j]$ and $v=2j-u$, then
\(a(-x,v)=(x,-v), \qquad t^{-j}(x,-v)=(x,u),\)
so
\((x,u)=t^{-j}a(-x,v).\)
Since $\rho(t^{-j}a)=-1$,
\(P_u(x)=-\alpha(2j-u)(-x).\)
Define
\(\Theta(h)=A\circ h\circ A.\)
Then
\(P_u=\Theta(\alpha(2j-u)) \qquad(2j-1\le u\le2j).\)
Since $q_d$ is odd, $\Theta(q_d)=q_d$, and therefore
\(P_h=q_d \qquad(h\in\ZZ).\)

Let
\(x=[\alpha]\in\pi_1(B_d,q_d).\)
On every interval $[2j,2j+1]$, the path $u\mapsto P_u$ represents $x$. On every interval $[2j-1,2j]$, it is $\Theta$ applied to the reverse of $\alpha$; the reverse represents $-x$, while \cref{thm:reflection-conjugation} gives $\Theta_*=-\id$, so the interval again contributes $x$. Hence a positive unit interval in the $v$-coordinate always contributes $x$.

Moreover,
\(V(0)=-\ell, \qquad V(1)=n-\ell, \qquad V(1)-V(0)=\det(\mathbf B),\)
and $P_{-\ell}=P_{n-\ell}=q_d$. Therefore
\(\left[u\longmapsto P_{V(u)}\right] =\det(\mathbf B)x.\)

Define
\[
\mathcal R_e:B_d\longrightarrow B_{de},
\qquad
\mathcal R_e(h)=h\circ q_e,
\]
and
\[
\mathcal L_d:B_e\longrightarrow B_{de},
\qquad
\mathcal L_d(k)=q_d\circ k.
\]
Since $q_e(N)=N$, $\mathcal R_e$ restricts to
\(r_e:F_d\longrightarrow F_{de}\)
and covers $\id_{S^2}$. Thus
\[
\begin{tikzcd}[column sep=2.3em,row sep=2.2em]
F_d \arrow[r,"r_e"] \arrow[d,"i_d"']
  & F_{de} \arrow[d,"i_{de}"] \\
B_d \arrow[r,"\mathcal R_e"] \arrow[d,"\operatorname{ev}_N"']
  & B_{de} \arrow[d,"\operatorname{ev}_N"] \\
S^2 \arrow[r,"\id_{S^2}"'] & S^2 .
\end{tikzcd}
\]
Naturality of the long exact sequences gives
\[
\begin{tikzcd}[column sep=2.3em,row sep=2.2em]
\pi_2(S^2,N) \arrow[r,"\partial_d"] \arrow[d,"\id"']
  & \pi_1(F_d,q_d) \arrow[r,"(i_d)_*"] \arrow[d,"(r_e)_*"']
  & \pi_1(B_d,q_d) \arrow[r] \arrow[d,"(\mathcal R_e)_*"'] & 0 \\
\pi_2(S^2,N) \arrow[r,"\partial_{de}"']
  & \pi_1(F_{de},q_{de}) \arrow[r,"(i_{de})_*"']
  & \pi_1(B_{de},q_{de}) \arrow[r] & 0 .
\end{tikzcd}
\]
Here $r_e(q_d)=q_d\circ q_e=q_{de}$. Choose a generator
$\iota_2\in\pi_2(S^2)$ and generators $\xi_d,\xi_{de}$ such that
\[
\partial_d(\iota_2)=2d\,\xi_d,
\qquad
\partial_{de}(\iota_2)=2de\,\xi_{de}.
\]
If $(r_e)_*(\xi_d)=a\xi_{de}$, the left square gives
\(2da\,\xi_{de}=2de\,\xi_{de},\)
so $a=e$. Therefore
\((r_e)_*(\xi_d)=e\xi_{de}.\)
Writing
\[
\bar\xi_d=(i_d)_*(\xi_d),
\qquad
\bar\xi_{de}=(i_{de})_*(\xi_{de}),
\]
the middle square gives
\((\mathcal R_e)_*(\bar\xi_d)=e\bar\xi_{de},\)
and hence
\[
\boxed{(\mathcal R_e)_*([n])=[en].}
\]

Similarly, since $q_d(N)=N$, $\mathcal L_d$ restricts to
\(\ell_d:F_e\longrightarrow F_{de}\)
and covers $q_d:S^2\to S^2$. Thus
\[
\begin{tikzcd}[column sep=2.3em,row sep=2.2em]
F_e \arrow[r,"\ell_d"] \arrow[d,"i_e"']
  & F_{de} \arrow[d,"i_{de}"] \\
B_e \arrow[r,"\mathcal L_d"] \arrow[d,"\operatorname{ev}_N"']
  & B_{de} \arrow[d,"\operatorname{ev}_N"] \\
S^2 \arrow[r,"q_d"'] & S^2 .
\end{tikzcd}
\]
and
\[
\begin{tikzcd}[column sep=2.3em,row sep=2.2em]
\pi_2(S^2,N) \arrow[r,"\partial_e"] \arrow[d,"(q_d)_*"']
  & \pi_1(F_e,q_e) \arrow[r,"(i_e)_*"] \arrow[d,"(\ell_d)_*"']
  & \pi_1(B_e,q_e) \arrow[r] \arrow[d,"(\mathcal L_d)_*"'] & 0 \\
\pi_2(S^2,N) \arrow[r,"\partial_{de}"']
  & \pi_1(F_{de},q_{de}) \arrow[r,"(i_{de})_*"']
  & \pi_1(B_{de},q_{de}) \arrow[r] & 0 .
\end{tikzcd}
\]
Since $\deg q_d=d$,
\((q_d)_*(\iota_2)=d\iota_2.\)
Choose generators $\xi_e,\xi_{de}$ such that
\[
\partial_e(\iota_2)=2e\xi_e,
\qquad
\partial_{de}(\iota_2)=2de\xi_{de}.
\]
If $(\ell_d)_*(\xi_e)=b\xi_{de}$, the left square gives
\(2eb\xi_{de}=2d^2e\xi_{de},\)
so $b=d^2$. Hence
\((\ell_d)_*(\xi_e)=d^2\xi_{de},\)
and the middle square yields
\[
\boxed{(\mathcal L_d)_*([n])=[d^2n].}
\]

Let $y=[\beta]\in\pi_1(B_e,q_e)$. Define
\[
\operatorname{Comp}:B_d\times B_e\longrightarrow B_{de},
\qquad
\operatorname{Comp}(h,k)=h\circ k.
\]
Since $q_d\circ q_e=q_{de}$, this is a based map. Put
\(j_1(h)=(h,q_e), \qquad j_2(k)=(q_d,k).\)
Then
\[
\operatorname{Comp}\circ j_1=\mathcal R_e,
\qquad
\operatorname{Comp}\circ j_2=\mathcal L_d,
\]
and we have
\[
\begin{tikzcd}[column sep=2.3em,row sep=2.2em]
B_d \arrow[r,"j_1"] \arrow[d,"\mathcal R_e"']
  & B_d\times B_e \arrow[d,"\operatorname{Comp}"] \\
B_{de} \arrow[r,equal] & B_{de}
\end{tikzcd}
\qquad
\begin{tikzcd}[column sep=2.3em,row sep=2.2em]
B_e \arrow[r,"j_2"] \arrow[d,"\mathcal L_d"']
  & B_d\times B_e \arrow[d,"\operatorname{Comp}"] \\
B_{de} \arrow[r,equal] & B_{de}.
\end{tikzcd}
\]
Under
\[
\pi_1(B_d\times B_e,(q_d,q_e))
\cong
\pi_1(B_d,q_d)\times\pi_1(B_e,q_e),
\]
the loop
\(u\longmapsto(P_{V(u)},\beta(u))\)
represents $(\det(\mathbf B)x,y)$. Since
\((\det(\mathbf B)x,y)=(\det(\mathbf B)x,0)+(0,y)\)
and $\operatorname{Comp}_*$ is a homomorphism,
\[
\begin{aligned}
[H]
&=\operatorname{Comp}_*(\det(\mathbf B)x,y)\\
&=(\mathcal R_e)_*(\det(\mathbf B)x)+(\mathcal L_d)_*(y)\\
&=e\det(\mathbf B)x+d^2y
\end{aligned}
\]
in
\[
\pi_1(B_{de},q_{de})\cong\ZZ/(2|de|)\ZZ.
\]

Finally, by \cref{thm:rotation-loop-class},
\[
\pi_0\mathcal H(O_{de},O_{de};B_{de})
\cong
\frac{\ZZ/(2|de|)\ZZ}{\langle|de|\rangle}
\cong
\ZZ/|de|\ZZ.
\]
If $x$ maps to $c$ and $y$ maps to $c'$, then
\(c_{\mathrm{new}} = e\det(\mathbf B)c+d^2c' \pmod{|de|}.\)
This is independent of the chosen integer representatives: changing a representative of $c$ by $|d|$ changes the expression by $e\det(\mathbf B)|d|$, a multiple of $|de|$; changing a representative of $c'$ by $|e|$ changes it by
$d^2|e|=|d||de|$, again a multiple of $|de|$.

Finally,
\(\id_G(a)=a=t^0a, \qquad \id_G(b)=b=t^{-1}a,\)
so $\id_G=\varphi_{I_2}$. The identity map induces the identity on $\pi_2(M)$, and $\ZZ/|1|\ZZ=0$. Thus
\([\id_M]=(I_2,1,0).\)
Since the operation $\star$ is induced on pointed homotopy classes by composition of pointed maps, and composition of maps is associative, the operation $\star$ is associative as well. Equivalently, for any pointed homotopy classes $[f],[g],[h]\in [M,M]_{*}$, one has
\[
([f]\star[g])\star[h]
=
[(f\circ g)\circ h]
=
[f\circ(g\circ h)]
=
[f]\star([g]\star[h]).
\]
Hence $\mathcal R$ is a unital submonoid.
\end{proof}

\begin{corollary}[Coordinate monoid]\label{cor:coordinate-monoid}
Put
\[
\mathscr P=\{(\mathbf A,d,c)\mid\mathbf A\in\mathscr M,\ d\in2\mathbb Z+1,\ c\in\mathbb Z/|d|\mathbb Z\},
\]
and define $\diamond$ by the formula in \cref{cor:double-reflection-composition}. Then $(\mathscr P,\diamond)$ is a monoid with identity $(I_2,1,0)$, and the coordinate map
\[
\Psi:(\mathcal R,\star)\longrightarrow(\mathscr P,\diamond),\qquad[f]\longmapsto(\mathbf A,d,c)
\]
is an isomorphism of monoids.
\end{corollary}

\begin{proof}
By \cref{thm:classification} and its proof, for every
\(\mathbf A\in\mathscr M, \qquad d\in 2\mathbb Z+1,\)
after fixing
\[
f_{*,1}=\varphi_{\mathbf A},
\qquad
f_{*,2}=d\,\operatorname{id}_{\mathbb Z},
\]
the corresponding set
\([M,M]_*^{\varphi_{\mathbf A},d}\)
of pointed homotopy classes carries a free and transitive $\mathbb Z/|d|\mathbb Z$-action.
After fixing the reference maps, compatible generators, and quotient
isomorphisms chosen in \cref{cor:double-reflection-composition}, every
\([f]\in [M,M]_*^{\varphi_{\mathbf A},d}\)
is represented uniquely by coordinates
\((\mathbf A,d,c), \qquad c\in\mathbb Z/|d|\mathbb Z.\)
Hence the coordinate map
\(\Psi:\mathcal R\longrightarrow\mathscr P\)
is a bijection.

We next examine composition in these coordinates.
Suppose
\[\Psi([f])=(\mathbf A,d,c), \qquad \Psi([g])=(\mathbf B,e,c').\]
By \cref{cor:double-reflection-composition}, the first two coordinates of
\([f]\star[g]=[f\circ g]\)
are
\(\mathbf A\mathbf B, \qquad de,\)
and its third coordinate is
\(e\det(\mathbf B)c+d^{2}c' \pmod{|de|}.\)
More precisely, if
\(\widetilde c,\widetilde c'\in\mathbb Z\)
represent
\[
c\in\mathbb Z/|d|\mathbb Z,
\qquad
c'\in\mathbb Z/|e|\mathbb Z,
\]
respectively, then the third coordinate is
\(\left[ e\det(\mathbf B)\widetilde c+d^{2}\widetilde c' \right]_{|de|}.\)
\cref{cor:double-reflection-composition} proves that this residue class is independent of the
choices of \(\widetilde c\) and \(\widetilde c'\).
Thus the operation
\(\diamond:\mathscr P\times\mathscr P\longrightarrow\mathscr P\)
is well defined, and
\(\Psi([f]\star[g]) = \Psi([f])\diamond\Psi([g]).\)
Hence
\(\Psi\)
preserves multiplication.

Furthermore, \cref{cor:double-reflection-composition} gives
\([\operatorname{id}_{M}] = (I_2,1,0).\)
Therefore
\(\Psi([\operatorname{id}_{M}]) = (I_2,1,0),\)
so
\((I_2,1,0)\)
is the identity element of
\((\mathscr P,\diamond).\)

It remains to verify associativity.
Let
\(x,y,z\in\mathscr P.\)
Since
\(\Psi:\mathcal R\longrightarrow\mathscr P\)
is bijective, there exist unique
\([f],[g],[h]\in \mathcal R\)
such that
\(x=\Psi([f]), \qquad y=\Psi([g]), \qquad z=\Psi([h]).\)
Since
\(\star\)
is induced by composition of maps,
\(([f]\star[g])\star[h] = [f]\star([g]\star[h]).\)
Using
\[
\Psi([u]\star[v])
=
\Psi([u])\diamond\Psi([v]),
\]
we obtain
\((x\diamond y)\diamond z = x\diamond(y\diamond z).\)
Thus
\(\diamond\)
is associative.

Therefore
\((\mathscr P,\diamond)\)
is a monoid, and the bijection
\(\Psi:\mathcal R\longrightarrow\mathscr P\)
preserves both multiplication and the identity element. Hence
\(\Psi\)
is an isomorphism of monoids:
\((\mathcal R,\star)\cong(\mathscr P,\diamond).\)
\end{proof}

\subsection{Units and pointed self-homotopy equivalences}

\begin{corollary}[Units and pointed homotopy equivalences]\label{cor:pointed-homotopy-equivalences}
A pointed self-map $f$ is a pointed homotopy equivalence if and only if $f_{*,1}$ and $f_{*,2}$ are isomorphisms, equivalently,
\[
[f]=(\mathbf A,d,0),\qquad\det\mathbf A=\pm1,\quad d=\pm1.
\]
In that case,
\[
(\mathbf A,d,0)^{-1}=(\mathbf A^{-1},d,0),\qquad[M,M]_*^\times\cong D_\infty\times C_2.
\]
Here $[M,M]_*^\times$ is the group of units, $D_\infty=G$, and $C_2=\mathbb Z/2\mathbb Z$.
\end{corollary}

\begin{proof}
A pointed homotopy class is invertible in $[M,M]_*$ if and only if its representative is a pointed homotopy equivalence. The trivial fundamental-group part contains no units: if $f_{*,1}=\varphi_{\mathrm{triv}}$, then for every $g$,
\[
(f\circ g)_{*,1}
=f_{*,1}\circ g_{*,1}
=\varphi_{\mathrm{triv}}
\neq\id_G.
\]
Hence every unit belongs to the double-reflection part.

Let $[f]=(\mathbf A,d,c)$ have inverse $[g]=(\mathbf B,e,c')$. The composition formula gives
\[
\mathbf A\mathbf B=I_2,\qquad de=1.
\]
Thus $d=e=\pm1$. Both secondary coordinate groups are trivial, so $c=c'=0$. Taking determinants gives $\det\mathbf A\det\mathbf B=1$. The two factors are integers, hence $\det\mathbf A=\pm1$.

Conversely, write $\mathbf A=A(m,k)$ with $m^2=1$, and suppose $d^2=1$. Set $\mathbf B=A(m,-mk)$. Direct multiplication gives
\[
\mathbf A\mathbf B=A(m^2,k-m^2k)=I_2,\qquad
\mathbf B\mathbf A=A(m^2,-mk+mk)=I_2.
\]
By \cref{thm:classification}, the class $[g]=(\mathbf B,d,0)$ exists. In both orders the degree multiplier is $d^2=1$ and the secondary coordinate vanishes. Thus
\[
[f]\star[g]=(I_2,1,0)=[g]\star[f].
\]
This proves the inverse formula. The map $(\mathbf A,d,0)\mapsto(\mathbf A,d)$ is a bijection onto $\mathscr M^\times\times\{\pm1\}$. The composition formula makes it multiplicative for the componentwise product. By \cref{prop:matrix-semigroup} the first factor is isomorphic to $D_\infty$, and the second is isomorphic to $C_2$.
\end{proof}

\subsection{Composition properties of classes inducing the trivial endomorphism}

\begin{corollary}[Composition properties of the trivial-endomorphism part]\label{prop:trivial-composition}
The class $0$ is a two-sided zero of $[M,M]_*$. For every $[f]\in\mathcal T$ and $[g]\in[M,M]_*$, both $[f]\circ[g]$ and $[g]\circ[f]$ belong to $\mathcal T$, and
\[
\mathcal T\circ\mathcal T=\{0\}.
\]
In particular, $0$ is the unique idempotent in $\mathcal T$.
\end{corollary}

\begin{proof}
We first prove the two-sided absorption property under composition.
Let $[f]\in \mathcal T$ and let $[g]\in[M,M]_{*}$ be arbitrary.
Since $f_{*,1}=\varphi_{\mathrm{triv}}$, for every $h\in G$,
\((f_{*,1}\circ g_{*,1})(h)=1,\)
and hence
\((f\circ g)_{*,1}=\varphi_{\mathrm{triv}}.\)
On the other hand,
\((g_{*,1}\circ f_{*,1})(h)=g_{*,1}(1)=1,\)
so
\((g\circ f)_{*,1}=\varphi_{\mathrm{triv}}.\)
Therefore
\[
\mathcal T\circ[M,M]_{*}\subseteq \mathcal T,
\qquad
[M,M]_{*}\circ \mathcal T\subseteq \mathcal T.
\]

We next prove that $[c]$ is a zero element.
For any pointed self-map $f$, since $c$ is constant,
\(c\circ f=c.\)
Since $f(x_{0})=x_{0}$,
\((f\circ c)(x)=f(x_{0})=x_{0}=c(x)\)
for every $x\in M$, and therefore
\(f\circ c=c.\)
Hence
\([c]\circ[f]=[c], \qquad [f]\circ[c]=[c].\)

Finally, we prove that
\(\mathcal T\circ \mathcal T=\{[c]\}.\)
Take arbitrary
\([f],[g]\in \mathcal T.\)
By the correspondence established above between pointed lifts and twisted
equivariant maps, together with
\[
\pi_{0}(\mathcal F_{\varphi_{\mathrm{triv}}})
\cong
\pi_{0}(\mathcal E_{\varphi_{\mathrm{triv}}}),
\]
we may choose standard representatives in the corresponding unbased
equivariant path components,
\[
F=(P,0),
\qquad
G=(Q,0):
S^{2}\times\mathbb R\longrightarrow S^{2}\times\mathbb R.
\]
Here the second coordinate is the constant $0$, which is the standard choice
in the trivial fundamental-group-homomorphism case of Lemma 5.1.

Strictly speaking, $F$ and $G$ need not fix the chosen lifted basepoint.
Choose corresponding based equivariant representatives $F^{b},G^{b}$ and
equivariant paths
\[
F_{t}:F^{b}\rightsquigarrow F,
\qquad
G_{t}:G^{b}\rightsquigarrow G.
\]
Then
\(K_{t}:=F_{t}\circ G_{t}\)
is a $\varphi_{\mathrm{triv}}$-equivariant path, since
\(K_{t}(hz) = F_{t}(G_{t}(hz)) = F_{t}(G_{t}(z)) = K_{t}(z).\)
Thus it suffices to determine the equivariant path component containing
$F\circ G$.

Since $\varphi_{\mathrm{triv}}(h)=1$,
\(\lambda_{\varphi_{\mathrm{triv}}}(h) = \rho(1) = 1,\)
and therefore both $P$ and $Q$ are $G$-invariant:
\(P(hz)=P(z), \qquad Q(hz)=Q(z).\)
For every $(x,u)\in S^{2}\times\mathbb R$,
\(G(x,u)=(Q(x,u),0),\)
and hence
\((F\circ G)(x,u) = \bigl(P(Q(x,u),0),0\bigr).\)
Define
\(P_{0}:S^{2}\longrightarrow S^{2}, \qquad P_{0}(y)=P(y,0).\)
Then
\((F\circ G)(x,u) = \bigl(P_{0}(Q(x,u)),0\bigr).\)

We now prove that $P_{0}$ is even.
The deck transformation $a$ acts by
\(a(y,v)=(-y,-v).\)
Since $P$ is $G$-invariant,
\(P(a(y,0))=P(y,0).\)
Since $a(y,0)=(-y,0)$, it follows that
\(P_{0}(-y)=P_{0}(y).\)
Thus $P_{0}\in E$. By \cref{lem:even-degree-zero}, every even map has ordinary degree zero,
so
\(\deg P_{0}=0.\)
Ordinary homotopy classes of self-maps of the sphere are classified by degree.
Hence there exist $y_{0}\in S^{2}$ and a homotopy
\(K:S^{2}\times I\longrightarrow S^{2}\)
such that
\(K(y,0)=P_{0}(y), \qquad K(y,1)=y_{0}.\)

Define
\[
H:(S^{2}\times\mathbb R)\times I
\longrightarrow
S^{2}\times\mathbb R
\]
by
\(H((x,u),t) = \bigl(K(Q(x,u),t),0\bigr).\)
At $t=0$,
\(H((x,u),0) = \bigl(P_{0}(Q(x,u)),0\bigr) = (F\circ G)(x,u),\)
whereas at $t=1$,
\(H((x,u),1)=(y_{0},0),\)
which is a constant map.

Since $Q$ is $G$-invariant, for every $h\in G$,
\[
\begin{aligned}
H(h(x,u),t)
&=
\bigl(K(Q(h(x,u)),t),0\bigr)\\
&=
\bigl(K(Q(x,u),t),0\bigr)\\
&=
H((x,u),t).
\end{aligned}
\]
The $\varphi_{\mathrm{triv}}$-twisted action on the target is trivial.
Therefore $H$ is a $\varphi_{\mathrm{triv}}$-equivariant homotopy.
Hence $F\circ G$ and a constant map lie in the same path component of
$\mathcal E_{\varphi_{\mathrm{triv}}}$.

Since $S^{2}\times\mathbb R$ is path connected, and every constant map is
equivariant for the trivial twist, all constant maps lie in the same path
component of $\mathcal E_{\varphi_{\mathrm{triv}}}$. In particular, $F\circ G$ lies in
the same path component as the constant map with value the chosen lifted
basepoint $\widetilde{x}_{0}$. The latter fixes the lifted basepoint and
descends to the pointed constant map $c$ on $M$. By
\[
\pi_{0}(\mathcal F_{\varphi_{\mathrm{triv}}})
\cong
\pi_{0}(\mathcal E_{\varphi_{\mathrm{triv}}}),
\]
we obtain
\([f\circ g]=[c].\)
Since $[f],[g]\in \mathcal T$ were arbitrary,
\(\mathcal T\circ \mathcal T\subseteq\{[c]\}.\)
Conversely, $[c]\in \mathcal T$ and
\([c]\circ[c]=[c],\)
so
\([c]\in \mathcal T\circ \mathcal T.\)
Therefore
\(\mathcal T\circ \mathcal T=\{[c]\}.\)

Finally, if $x\in \mathcal T$ and $x\circ x=x$, then
$\mathcal T\circ \mathcal T=\{[c]\}$ gives $x\circ x=[c]$, and hence
$x=[c]$. Thus $[c]$ is the unique idempotent element of $\mathcal T$.
\end{proof}

\subsection{Secondary invariants and induced homotopy homomorphisms}

\begin{corollary}[Failure of detection by homotopy groups]\label{cor:homotopy-groups-do-not-detect}
For every $\mathbf A\in\mathscr M$ and odd integer $d$, the set $[M,M]_*^{\varphi_{\mathbf A},d}$ has exactly $|d|$ elements, and any two classes $[f],[g]$ in it satisfy
\[
f_{*,n}=g_{*,n}:\pi_n(M,x_0)\longrightarrow\pi_n(M,x_0)\qquad(n\ge1).
\]
\end{corollary}

\begin{proof}
By \cref{thm:classification} and its proof, $[M,M]_*^{\varphi_{\mathbf A},d}$ is nonempty and carries a free and transitive $\mathbb Z/|d|\mathbb Z$-action. Choosing a reference homotopy class gives a bijection with this cyclic group, so its cardinality is $|d|$ for every odd $d$, including $d=\pm1$.
Let $[f],[g]\in[M,M]_*^{\varphi_{\mathbf A},d}$ be arbitrary; they need not be distinct. By definition,

\[
f_{*,1}=g_{*,1}=\varphi_{\mathbf A},\qquad
f_{*,2}=g_{*,2}=d\,\mathrm{id}_{\mathbb Z}.
\]
It remains to compare the induced maps for $n\ge2$.
Fix $\widetilde x_0=(x_b,u_b)$ over $x_0$. Let $F,G$ be the unique pointed lifts of $fp,gp$, respectively. They satisfy the strictly commutative diagrams

\[
\begin{tikzcd}[column sep=large,row sep=large]
(\widetilde M,\widetilde x_0)\arrow[r,"F"]\arrow[d,"p"'] & (\widetilde M,\widetilde x_0)\arrow[d,"p"]\\
(M,x_0)\arrow[r,"f"'] & (M,x_0)
\end{tikzcd}
\]
\[
\begin{tikzcd}[column sep=large,row sep=large]
(\widetilde M,\widetilde x_0)\arrow[r,"G"]\arrow[d,"p"'] & (\widetilde M,\widetilde x_0)\arrow[d,"p"]\\
(M,x_0)\arrow[r,"g"'] & (M,x_0).
\end{tikzcd}
\]
Define pointed maps $q(x,u)=x$ and $i(x)=(x,u_b)$. Then $qi=\mathrm{id}_{S^2}$. The homotopy

\[
H((x,u),t)=(x,(1-t)u+tu_b)
\]
satisfies $H(-,0)=\mathrm{id}_{\widetilde M}$, $H(-,1)=iq$, and $H(\widetilde x_0,t)=\widetilde x_0$ for all $t$. Thus $i,q$ are pointed homotopy inverses. Consequently, for every $n\ge2$,

\[
q_*i_*=\mathrm{id}_{\pi_n(S^2,x_b)},\qquad
i_*q_*=\mathrm{id}_{\pi_n(\widetilde M,\widetilde x_0)}.
\]
Set $h_f=qFi$ and $h_g=qGi$. These sphere maps are pointed because $i,F,G,q$ are pointed. The cover induces an isomorphism $p_*$ in every degree $n\ge2$. For $n=2$, the first square above gives

\[
F_*=p_*^{-1}f_{*,2}p_*.
\]
\[
\begin{aligned}
(h_f)_*&=q_*F_*i_*\\
&=q_*p_*^{-1}f_{*,2}p_*i_*\\
&=q_*p_*^{-1}(d\,\mathrm{id})p_*i_*\\
&=d\,q_*i_*=d\,\mathrm{id}_{\pi_2(S^2,x_b)}.
\end{aligned}
\]
Here the integer identifications are exactly those fixed in the classification theorem. Replacing $F,f$ by $G,g$ gives $(h_g)_*=d\,\mathrm{id}$ as well. The degree agrees with this multiplier on $\pi_2(S^2)$ by the Hurewicz isomorphism, so $\deg h_f=\deg h_g=d$.
The pointed degree classification $[S^2,S^2]_*\cong\mathbb Z$ therefore gives $h_f\simeq_*h_g$. Homotopy invariance now implies, in every degree $n\ge2$,

\[
q_*F_*i_*=(h_f)_*=(h_g)_*=q_*G_*i_*.
\]
The relationship with the chosen sphere model can also be displayed as follows; the top square commutes on homotopy groups because $i_*q_*=\mathrm{id}$, and the bottom square comes from $pF=fp$:

\[
\begin{tikzcd}[column sep=large,row sep=large]
\pi_n(S^2,x_b)\arrow[r,"(h_f)_*"]\arrow[d,"i_*"',"\cong"] & \pi_n(S^2,x_b)\arrow[d,"i_*","\cong"']\\
\pi_n(\widetilde M,\widetilde x_0)\arrow[r,"F_*"]\arrow[d,"p_*"',"\cong"] & \pi_n(\widetilde M,\widetilde x_0)\arrow[d,"p_*","\cong"']\\
\pi_n(M,x_0)\arrow[r,"f_{*,n}"'] & \pi_n(M,x_0).
\end{tikzcd}
\]
The same diagram holds with $h_g,G,g$. Multiplying the equality $q_*F_*i_*=q_*G_*i_*$ on the left by $i_*$ and on the right by $q_*$ gives, explicitly,

\[
F_*=(i_*q_*)F_*(i_*q_*)
=i_*(q_*F_*i_*)q_*
=i_*(q_*G_*i_*)q_*=G_*.
\]
Finally the covering squares give

\[
f_{*,n}=p_*F_*p_*^{-1}=p_*G_*p_*^{-1}=g_{*,n}\qquad(n\ge2).
\]
Together with the equality for $n=1$, this proves the assertion for the arbitrary elements originally chosen. Distinctness of $[f]$ and $[g]$ was never required; when $|d|=1$, the classification set has one element and the assertion still holds.
\end{proof}

\begin{remark}
In particular, if $|d|>1$, these are $|d|$ distinct pointed homotopy classes with the same induced homomorphisms in every degree. Hence the map
\[
\mathcal D:[M,M]_*\longrightarrow\prod_{n\ge1}\operatorname{End}(\pi_n(M,x_0)),\qquad
\mathcal D([f])=(f_{*,n})_{n\ge1},
\]
is not injective.
\end{remark}

\begin{remark}[Secondary invariant]\label{stmt:original-54}
The coordinate $c\in\ZZ/|d|\ZZ$ distinguishes pointed homotopy classes within $[M,M]_*^{\varphi_{\mathbf A},d}$, but is undetected by every induced homomorphism $f_{*,n}$. This assertion concerns the coordinate after the choices in \cref{def:secondary-coordinates}; the underlying homotopy-class set described in \cref{thm:classification} is independent of these choices.
\end{remark}

\section{Free homotopy classification of self-maps of \texorpdfstring{\(\RP^3\#\RP^3\)}{RP3 \# RP3}}
\label{sec:free-homotopy}

Using the pointed classification and composition formula established above, we compute the action of the target fundamental group on the pointed homotopy classes of self-maps, thereby obtaining the free homotopy classification and determining the group of self-homotopy equivalences up to free homotopy.

\subsection{Free homotopy classes and the fundamental group action}

\begin{definition}[Forgetful map]\label{def:free-coordinates}
Retain the basepoint $x_0\in M$ and let $[M,M]$ denote the set of free homotopy classes of self-maps. To distinguish the two kinds of classes, in this section write $[f]_*$ for the pointed homotopy class (denoted $[f]$ in the preceding sections) and $[f]$ for the free homotopy class.
Define the forgetful map
\[
Q:[M,M]_*\longrightarrow[M,M],\qquad Q([f]_*)=[f].
\]
\end{definition}

We retain the notation $G,t,a,b,\mathcal T,\mathcal R$ and the classification coordinates of \cref{thm:classification,def:secondary-coordinates}.
Put $I=[0,1]$, write $p:\widetilde M=S^2\times\RR\to M$ for the universal covering, and choose the lifted basepoint $\widetilde x_0$.
By \cref{def:cover-action}, for every $r\in\ZZ$ we have
\begin{equation}\label{eq:free-deck}
 t^r(x,v)=(x,v-2r),\qquad t^ra(x,v)=(-x,-v-2r).
\end{equation}
Put $\rho(t^ra^\delta)=(-1)^\delta$, where $\delta\in\{0,1\}$.

\begin{lemma}[Orbit description of free homotopy classes]\label{lem:free-orbit}
Let $(X,x_0)$ be well-pointed, meaning that the inclusion $\{x_0\}\hookrightarrow X$ is a closed cofibration,
and let $(Y,y_0)$ be a path-connected pointed space. The fundamental group $\pi_1(Y,y_0)$ of the target has
a natural left action on the set $[X,Y]_*$ of pointed homotopy classes,
\[
\pi_1(Y,y_0)\times[X,Y]_*
\longrightarrow[X,Y]_*,\qquad (g,[f]_*)\longmapsto g\cdot[f]_*.
\]
Explicitly, for any $g\in\pi_1(Y,y_0)$ and $[f]_*\in[X,Y]_*$,
choose a based loop representing $g^{-1}$,
$\lambda:(I,\partial I)\longrightarrow(Y,y_0)$
and a homotopy $H:X\times I\to Y$ satisfying
\[
H(x,0)=f(x),\qquad H(x_0,s)=\lambda(s).
\]
Define
$g\cdot[f]_*:=[H(-,1)]_*.$
The action is well defined, and the forgetful map
\[
Q:[X,Y]_*\longrightarrow[X,Y],\qquad [f]_*\longmapsto[f]
\]
is surjective. For any $[f_0]_*,[f_1]_*\in[X,Y]_*$, we have
\[
Q([f_0]_*)=Q([f_1]_*)\quad\Longleftrightarrow\quad
\exists\,g\in\pi_1(Y,y_0),\ [f_1]_* = g\cdot[f_0]_*.
\]
Consequently, writing $[X,Y]_*/\pi_1(Y,y_0)$ for the orbit set of this left action, $Q$ induces a bijection of sets
\[
[X,Y]_*/\pi_1(Y,y_0)
\xrightarrow{\ \cong\ }[X,Y].
\]
\end{lemma}

\begin{proof}
See Hatcher~\cite[\S4.A, Propositions~4A.1--4A.2, pp.~421--422]{hatcher}.
\end{proof}

In particular, $M$ is a path-connected CW complex and admits a CW structure with $x_0$ as a $0$-cell.
Taking $X=Y=M$ gives the left action of $G=\pi_1(M,x_0)$ on $[M,M]_*$ and a bijection of sets
\begin{equation}\label{eq:free-orbit}
[M,M]_*/G\xrightarrow{\ \cong\ }[M,M].
\end{equation}
We use this left-action convention throughout, and write $S/G$ for the orbit set of a left action of $G$ on a set $S$.

\begin{lemma}[The change-of-basepoint action on the universal cover]\label{lem:free-deck-action}
Let $f:(M,x_0)\to(M,x_0)$ be a pointed map, and put
$\varphi=f_{*,1}:G\longrightarrow G,$
and let
$F_f:(\widetilde M,\widetilde x_0)\longrightarrow(\widetilde M,\widetilde x_0)$
be the unique pointed lift of $f\circ p$.
For any $g\in G$, write
\[
\operatorname{Ad}_g:G\longrightarrow G,\qquad \operatorname{Ad}_g(h)=ghg^{-1}.
\]
The following statements hold.
\begin{enumerate}
\item In the equivariant classification given by \cref{lem:equivariance-pi1,thm:remove-lift-basepoint}, the pointed homotopy class
\[
g\cdot[f]_*
\]
is represented by the $(\operatorname{Ad}_g\circ\varphi)$-equivariant homotopy class of the map $g\circ F_f:\widetilde M\longrightarrow\widetilde M$.
In particular, the action of the target fundamental group changes the induced fundamental-group homomorphism $\varphi$ to
$\operatorname{Ad}_g\circ\varphi,$
and changes the multiplier $d$ of the induced homomorphism on the second homotopy group to
$\rho(g)d.$
\item Put
\[
\beta_g=g\cdot[\id_M]_*.
\]
Retaining the composition convention in \cref{def:secondary-coordinates},
$[u]_*\star[v]_*=[u\circ v]_*,$
we have
\begin{equation}\label{eq:free-action-composition}
 g\cdot[f]_* = \beta_g\star[f]_*.
\end{equation}
\end{enumerate}
\end{lemma}
\begin{proof}
We first prove (1).
For any group endomorphism $\psi:G\to G$, following Definition~\ref{stmt:original-7}, write
\[
\mathcal E_\psi=\{F:\widetilde M\to\widetilde M\mid F(hz)=\psi(h)F(z),\ h\in G,\ z\in\widetilde M\}.
\]
Write its pointed subspace as
$\mathcal F_\psi=\{F\in\mathcal E_\psi\mid F(\widetilde x_0)=\widetilde x_0\}.$
By Lemma~\ref{lem:equivariance-pi1}, lifting pointed maps gives a bijection
$[M,M]_*^\psi\cong\pi_0(\mathcal F_\psi).$
By Theorem~\ref{thm:remove-lift-basepoint}, the inclusion
$\mathcal F_\psi\hookrightarrow\mathcal E_\psi$
induces a bijection on path components
$\pi_0(\mathcal F_\psi)\xrightarrow{\cong}\pi_0(\mathcal E_\psi).$
Composing these two bijections gives
\begin{equation}\label{eq:free-equivariant-components}
[M,M]_*^\psi\cong\pi_0(\mathcal E_\psi).
\end{equation}
Thus, when computing pointed homotopy classes, equivariant maps that do not fix the lifted basepoint may be used as representatives of the corresponding path components.
Now let $g\in G$.
First we verify the equivariance of $g\circ F_f$.
Since $F_f$ is $\varphi$-equivariant, for any $h\in G$ and $z\in\widetilde M$ we have
$F_f(hz)=\varphi(h)F_f(z).$
Therefore
\[
\begin{aligned}
(g\circ F_f)(hz)&=g\bigl(F_f(hz)\bigr)\\
&=g\bigl(\varphi(h)F_f(z)\bigr)\\
&=(g\varphi(h))F_f(z)\\
&=(g\varphi(h)g^{-1})\bigl(gF_f(z)\bigr)\\
&=(\operatorname{Ad}_g\circ\varphi)(h)(g\circ F_f)(z).
\end{aligned}
\]
Hence
\begin{equation}\label{eq:free-deck-equivariance}
g\circ F_f\in\mathcal E_{\operatorname{Ad}_g\circ\varphi}.
\end{equation}
However, since
$(g\circ F_f)(\widetilde x_0)=g\widetilde x_0,$
when $g\ne1_G$, $g\circ F_f$ is generally not a pointed lift.
We now prove that its equivariant path component corresponds precisely to $g\cdot[f]_*$.
Choose a based loop representing $g^{-1}$,
$\lambda:I\longrightarrow M.$
By the definition of the action in Lemma~\ref{lem:free-orbit}, choose a homotopy
$H:M\times I\longrightarrow M$
satisfying
\[
H(x,0)=f(x),\qquad H(x_0,s)=\lambda(s).
\]
Write
$H_s(x)=H(x,s).$
Then
\begin{equation}\label{eq:free-action-endpoint}
g\cdot[f]_*=[H_1]_*.
\end{equation}
Consider the continuous map
\[
\overline H:\widetilde M\times I\longrightarrow M,
\]
\[
\overline H(z,s)=H(p(z),s).
\]
At $s=0$,
\[
\begin{aligned}
\overline H(z,0)&=H(p(z),0)\\
&=f(p(z))\\
&=p(F_f(z)).
\end{aligned}
\]
Thus $F_f$ is a lift of $\overline H(-,0)$.
By the homotopy lifting property for covering maps~\cite[\S1.3, Proposition~1.30, p.~60]{hatcher}, there exists a unique continuous map
$\widetilde H:\widetilde M\times I\longrightarrow\widetilde M$
satisfying
\begin{equation}\label{eq:free-homotopy-initial}
\widetilde H(z,0)=F_f(z),
\end{equation}
and
\begin{equation}\label{eq:free-homotopy-projection}
p\bigl(\widetilde H(z,s)\bigr)=H(p(z),s).
\end{equation}
Write
$\widetilde H_s(z)=\widetilde H(z,s).$
Then for all $s\in I$,
$p\circ\widetilde H_s=H_s\circ p.$
We now verify that the homotopy $\widetilde H$ remains $\varphi$-equivariant throughout.
Fix $h\in G$ and define continuous maps
$U_h,V_h:\widetilde M\times I\longrightarrow\widetilde M$
by
\[
U_h(z,s)=\widetilde H(hz,s),
\]
\[
V_h(z,s)=\varphi(h)\widetilde H(z,s).
\]
Since $h$ is a deck transformation,
$p(hz)=p(z).$
Therefore
\[
\begin{aligned}
pU_h(z,s)&=p\bigl(\widetilde H(hz,s)\bigr)\\
&=H(p(hz),s)\\
&=H(p(z),s)\\
&=\overline H(z,s).
\end{aligned}
\]
On the other hand, since $\varphi(h)$ is also a deck transformation,
\[
\begin{aligned}
pV_h(z,s)&=p\bigl(\varphi(h)\widetilde H(z,s)\bigr)\\
&=p\bigl(\widetilde H(z,s)\bigr)\\
&=H(p(z),s)\\
&=\overline H(z,s).
\end{aligned}
\]
Hence both $U_h$ and $V_h$ are lifts of $\overline H$.
At $s=0$, by \eqref{eq:free-homotopy-initial} and the $\varphi$-equivariance of $F_f$,
\[
\begin{aligned}
U_h(z,0)&=\widetilde H(hz,0)\\
&=F_f(hz)\\
&=\varphi(h)F_f(z)\\
&=\varphi(h)\widetilde H(z,0)\\
&=V_h(z,0).
\end{aligned}
\]
Since $\widetilde M\times I$ is connected and the two lifts agree on $\widetilde M\times\{0\}$, uniqueness of lifts gives
$U_h=V_h.$
Since $h\in G$ is arbitrary,
\begin{equation}\label{eq:free-homotopy-equivariance}
\widetilde H_s(hz)=\varphi(h)\widetilde H_s(z)
\end{equation}
holds for every $h\in G$, $z\in\widetilde M$ and $s\in I$.
Next we compute the value of $\widetilde H_1$ at the lifted basepoint.
Consider the path
\[
\widetilde\lambda:I\longrightarrow\widetilde M,\qquad\widetilde\lambda(s)=\widetilde H(\widetilde x_0,s).
\]
By \eqref{eq:free-homotopy-projection},
\[
\begin{aligned}
p(\widetilde\lambda(s))&=p\bigl(\widetilde H(\widetilde x_0,s)\bigr)\\
&=H(p(\widetilde x_0),s)\\
&=H(x_0,s)\\
&=\lambda(s).
\end{aligned}
\]
Moreover,
\[
\widetilde\lambda(0)=\widetilde H(\widetilde x_0,0)=F_f(\widetilde x_0)=\widetilde x_0.
\]
Thus $\widetilde\lambda$ is the unique lift of the loop $\lambda$ starting at $\widetilde x_0$.
Retaining the convention in \eqref{eq:path-end-deck} relating the fundamental group and the deck transformation group, since
$[\lambda]=g^{-1},$
we have
$\widetilde\lambda(1)=g^{-1}\widetilde x_0.$
Consequently,
\begin{equation}\label{eq:free-lift-endpoint}
\widetilde H_1(\widetilde x_0)=g^{-1}\widetilde x_0.
\end{equation}
Define
$K_s=g\circ\widetilde H_s.$
We now verify that $K_1$ is the pointed lift of $H_1\circ p$.
Since $g$ is a deck transformation,
$p\circ g=p.$
Therefore
\[
\begin{aligned}
p\circ K_1&=p\circ g\circ\widetilde H_1\\
&=p\circ\widetilde H_1\\
&=H_1\circ p.
\end{aligned}
\]
Also, by \eqref{eq:free-lift-endpoint},
\[
\begin{aligned}
K_1(\widetilde x_0)&=g\widetilde H_1(\widetilde x_0)\\
&=g(g^{-1}\widetilde x_0)\\
&=\widetilde x_0.
\end{aligned}
\]
Hence, by uniqueness of the pointed lift,
\begin{equation}\label{eq:free-pointed-lift-endpoint}
K_1=F_{H_1}.
\end{equation}
We now verify that $K_s$ is $(\operatorname{Ad}_g\circ\varphi)$-equivariant.
By \eqref{eq:free-homotopy-equivariance}, for any $h\in G$ and $z\in\widetilde M$,
\[
\begin{aligned}
K_s(hz)&=g\widetilde H_s(hz)\\
&=g\varphi(h)\widetilde H_s(z)\\
&=(g\varphi(h)g^{-1})g\widetilde H_s(z)\\
&=(\operatorname{Ad}_g\circ\varphi)(h)K_s(z).
\end{aligned}
\]
Thus for every $s\in I$,
$K_s\in\mathcal E_{\operatorname{Ad}_g\circ\varphi}.$
At the two endpoints,
\[
K_0=g\circ F_f,
\]
\[
K_1=F_{H_1}.
\]
Hence $s\mapsto K_s$ gives
$g\circ F_f\simeq_{\operatorname{Ad}_g\circ\varphi}F_{H_1}.$
That is,
$[g\circ F_f]=[F_{H_1}]$
in
$\pi_0(\mathcal E_{\operatorname{Ad}_g\circ\varphi}).$
By \eqref{eq:free-equivariant-components} and \eqref{eq:free-action-endpoint}, the equivariant path component containing $g\circ F_f$ corresponds to the pointed homotopy class
$[H_1]_* = g\cdot[f]_*.$
Moreover, since $F_{H_1}=K_1$ is an $(\operatorname{Ad}_g\circ\varphi)$-equivariant map fixing the lifted basepoint, Lemma~\ref{lem:equivariance-pi1} gives
$(H_1)_{*,1}=\operatorname{Ad}_g\circ\varphi.$
This proves the transformation formula for the induced homomorphism on the fundamental group.
We now compute the multiplier on the second homotopy group.
Write $F_f$ as
$F_f(x,v)=(P_f(x,v),U_f(x,v)).$
Fix $v_0\in(0,1)$ and define the sphere slice map
\[
(P_f)_{v_0}:S^2\longrightarrow S^2,\qquad (P_f)_{v_0}(x)=P_f(x,v_0).
\]
By the interior slice computation in the proof of Theorem~\ref{thm:classification},
$\deg((P_f)_{v_0})=d,$
where
$f_{*,2}=d\id_{\ZZ}.$
By \eqref{eq:free-deck}, the sphere-coordinate component of $g\circ F_f$ is
$(x,v)\longmapsto\rho(g)P_f(x,v).$
Thus its sphere slice at $v_0$ is
\[
S^2\longrightarrow S^2,\qquad x\longmapsto\rho(g)(P_f)_{v_0}(x).
\]
If $\rho(g)=1$, this map is $(P_f)_{v_0}$ itself and its degree is $d$.
If $\rho(g)=-1$, let
\[
\mathcal A:S^2\longrightarrow S^2,\qquad\mathcal A(x)=-x
\]
be the antipodal map. Then the slice map is
$\mathcal A\circ(P_f)_{v_0}.$
Since
$\deg\mathcal A=(-1)^{2+1}=-1,$
we have
\[
\begin{aligned}
\deg(\mathcal A\circ(P_f)_{v_0})&=\deg\mathcal A\cdot\deg((P_f)_{v_0})\\
&=-d.
\end{aligned}
\]
In summary, the degree of the sphere slice of $g\circ F_f$ is
$\rho(g)d.$
Since $K_s$ is an equivariant homotopy from $g\circ F_f$ to $F_{H_1}$, restricting it to
$S^2\times\{v_0\}$
and projecting onto the sphere coordinate gives an ordinary homotopy between the corresponding sphere slice maps.
By homotopy invariance of degree, the degree of the sphere slice of $F_{H_1}$ also equals
$\rho(g)d.$
Applying the interior slice computation of Theorem~\ref{thm:classification} again gives
$(H_1)_{*,2}=\rho(g)d\id_{\ZZ}.$
Statement (1) follows.
We now prove (2).
Choose a loop representing $g^{-1}$,
$\lambda:I\longrightarrow M.$
By the construction in Lemma~\ref{lem:free-orbit}, there exists a homotopy
$J:M\times I\longrightarrow M$
satisfying
\[
J(x,0)=x,\qquad J(x_0,s)=\lambda(s).
\]
Write
$J_s(x)=J(x,s).$
Then
$J_0=\id_M.$
By the definition of the target fundamental group action,
$[J_1]_* = g\cdot[\id_M]_* = \beta_g.$
Define a homotopy
$L:M\times I\longrightarrow M$
by
$L(x,s)=J(f(x),s).$
This map is continuous because it is the composite of continuous maps
$M\times I\xrightarrow{f\times\id_I}M\times I\xrightarrow{J}M$
.
At $s=0$,
$L(x,0)=J(f(x),0)=f(x).$
At $s=1$,
\[
L(x,1)=J_1(f(x))=(J_1\circ f)(x).
\]
Since $f(x_0)=x_0$, we also have
\[
\begin{aligned}
L(x_0,s)&=J(f(x_0),s)\\
&=J(x_0,s)\\
&=\lambda(s).
\end{aligned}
\]
Thus $L$ is a free homotopy starting at $f$ whose basepoint track represents $g^{-1}$ and whose endpoint is the pointed map $J_1\circ f$.
By the definition of the action,
$g\cdot[f]_*=[J_1\circ f]_*.$
By the definition of composition of pointed homotopy classes,
$[J_1\circ f]_*=[J_1]_*\star[f]_*.$
Substituting $[J_1]_*=\beta_g$ gives
$g\cdot[f]_* = \beta_g\star[f]_*.$
Statement (2) follows.
\end{proof}

\subsection{Free homotopy classification of the trivial fundamental-group type}
We now study the action of the target fundamental group $G$ on $\mathcal T$.
By Lemma~\ref{lem:free-deck-action}, if $f_{*,1}=\varphi_{\mathrm{triv}}$, then for any $g\in G$,
$\operatorname{Ad}_g\circ\varphi_{\mathrm{triv}}=\varphi_{\mathrm{triv}}.$
Indeed, for any $h\in G$,
\[
\begin{aligned}
(\operatorname{Ad}_g\circ\varphi_{\mathrm{triv}})(h)&=g\varphi_{\mathrm{triv}}(h)g^{-1}\\
&=g1_Gg^{-1}\\
&=1_G\\
&=\varphi_{\mathrm{triv}}(h).
\end{aligned}
\]
Thus $\mathcal T$ is invariant under the action of $G$.
By Theorem~\ref{thm:classification}(1), each element of $\mathcal T$ has coordinates
\[
(\alpha,n),\qquad \alpha\in H^2(M;\ZZ),\qquad n\in\ZZ.
\]
To complete the free homotopy classification, we must determine the simultaneous action of $G$ on these two coordinates.
We first prove a lemma concerning sphere-valued maps and postcomposition with the antipodal map.
\begin{lemma}[A sphere-valued map and its antipodal postcomposition]\label{lem:free-antipodal}
Let $X$ be a CW complex, let $h:X\to S^2$ be continuous, and let
\[
\mathcal A:S^2\longrightarrow S^2,\qquad \mathcal A(y)=-y
\]
be the antipodal map.
Fix a generator $u$ of $H^2(S^2;\ZZ)$. Then
\begin{equation}\label{eq:free-antipodal-criterion}
h\simeq\mathcal A\circ h\quad\Longleftrightarrow\quad 2h^*u=0\quad\text{in }H^2(X;\ZZ).
\end{equation}
In particular, for
$M=\RP^3\#\RP^3,$
every continuous map $h:M\to S^2$ satisfies
$h\simeq\mathcal A\circ h.$
\end{lemma}
\begin{proof}
First we prove necessity.
Suppose
$h\simeq\mathcal A\circ h.$
By homotopy invariance of cohomology,
\[
h^*=(\mathcal A\circ h)^*:H^2(S^2;\ZZ)\longrightarrow H^2(X;\ZZ).
\]
Since the degree of the antipodal map $\mathcal A:S^2\to S^2$ is
$\deg\mathcal A=-1,$
we have
$\mathcal A^*u=-u.$
By contravariant functoriality of cohomology,
$(\mathcal A\circ h)^*=h^*\circ\mathcal A^*.$
Therefore
\[
\begin{aligned}
h^*u&=(\mathcal A\circ h)^*u\\
&=h^*(\mathcal A^*u)\\
&=h^*(-u)\\
&=-h^*u.
\end{aligned}
\]
Adding $h^*u$ to both sides gives
$2h^*u=0.$
Necessity follows.
We now prove sufficiency.
Suppose
$2h^*u=0.$
Consider the tangent bundle of the sphere
$\pi:TS^2\longrightarrow S^2$
and its pullback along $h$,
$h^*TS^2\longrightarrow X.$
Regard $S^2$ as the unit sphere in $\RR^3$, that is,
$S^2=\{y\in\RR^3\mid\|y\|=1\}.$
For any $y\in S^2$, we have
$T_yS^2=\{v\in\RR^3\mid\langle y,v\rangle=0\}.$
On each tangent plane define a linear map
$J_y:T_yS^2\longrightarrow T_yS^2$
by
$J_y(v)=y\times v,$
where $\times$ denotes the cross product in $\RR^3$.
Since the cross product $y\times v$ is orthogonal to $y$,
$\langle y,y\times v\rangle=0.$
Hence
$J_y(v)\in T_yS^2,$
so $J_y$ is well defined.
By the vector triple product formula,
\[
y\times(y\times v)=\langle y,v\rangle y-\langle y,y\rangle v.
\]
For $v\in T_yS^2$, we have
$\langle y,v\rangle=0,\qquad\langle y,y\rangle=1.$
Therefore
\[
\begin{aligned}
J_y^2(v)&=y\times(y\times v)\\
&=\langle y,v\rangle y-\langle y,y\rangle v\\
&=0\cdot y-1\cdot v\\
&=-v.
\end{aligned}
\]
Consequently,
$J_y^2=-\id_{T_yS^2}.$
Define the action of complex numbers on $T_yS^2$ by
\[
(a+bi)v=av+bJ_y(v),\qquad a,b\in\RR.
\]
Since $J_y^2=-\id$, for arbitrary
$a+bi,\ c+di\in\mathbb C$
and $v\in T_yS^2$, we have
\[
\begin{aligned}
(a+bi)((c+di)v)&=(a+bJ_y)(cv+dJ_yv)\\
&=acv+adJ_yv+bcJ_yv+bdJ_y^2v\\
&=(ac-bd)v+(ad+bc)J_yv\\
&=((a+bi)(c+di))v.
\end{aligned}
\]
The remaining complex vector space axioms follow directly from the real linearity of $J_y$.
Thus each $T_yS^2$ acquires a one-dimensional complex vector space structure.
Moreover, since
$(y,v)\longmapsto y\times v$
is continuous, these complex structures vary continuously with $y$.
More concretely, on any open subset of $S^2$ admitting a nowhere-zero local tangent vector field $w$, the vector fields $w(y)$ and $J_yw(y)$ form a real basis of the tangent plane, so that
\[
(y,a+bi)\longmapsto aw(y)+bJ_yw(y)
\]
gives a local complex line bundle trivialization. Thus $TS^2$ becomes a complex line bundle.
Therefore its pullback $h^*TS^2$ along $h$ is also a complex line bundle.
We now compute its first Chern class.
By the relation between the top Chern class of a complex vector bundle and the Euler class of its underlying oriented real vector bundle~\cite[Proposition~3.13(c), p.~91]{hatcher-vb},
$c_1(TS^2)=e(TS^2),$
where the Euler class uses the orientation given by the complex structure above.
By the computation of the Euler class of the tangent bundle of a sphere~\cite[Proposition~3.14, p.~92]{hatcher-vb},
$\langle e(TS^2),[S^2]\rangle=2$
holds with compatible orientations.
Since
$H^2(S^2;\ZZ)\cong\ZZ$
and $u$ is a generator, there is an integer $j$ such that
$e(TS^2)=ju.$
Since $u$ is a generator,
$\langle u,[S^2]\rangle=\pm1.$
Therefore
\[
\begin{aligned}
2&=\langle e(TS^2),[S^2]\rangle\\
&=\langle ju,[S^2]\rangle\\
&=j\langle u,[S^2]\rangle,
\end{aligned}
\]
so
$j=\pm2.$
Hence
$e(TS^2)=\pm2u.$
By naturality of the Euler class and the relation between the first Chern class and the Euler class~\cite[Proposition~3.13(a), (c), p.~91]{hatcher-vb},
\[
\begin{aligned}
c_1(h^*TS^2)&=e(h^*TS^2)\\
&=h^*e(TS^2)\\
&=h^*(\pm2u)\\
&=\pm2h^*u.
\end{aligned}
\]
By the assumption
$2h^*u=0,$
we obtain
\begin{equation}\label{eq:free-vanishing-chern}
c_1(h^*TS^2)=0.
\end{equation}
Since $X$ is a CW complex, by the classification of complex line bundles by the first Chern class~\cite[Proposition~3.10, p.~86]{hatcher-vb}, the assignment
$L\longmapsto c_1(L)$
gives a bijection between the set of isomorphism classes of complex line bundles over $X$ and
$H^2(X;\ZZ)$
.
Since the trivial complex line bundle
$X\times\mathbb C\longrightarrow X$
has zero first Chern class, by \eqref{eq:free-vanishing-chern} there exists a complex line bundle isomorphism
$\Phi:X\times\mathbb C\xrightarrow{\cong}h^*TS^2.$
Write the pullback bundle explicitly as
\[
h^*TS^2=\{(x,w)\in X\times\RR^3\mid\langle h(x),w\rangle=0\}.
\]
Since $\Phi$ covers the identity map on $X$, for each $x\in X$ we can write
$\Phi(x,1)=(x,w(x)),$
where
$w(x)\in T_{h(x)}S^2.$
Since $\Phi$ is a complex linear isomorphism on each fiber and $1\in\mathbb C$ is nonzero,
$w(x)\ne0$
for all $x\in X$.
The map
$w:X\longrightarrow\RR^3$
is continuous because it is the composite of the continuous map
\[
X\longrightarrow X\times\mathbb C,\qquad x\longmapsto(x,1)
\]
with $\Phi$ and the projection onto the second coordinate.
Define
\[
v:X\longrightarrow\RR^3,\qquad v(x)=\frac{w(x)}{\|w(x)\|}.
\]
Since $\|w(x)\|>0$, the above map is well defined and continuous.
Since $w(x)\in T_{h(x)}S^2$, we also have
$\langle h(x),w(x)\rangle=0.$
Consequently,
\[
\begin{aligned}
\langle h(x),v(x)\rangle&=\left\langle h(x),\frac{w(x)}{\|w(x)\|}\right\rangle\\
&=\frac{\langle h(x),w(x)\rangle}{\|w(x)\|}\\
&=0.
\end{aligned}
\]
Moreover,
\[
\|v(x)\|=\frac{\|w(x)\|}{\|w(x)\|}=1.
\]
Thus $h(x)$ and $v(x)$ are two mutually orthogonal unit vectors in $\RR^3$.
Now define
$L:X\times I\longrightarrow\RR^3$
by
\begin{equation}\label{eq:free-antipodal-homotopy}
L(x,s)=\cos(\pi s)h(x)+\sin(\pi s)v(x).
\end{equation}
Since $h,v$ are continuous and scalar multiplication and vector addition are continuous, $L$ is continuous.
We now verify that the values of $L$ always lie in $S^2$.
By bilinearity of the inner product,
\[
\begin{aligned}
\|L(x,s)\|^2&=\|\cos(\pi s)h(x)+\sin(\pi s)v(x)\|^2\\
&=\cos^2(\pi s)\|h(x)\|^2\\
&\quad+2\cos(\pi s)\sin(\pi s)\langle h(x),v(x)\rangle\\
&\quad+\sin^2(\pi s)\|v(x)\|^2.
\end{aligned}
\]
Since
\[
\|h(x)\|=\|v(x)\|=1,\qquad\langle h(x),v(x)\rangle=0,
\]
we obtain
\[
\begin{aligned}
\|L(x,s)\|^2&=\cos^2(\pi s)+\sin^2(\pi s)\\
&=1.
\end{aligned}
\]
Hence
$L(x,s)\in S^2.$
Thus \eqref{eq:free-antipodal-homotopy} actually defines a continuous map
$L:X\times I\longrightarrow S^2.$
At $s=0$,
\[
\begin{aligned}
L(x,0)&=\cos0\,h(x)+\sin0\,v(x)\\
&=h(x).
\end{aligned}
\]
At $s=1$,
\[
\begin{aligned}
L(x,1)&=\cos\pi\,h(x)+\sin\pi\,v(x)\\
&=-h(x)\\
&=\mathcal A(h(x)).
\end{aligned}
\]
Thus $L$ gives a homotopy between $h$ and $\mathcal A\circ h$, that is,
$h\simeq\mathcal A\circ h.$
Sufficiency follows.
Finally, take
$X=M=\RP^3\#\RP^3.$
By Lemma~\ref{lem:cohomology-endpoints},
$H^2(M;\ZZ)\cong\ZZ/2\oplus\ZZ/2.$
Thus every element $\alpha$ of this cohomology group satisfies
$2\alpha=0.$
In particular, for any continuous map $h:M\to S^2$,
$h^*u\in H^2(M;\ZZ)$
satisfies
$2h^*u=0.$
Applying the sufficiency already proved gives
$h\simeq\mathcal A\circ h.$
\end{proof}

\begin{theorem}[Free homotopy classification of the trivial fundamental-group type]\label{thm:free-trivial}
For any
$[f]_*\in\mathcal T$
and any $g\in G$, we have
\begin{equation}\label{eq:free-trivial-action}
g\cdot[f]_*=[f]_*.
\end{equation}
Thus the action of the target fundamental group $G$ on $\mathcal T$ is trivial.
The restriction of the forgetful map
$Q|_{\mathcal T}:\mathcal T\longrightarrow Q(\mathcal T)$
is a bijection, and hence there is a bijection of sets
\begin{equation}\label{eq:free-trivial}
Q(\mathcal T)\cong\bigsqcup_{\alpha\in H^2(M;\ZZ)}\ZZ\cong\bigsqcup_{\alpha\in(\ZZ/2)^2}\ZZ.
\end{equation}
In particular, for two pointed maps $f,f'$ inducing trivial fundamental-group homomorphisms,
$f\simeq f'\quad\Longleftrightarrow\quad f\simeq_* f'.$
\end{theorem}
\begin{proof}
Take any $[f]_*\in\mathcal T$ and choose a pointed representative
$f:(M,x_0)\longrightarrow(M,x_0).$
Since
$f_{*,1}=\varphi_{\mathrm{triv}},$
by Lemma~\ref{lem:cohomology-endpoints}, there exists a unique pointed lift
$\widehat f:M\longrightarrow\widetilde M$
satisfying
$p\circ\widehat f=f,\qquad\widehat f(x_0)=\widetilde x_0.$
On the other hand, write
$F_f:\widetilde M\longrightarrow\widetilde M$
for the unique pointed lift of $f\circ p$.
Note that $\widehat f$ and $F_f$ have different domains.
We now prove
\begin{equation}\label{eq:free-two-lifts}
F_f=\widehat f\circ p.
\end{equation}
Indeed,
\[
\begin{aligned}
p\circ(\widehat f\circ p)&=(p\circ\widehat f)\circ p\\
&=f\circ p.
\end{aligned}
\]
Moreover,
\[
\begin{aligned}
(\widehat f\circ p)(\widetilde x_0)&=\widehat f(p(\widetilde x_0))\\
&=\widehat f(x_0)\\
&=\widetilde x_0.
\end{aligned}
\]
Thus $\widehat f\circ p$ is a pointed lift of $f\circ p$.
Uniqueness of the pointed lift gives \eqref{eq:free-two-lifts}.
Since
$\widetilde M=S^2\times\RR,$
we can write $\widehat f$ uniquely as
$\widehat f(x)=(h(x),V(x)),$
where
\[
h:M\longrightarrow S^2,\qquad V:M\longrightarrow\RR
\]
are both continuous maps.
We now prove that for every $g\in G$, there is an ordinary homotopy
$\widehat f\simeq g\circ\widehat f$
as maps from $M$ to $\widetilde M$.
By the normal form for elements of $G$, it suffices to consider the following two cases.
\par\addvspace{5pt plus 1pt minus 1pt}\noindent\textbf{Case one: $g=t^r$, where $r\in\ZZ$.}\par\nobreak\vspace{2pt}\nobreak
By \eqref{eq:free-deck},
$t^r(y,v)=(y,v-2r).$
Therefore
\[
\begin{aligned}
(g\circ\widehat f)(x)&=t^r(h(x),V(x))\\
&=(h(x),V(x)-2r).
\end{aligned}
\]
Define
$K:M\times I\longrightarrow S^2\times\RR$
by
$K(x,s)=(h(x),V(x)-2rs).$
Since both coordinate components are continuous, $K$ is continuous.
At $s=0$,
\[
K(x,0)=(h(x),V(x))=\widehat f(x).
\]
At $s=1$,
\[
\begin{aligned}
K(x,1)&=(h(x),V(x)-2r)\\
&=(g\circ\widehat f)(x).
\end{aligned}
\]
Hence
$\widehat f\simeq g\circ\widehat f.$
\par\addvspace{5pt plus 1pt minus 1pt}\noindent\textbf{Case two: $g=t^ra$, where $r\in\ZZ$.}\par\nobreak\vspace{2pt}\nobreak
By \eqref{eq:free-deck},
$t^ra(y,v)=(-y,-v-2r).$
Therefore
\[
\begin{aligned}
(g\circ\widehat f)(x)&=t^ra(h(x),V(x))\\
&=(-h(x),-V(x)-2r).
\end{aligned}
\]
By Lemma~\ref{lem:free-antipodal}, every map $h:M\to S^2$ is homotopic to its antipodal postcomposition.
Thus there is a continuous map
$L:M\times I\longrightarrow S^2$
satisfying
\[
L(x,0)=h(x),\qquad L(x,1)=-h(x).
\]
Define
$K:M\times I\longrightarrow S^2\times\RR$
by
\[
K(x,s)=\bigl(L(x,s),(1-s)V(x)+s(-V(x)-2r)\bigr).
\]
Since $L,V$ are continuous, $K$ is continuous.
At $s=0$,
\[
\begin{aligned}
K(x,0)&=(L(x,0),V(x))\\
&=(h(x),V(x))\\
&=\widehat f(x).
\end{aligned}
\]
At $s=1$,
\[
\begin{aligned}
K(x,1)&=(L(x,1),-V(x)-2r)\\
&=(-h(x),-V(x)-2r)\\
&=(g\circ\widehat f)(x).
\end{aligned}
\]
Hence again
$\widehat f\simeq g\circ\widehat f.$
In summary, for any $g\in G$ there exists a homotopy
$K:M\times I\longrightarrow\widetilde M$
satisfying
\begin{equation}\label{eq:free-trivial-homotopy}
K(x,0)=\widehat f(x),\qquad K(x,1)=g\widehat f(x).
\end{equation}
We now turn this into an equivariant homotopy corresponding to the trivial fundamental-group homomorphism.
Define
$\widetilde K:\widetilde M\times I\longrightarrow\widetilde M$
by
$\widetilde K(z,s)=K(p(z),s).$
This map is continuous because it is the composite of
$\widetilde M\times I\xrightarrow{p\times\id_I}M\times I\xrightarrow{K}\widetilde M$
.
For any $h_0\in G$, $z\in\widetilde M$ and $s\in I$, since
$p(h_0z)=p(z),$
we have
\[
\begin{aligned}
\widetilde K(h_0z,s)&=K(p(h_0z),s)\\
&=K(p(z),s)\\
&=\widetilde K(z,s).
\end{aligned}
\]
But
$\varphi_{\mathrm{triv}}(h_0)=1_G,$
hence
\[
\begin{aligned}
\widetilde K(h_0z,s)&=\widetilde K(z,s)\\
&=\varphi_{\mathrm{triv}}(h_0)\widetilde K(z,s).
\end{aligned}
\]
Thus $\widetilde K$ is a $\varphi_{\mathrm{triv}}$-equivariant homotopy.
By \eqref{eq:free-two-lifts} and \eqref{eq:free-trivial-homotopy}, its initial endpoint satisfies
\[
\begin{aligned}
\widetilde K(z,0)&=K(p(z),0)\\
&=\widehat f(p(z))\\
&=F_f(z).
\end{aligned}
\]
Its final endpoint satisfies
\[
\begin{aligned}
\widetilde K(z,1)&=K(p(z),1)\\
&=g\widehat f(p(z))\\
&=gF_f(z).
\end{aligned}
\]
Therefore
$F_f\simeq_{\varphi_{\mathrm{triv}}}g\circ F_f.$
Thus
$[F_f]=[g\circ F_f]$
in
\begin{equation}\label{eq:free-trivial-components}
\pi_0(\mathcal E_{\varphi_{\mathrm{triv}}}).
\end{equation}
On the other hand, by Lemma~\ref{lem:free-deck-action}, since
$\operatorname{Ad}_g\circ\varphi_{\mathrm{triv}}=\varphi_{\mathrm{triv}},$
under the bijection
$[M,M]_*^{\varphi_{\mathrm{triv}}}\cong\pi_0(\mathcal E_{\varphi_{\mathrm{triv}}})$
the pointed homotopy class
$[f]_*$
corresponds to $[F_f]$, whereas
$g\cdot[f]_*$
corresponds to $[g\circ F_f]$.
By \eqref{eq:free-trivial-components} and injectivity of this correspondence, we obtain
$g\cdot[f]_*=[f]_*.$
Since $[f]_*\in\mathcal T$ and $g\in G$ were arbitrary, the action of $G$ on $\mathcal T$ is trivial.
We now prove that
$Q|_{\mathcal T}:\mathcal T\longrightarrow Q(\mathcal T)$
is a bijection.
By the definition of the image set $Q(\mathcal T)$, the restriction is surjective.
To prove injectivity, suppose
$[f]_*,[f']_*\in\mathcal T$
satisfy
$Q([f]_*)=Q([f']_*).$
By Lemma~\ref{lem:free-orbit}, there exists $g\in G$ such that
$[f']_* = g\cdot[f]_*.$
However, by the triviality of the action already proved,
$g\cdot[f]_*=[f]_*.$
Therefore
$[f']_*=[f]_*.$
Thus $Q|_{\mathcal T}$ is injective and hence bijective.
By Theorem~\ref{thm:classification}(1),
$\mathcal T\cong\bigsqcup_{\alpha\in H^2(M;\ZZ)}\ZZ.$
Combining this with the isomorphism given by Lemma~\ref{lem:cohomology-endpoints},
$H^2(M;\ZZ)\cong(\ZZ/2)^2,$
we finally obtain
$Q(\mathcal T)\cong\bigsqcup_{\alpha\in(\ZZ/2)^2}\ZZ.$
\end{proof}
\smallskip\noindent\textbf{Remark.}
Theorem~\ref{thm:free-trivial} not only shows that the primary invariant
$\alpha\in H^2(M;\ZZ)$
is invariant under the action of the target fundamental group, but also shows that every integer secondary coordinate $n$ remains unchanged.
More precisely, if
\[
[f]_*=(\alpha,n),\qquad [f']_*=(\alpha',n'),
\]
then
\[
f\simeq f'\quad\Longleftrightarrow\quad\alpha=\alpha'\text{ and }n=n'.
\]
Thus in the trivial fundamental-group type, forgetting the basepoint introduces no new identifications of homotopy classes.

\subsection{Free homotopy classification of the double-reflection type}
We now compute the action of the target fundamental group $G$ on the classification set $\mathcal R$ of the double-reflection type.
Throughout we use the compatible reference maps and generators chosen in Lemma~\ref{lem:compatible-odd-reference-family} and Notation~\ref{def:secondary-coordinates}, so that the composition formula of Theorem~\ref{cor:double-reflection-composition} holds.
\begin{lemma}[The change-of-basepoint action on the double-reflection type]\label{lem:free-reflections}
For any
\[
r,m,k\in\ZZ,\qquad d\in2\ZZ+1,\qquad c\in\ZZ/|d|\ZZ,
\]
we have
\begin{equation}\label{eq:free-generator-action}
\begin{aligned}
t^r\cdot(A(m,k),d,c)&=(A(m,k+2r),d,c),\\
t^ra\cdot(A(m,k),d,c)&=(A(-m,2r-k),-d,c).
\end{aligned}
\end{equation}
In particular, the action of the target fundamental group preserves $\mathcal R$ and leaves its secondary coordinate $c$ unchanged.
\end{lemma}
\begin{proof}
Take any $g\in G$.
By Lemma~\ref{lem:dihedral-normal-form}, there are unique
$r\in\ZZ,\qquad\delta\in\{0,1\},$
such that
$g=t^ra^\delta.$
Put
$\varepsilon=(-1)^\delta.$
By the definition of the sign homomorphism $\rho$,
$\rho(g)=\varepsilon.$
First compute the inner automorphism
$\operatorname{Ad}_g:G\longrightarrow G.$
Since
$g=t^ra^\delta,$
we have
$g^{-1}=a^{-\delta}t^{-r}.$
Since $a^2=1_G$, we have
$a^{-1}=a.$
When $\delta=0$,
$a^\delta ta^{-\delta}=t.$
When $\delta=1$,
$a^\delta ta^{-\delta}=ata=t^{-1}.$
Thus in both cases
$a^\delta ta^{-\delta}=t^\varepsilon.$
Consequently,
\begin{equation}\label{eq:free-inner-t}
\begin{aligned}
\operatorname{Ad}_g(t)&=gtg^{-1}\\
&=t^ra^\delta ta^{-\delta}t^{-r}\\
&=t^rt^\varepsilon t^{-r}\\
&=t^{r+\varepsilon-r}\\
&=t^\varepsilon.
\end{aligned}
\end{equation}
We now compute $\operatorname{Ad}_g(a)$.
Since $a^\delta$ commutes with $a$,
$a^\delta aa^{-\delta}=a.$
Since the relation
$at^j=t^{-j}a$
holds for every $j\in\ZZ$, taking $j=-r$ gives
$at^{-r}=t^ra.$
Therefore
\begin{equation}\label{eq:free-inner-a}
\begin{aligned}
\operatorname{Ad}_g(a)&=gag^{-1}\\
&=t^ra^\delta aa^{-\delta}t^{-r}\\
&=t^rat^{-r}\\
&=t^r(t^ra)\\
&=t^{2r}a.
\end{aligned}
\end{equation}
Since $t,a$ generate $G$, the group homomorphism $\operatorname{Ad}_g$ is uniquely determined by its values on $t,a$.
By \eqref{eq:free-inner-t}, \eqref{eq:free-inner-a}, and the definition of the group endomorphism corresponding to the matrix parameter,
\begin{equation}\label{eq:free-inner}
\operatorname{Ad}_g=\varphi_{A(\varepsilon,2r)}.
\end{equation}
Next compute
$\beta_g=g\cdot[\id_M]_*.$
The identity map induces the identity homomorphism on the fundamental group,
$(\id_M)_{*,1}=\id_G.$
It also induces the identity homomorphism on the second homotopy group, namely
$(\id_M)_{*,2}=\id_{\ZZ}.$
By Lemma~\ref{lem:free-deck-action}(1), the fundamental-group homomorphism induced by $\beta_g$ is
$\operatorname{Ad}_g\circ\id_G=\operatorname{Ad}_g=\varphi_{A(\varepsilon,2r)}.$
and the multiplier of its induced homomorphism on the second homotopy group is
$\rho(g)\cdot1=\varepsilon.$
Since
$\varepsilon\in\{1,-1\},$
we have
$|\varepsilon|=1.$
By Theorem~\ref{thm:classification}(2), when the induced fundamental-group homomorphism is fixed to be
$\varphi_{A(\varepsilon,2r)}$
and the multiplier of the induced homomorphism on the second homotopy group is fixed to be $\varepsilon$, the corresponding set of pointed homotopy classes is in bijection with
$\ZZ/|\varepsilon|\ZZ=\ZZ/\ZZ=\{0\}$
.
Hence there is exactly one such pointed homotopy class, and its secondary coordinate must be $0$.
Consequently,
\begin{equation}\label{eq:free-identity-orbit}
\beta_g=(A(\varepsilon,2r),\varepsilon,0).
\end{equation}
Now take any
$[f]_*=(A(m,k),d,c)\in\mathcal R.$
By Lemma~\ref{lem:free-deck-action}(2),
$g\cdot[f]_* =\beta_g\star[f]_*.$
Substituting \eqref{eq:free-identity-orbit} gives
\begin{equation}\label{eq:free-action-substitution}
g\cdot(A(m,k),d,c)=(A(\varepsilon,2r),\varepsilon,0)\star(A(m,k),d,c).
\end{equation}
To compute the right-hand side, apply Theorem~\ref{cor:double-reflection-composition} directly.
Write the coordinates of the outer map in Theorem~\ref{cor:double-reflection-composition} as
$(\mathbf C,e,b),$
and those of the inner map as
$(\mathbf B,d,c).$
Then that theorem gives the composition formula
\begin{equation}\label{eq:free-composition-recalled}
\begin{aligned}
(\mathbf C,e,b)\star(\mathbf B,d,c)
&=\left(\mathbf C\mathbf B,ed,\left[d\det(\mathbf B)\widetilde b+e^2\widetilde c\right]_{|ed|}\right),
\end{aligned}
\end{equation}
where
$\widetilde b,\widetilde c\in\ZZ$
are integer representatives of
$b\in\ZZ/|e|\ZZ,\qquad c\in\ZZ/|d|\ZZ$
respectively, and $[j]_N$ denotes the residue class of the integer $j$ modulo $N$.
In \eqref{eq:free-action-substitution},
\[
\mathbf C=A(\varepsilon,2r),\qquad e=\varepsilon,\qquad b=0,
\]
whereas
$\mathbf B=A(m,k).$
First compute the matrix coordinate:
\[
\begin{aligned}
\mathbf C\mathbf B&=A(\varepsilon,2r)A(m,k)\\
&=\begin{pmatrix}\varepsilon&2r\\0&1\end{pmatrix}\begin{pmatrix}m&k\\0&1\end{pmatrix}\\
&=\begin{pmatrix}\varepsilon m+2r\cdot0&\varepsilon k+2r\cdot1\\0\cdot m+1\cdot0&0\cdot k+1\cdot1\end{pmatrix}\\
&=\begin{pmatrix}\varepsilon m&\varepsilon k+2r\\0&1\end{pmatrix}\\
&=A(\varepsilon m,\varepsilon k+2r).
\end{aligned}
\]
Thus the matrix coordinate after composition is
\begin{equation}\label{eq:free-composed-matrix}
A(\varepsilon m,\varepsilon k+2r).
\end{equation}
Next, the multiplier of the induced homomorphism on the second homotopy group is
\begin{equation}\label{eq:free-composed-degree}
ed=\varepsilon d.
\end{equation}
Since $d$ is odd and $\varepsilon=\pm1$, $\varepsilon d$ is still odd.
Finally compute the secondary coordinate.
Since $b=0\in\ZZ/|\varepsilon|\ZZ$, we may choose its integer representative to be
$\widetilde b=0.$
At the same time,
$\det(\mathbf B)=\det A(m,k)=m.$
Thus by \eqref{eq:free-composition-recalled}, the secondary coordinate after composition is
\[
\begin{aligned}
c_{\mathrm{new}}&=\left[d\det(\mathbf B)\widetilde b+e^2\widetilde c\right]_{|ed|}\\
&=\left[dm\cdot0+\varepsilon^2\widetilde c\right]_{|\varepsilon d|}.
\end{aligned}
\]
Since
$\varepsilon^2=1,\qquad|\varepsilon d|=|d|,$
we have
\begin{equation}\label{eq:free-composed-secondary}
\begin{aligned}
c_{\mathrm{new}}&=[\widetilde c]_{|d|}\\
&=c.
\end{aligned}
\end{equation}
There is no issue of identifying different residue class groups here, since
\[
\ZZ/|\varepsilon d|\ZZ=\ZZ/|d|\ZZ.
\]
Combining \eqref{eq:free-action-substitution}, \eqref{eq:free-composed-matrix}, \eqref{eq:free-composed-degree}, and \eqref{eq:free-composed-secondary}, we obtain
\begin{equation}\label{eq:free-full-action}
g\cdot(A(m,k),d,c)=(A(\varepsilon m,\varepsilon k+2r),\varepsilon d,c).
\end{equation}
Since $m,k,r\in\ZZ$, $\varepsilon d$ is odd, and
$c\in\ZZ/|\varepsilon d|\ZZ,$
the right-hand side still belongs to $\mathcal R$. This also proves that $\mathcal R$ is invariant under the action of $G$.
When $\delta=0$,
$g=t^r,\qquad\varepsilon=1.$
Substituting into \eqref{eq:free-full-action} gives
\[
t^r\cdot(A(m,k),d,c)=(A(m,k+2r),d,c).
\]
When $\delta=1$,
$g=t^ra,\qquad\varepsilon=-1.$
Substituting into \eqref{eq:free-full-action} gives
\[
t^ra\cdot(A(m,k),d,c)=(A(-m,-k+2r),-d,c).
\]
That is,
\[
t^ra\cdot(A(m,k),d,c)=(A(-m,2r-k),-d,c).
\]
\end{proof}
\subsection{The free homotopy classification theorem}
By Lemma~\ref{lem:free-orbit}, the set of free homotopy classes of self-maps of $M$ is the orbit set for the action of the fundamental group $G$ on $[M,M]_*$. 
Theorem~\ref{thm:free-trivial} and Lemma~\ref{lem:free-reflections} compute the restrictions of this action to $\mathcal T$ and $\mathcal R$, respectively. Thus we can give the complete free homotopy classification.
\begin{theorem}[Free homotopy classification of self-maps of \(\RP^3\#\RP^3\)]\label{thm:free-classification}
Let $M=\mathbb{RP}^3\#\mathbb{RP}^3$. For any homomorphism $\varphi:G\to G$ and integer $d$, write
\[
[M,M]^{\varphi,d}=\left\{[f]\in[M,M]\;\middle|\;
\begin{aligned}
&[f]\text{ has a pointed representative } f_0,\\
&(f_0)_{*,1}=\varphi,\\
&(f_0)_{*,2}=d\,\mathrm{id}_{\pi_2(M,x_0)}
\end{aligned}\right\},
\]
where $[M,M]$ is the set of free homotopy classes of self-maps of $M$. The classification is as follows.
\begin{enumerate}[label=\textnormal{(\arabic*)}]
\item If $\varphi$ is trivial and $d=0$, then, after choosing a reference homotopy class in each of the four families distinguished by $\alpha\in H^2(M;\ZZ)\cong(\ZZ/2\ZZ)^2$, there is a bijection of sets
\[
[M,M]^{\varphi,0}\cong\bigsqcup_{\alpha\in(\ZZ/2\ZZ)^2}\ZZ.
\]
\item For any integers $m,k$ and odd integer $d$, let $\varphi:G\to G$ be the homomorphism determined by
\[
\varphi(t)=t^m,\qquad \varphi(a)=t^ka.
\]
After choosing a reference homotopy class, there is a bijection of sets
\[
[M,M]^{\varphi,d}\cong\ZZ/|d|\ZZ.
\]
Furthermore, suppose that $\varphi':G\to G$ is determined by $\varphi'(t)=t^{m'}$ and $\varphi'(a)=t^{k'}a$, and that $d'$ is odd. Then
\[
[M,M]^{\varphi,d}=[M,M]^{\varphi',d'}
\]
if and only if
\[
\begin{cases}
\operatorname{sgn}(d)m=\operatorname{sgn}(d')m',\\
k\equiv k'\pmod 2,\\
|d|=|d'|.
\end{cases}
\]
Otherwise, these two sets are disjoint.
\end{enumerate}
In all other cases, $[M,M]^{\varphi,d}=\varnothing$.

Write $\ZZ_{\mathrm{odd}}^{>0}=\{1,3,5,\ldots\}$. The unique pointed normal form of a free homotopy class of double-reflection type is $(A(\mu,\kappa),D,c)$, where $c\in\ZZ/D\ZZ$ is the original secondary coordinate and
\begin{equation}\label{eq:free-normal-parameters}
\mu=\operatorname{sgn}(d)m,\qquad
\kappa\equiv k\pmod 2,\quad \kappa\in\{0,1\},\qquad D=|d|,
\end{equation}
and hence we obtain a bijection of sets
\begin{equation}\label{eq:free-classification}
[M,M]\cong\left(\bigsqcup_{\alpha\in(\ZZ/2)^2}\ZZ\right)
\sqcup\left(\bigsqcup_{\substack{\mu\in\ZZ,\ \kappa\in\{0,1\}\\D\in\ZZ_{\mathrm{odd}}^{>0}}}\ZZ/D\ZZ\right).
\end{equation}
The displayed bijection depends on the choices of reference homotopy classes and generators.
\end{theorem}
\begin{proof}
By Theorem~\ref{thm:classification},
$[M,M]_* = \mathcal T\sqcup\mathcal R.$
By Lemma~\ref{lem:free-orbit}, the forgetful map
$Q:[M,M]_*\longrightarrow[M,M]$
is surjective, so
\[
\begin{aligned}
[M,M]&=Q([M,M]_*)\\
&=Q(\mathcal T\sqcup\mathcal R)\\
&=Q(\mathcal T)\cup Q(\mathcal R).
\end{aligned}
\]
First we prove that the above union is disjoint.
Suppose there exists
$[f]\in Q(\mathcal T)\cap Q(\mathcal R).$
By the definition of an image set, there are
$[f_0]_*\in\mathcal T,\qquad[f_1]_*\in\mathcal R$
satisfying
$Q([f_0]_*)=[f]=Q([f_1]_*).$
By Lemma~\ref{lem:free-orbit}, there exists $g\in G$ such that
$[f_1]_* = g\cdot[f_0]_*.$
However, by Theorem~\ref{thm:free-trivial},
$g\cdot[f_0]_*=[f_0]_*.$
Thus
$[f_1]_*=[f_0]_*.$
This means that
$[f_0]_*\in\mathcal T\cap\mathcal R,$
contradicting the decomposition in Theorem~\ref{thm:classification},
$[M,M]_* = \mathcal T\sqcup\mathcal R.$

Therefore
$Q(\mathcal T)\cap Q(\mathcal R)=\varnothing,$
and hence
\begin{equation}\label{eq:free-disjoint-images}
[M,M]=Q(\mathcal T)\sqcup Q(\mathcal R).
\end{equation}
By Theorem~\ref{thm:free-trivial},
\begin{equation}\label{eq:free-trivial-images}
Q(\mathcal T)\cong\bigsqcup_{\alpha\in(\ZZ/2)^2}\ZZ.
\end{equation}
Hence it remains only to compute $Q(\mathcal R)$.
By Lemma~\ref{lem:free-reflections}, $\mathcal R$ is invariant under the action of $G$. Thus Lemma~\ref{lem:free-orbit} gives a bijection
\begin{equation}\label{eq:free-reflection-orbits}
\mathcal R/G\xrightarrow{\cong}Q(\mathcal R).
\end{equation}
We now prove that every orbit in $\mathcal R$ contains exactly one pointed homotopy class having the stated normal form.

Take any
$(A(m,k),d,c)\in\mathcal R.$
Since
$d\in2\ZZ+1,$
we have $d\ne0$.
Thus we can define
\[
s=\operatorname{sgn}(d)=\begin{cases}1,&d>0,\\-1,&d<0.\end{cases}
\]
Put
$D=|d|.$
Since $d$ is odd, $D$ is a positive odd integer.
By the definitions of absolute value and the sign function,
\begin{equation}\label{eq:free-positive-degree}
sd=|d|=D.
\end{equation}
Take the unique
$\kappa\in\{0,1\}$
such that
$\kappa\equiv k\pmod2.$
Since $s\in\{1,-1\}$,
$s\equiv1\pmod2.$
Consequently,
$sk\equiv k\equiv\kappa\pmod2.$
Thus
$\kappa-sk\equiv0\pmod2.$
This means that there is an integer $r$ such that
$2r=\kappa-sk.$
Explicitly, we may take
\begin{equation}\label{eq:free-normal-translation}
r=\frac{\kappa-sk}{2}\in\ZZ.
\end{equation}
Define the group element
\[
g=\begin{cases}t^r,&s=1,\\t^ra,&s=-1.\end{cases}
\]
By the definition of the sign homomorphism,
$\rho(g)=s.$
Applying the unified action formula \eqref{eq:free-full-action} in Lemma~\ref{lem:free-reflections} gives
\[
g\cdot(A(m,k),d,c)=(A(sm,sk+2r),sd,c).
\]
By \eqref{eq:free-positive-degree},
$sd=D.$
By \eqref{eq:free-normal-translation},
\[
\begin{aligned}
sk+2r&=sk+(\kappa-sk)\\
&=\kappa.
\end{aligned}
\]
Put
$\mu=sm=\operatorname{sgn}(d)m.$
Consequently,
\begin{equation}\label{eq:free-normal-existence}
g\cdot(A(m,k),d,c)=(A(\mu,\kappa),D,c).
\end{equation}
Since
\[
\mu\in\ZZ,\qquad\kappa\in\{0,1\},\qquad D\in\ZZ_{\mathrm{odd}}^{>0},
\]
and
$\ZZ/|d|\ZZ=\ZZ/D\ZZ,$
the right-hand side has the required normal form.
By \eqref{eq:free-normal-existence}, it belongs to the same $G$-orbit as the original pointed homotopy class.
Thus every orbit contains at least one pointed homotopy class having the stated normal form.

Suppose two pointed homotopy classes in normal form,
\[
(A(\mu,\kappa),D,c),
\]
\[
(A(\mu',\kappa'),D',c')
\]
belong to the same $G$-orbit, where
\[
\mu,\mu'\in\ZZ,
\]
\[
\kappa,\kappa'\in\{0,1\},
\]
\[
D,D'\in\ZZ_{\mathrm{odd}}^{>0},
\]
and
$c\in\ZZ/D\ZZ,\qquad c'\in\ZZ/D'\ZZ.$
By the definition of belonging to the same orbit, there exists $g\in G$ such that
\begin{equation}\label{eq:free-normal-same-orbit}
g\cdot(A(\mu,\kappa),D,c)=(A(\mu',\kappa'),D',c').
\end{equation}
Write $g$ uniquely as
\[
g=t^ra^\delta,\qquad r\in\ZZ,\quad\delta\in\{0,1\},
\]
and put
$\varepsilon=(-1)^\delta.$
By Lemma~\ref{lem:free-reflections},
\[
\begin{aligned}
g\cdot(A(\mu,\kappa),D,c)&=(A(\varepsilon\mu,\varepsilon\kappa+2r),\varepsilon D,c).
\end{aligned}
\]
Substituting into \eqref{eq:free-normal-same-orbit} gives
\[
\begin{aligned}
(A(\varepsilon\mu,\varepsilon\kappa+2r),\varepsilon D,c)&=(A(\mu',\kappa'),D',c').
\end{aligned}
\]
By uniqueness of the pointed classification coordinates,
\begin{equation}\label{eq:free-unique-degree-sign}
D'=\varepsilon D.
\end{equation}
Since
$D>0,\qquad D'>0,$
if $\varepsilon=-1$, then \eqref{eq:free-unique-degree-sign} gives
$D'=-D<0,$
a contradiction.
Therefore
$\varepsilon=1.$
Consequently,
\begin{equation}\label{eq:free-unique-degree}
D'=D.
\end{equation}
Comparing matrix coordinates gives
$A(\mu,\kappa+2r)=A(\mu',\kappa').$
Expanding the matrices,
\[
\begin{pmatrix}\mu&\kappa+2r\\0&1\end{pmatrix}=\begin{pmatrix}\mu'&\kappa'\\0&1\end{pmatrix}.
\]
Hence
\begin{equation}\label{eq:free-unique-mu}
\mu'=\mu
\end{equation}
and
\begin{equation}\label{eq:free-unique-kappa-shift}
\kappa'=\kappa+2r.
\end{equation}
By \eqref{eq:free-unique-kappa-shift},
$\kappa'-\kappa=2r.$
Since
$\kappa,\kappa'\in\{0,1\},$
we have
$\kappa'-\kappa\in\{-1,0,1\}.$
Moreover, since $r\in\ZZ$, $2r$ is even.
The only even number in $\{-1,0,1\}$ is $0$, so
$\kappa'-\kappa=0.$
Therefore
\begin{equation}\label{eq:free-unique-kappa}
\kappa'=\kappa,\qquad r=0.
\end{equation}
Finally, comparing the secondary coordinates gives
$c'=c$
in
\begin{equation}\label{eq:free-unique-secondary}
\ZZ/D\ZZ.
\end{equation}
Here both secondary coordinates belong to the same residue class group because \eqref{eq:free-unique-degree} has already given $D'=D$.
Combining \eqref{eq:free-unique-degree}, \eqref{eq:free-unique-mu}, \eqref{eq:free-unique-kappa}, and \eqref{eq:free-unique-secondary} gives
\[
\mu'=\mu,\qquad\kappa'=\kappa,\qquad D'=D,\qquad c'=c.
\]
Thus the two normal forms are identical.
This proves uniqueness of the normal form.
In particular, this proof does not use the condition $\mu\ne0$. Thus even when $\mu=0$, the above uniqueness remains valid and no additional free homotopy identifications occur.

Take any
\[
\mu\in\ZZ,\qquad\kappa\in\{0,1\},
\]
\[
D\in\ZZ_{\mathrm{odd}}^{>0},\qquad c\in\ZZ/D\ZZ.
\]
Put
$A=A(\mu,\kappa).$
Since $\mu,\kappa\in\ZZ$, $A$ belongs to the matrix semigroup $\mathscr M$ defined in \cref{def:secondary-coordinates}.
Since $D$ is also odd, by Theorem~\ref{thm:classification}(2) the set of pointed homotopy classes
$[M,M]_*^{\varphi_A,D}$
is nonempty, and after choosing a reference homotopy class there is a bijection of sets
$[M,M]_*^{\varphi_A,D}\cong\ZZ/D\ZZ.$
Using the reference maps and generators chosen in Notation~\ref{def:secondary-coordinates}, there exists a unique pointed homotopy class corresponding to the given residue class $c$.
Thus there exists a pointed map
$f:(M,x_0)\longrightarrow(M,x_0)$
satisfying
$[f]_*=(A(\mu,\kappa),D,c).$
Consequently, the $G$-orbit containing this pointed homotopy class realizes the given normal-form parameters.
Hence all the listed normal-form parameters are realizable.

Put
\[
\mathcal N=\left\{(\mu,\kappa,D,c)\ \middle|\ \begin{aligned}
&\mu\in\ZZ,\ \kappa\in\{0,1\},\\
&D\in\ZZ_{\mathrm{odd}}^{>0},\\
&c\in\ZZ/D\ZZ
\end{aligned}\right\}.
\]
Define the map
$\Psi:\mathcal R/G\longrightarrow\mathcal N$
by
\[
\begin{aligned}
\Psi\bigl(G\cdot(A(m,k),d,c)\bigr)&=\bigl(\operatorname{sgn}(d)m,\kappa,|d|,c\bigr),
\end{aligned}
\]
where $\kappa\in\{0,1\}$ satisfies
$\kappa\equiv k\pmod2.$
By existence of the normal form, every orbit has the above normal form.
If two different pointed homotopy classes are chosen as representatives of the same orbit and each is reduced to normal form, the resulting normal forms still belong to the same orbit.
By uniqueness of the normal form, these two normal forms are identical.
Thus $\Psi$ is independent of the choice of orbit representative and is well defined.
If two orbits have the same image under $\Psi$, then they contain the same normal form. Since two orbits are either equal or disjoint, these orbits are equal.
Thus $\Psi$ is injective.
By the realizability established above, every normal-form parameter is realized by a pointed homotopy class, so $\Psi$ is surjective.
Hence
$\mathcal R/G\cong\mathcal N.$
By the definition of $\mathcal N$,
\[
\mathcal N\cong\bigsqcup_{\substack{\mu\in\ZZ,\ \kappa\in\{0,1\}\\D\in\ZZ_{\mathrm{odd}}^{>0}}}\ZZ/D\ZZ.
\]
Combining this with \eqref{eq:free-reflection-orbits} gives
\begin{equation}\label{eq:free-reflection-images}
Q(\mathcal R)\cong\bigsqcup_{\substack{\mu\in\ZZ,\ \kappa\in\{0,1\}\\D\in\ZZ_{\mathrm{odd}}^{>0}}}\ZZ/D\ZZ.
\end{equation}
Finally, substituting \eqref{eq:free-trivial-images} and \eqref{eq:free-reflection-images} into \eqref{eq:free-disjoint-images} gives
\[
[M,M]\cong\left(\bigsqcup_{\alpha\in(\ZZ/2)^2}\ZZ\right)
\sqcup\left(\bigsqcup_{\substack{\mu\in\ZZ,\ \kappa\in\{0,1\}\\D\in\ZZ_{\mathrm{odd}}^{>0}}}\ZZ/D\ZZ\right).
\]
\end{proof}

By Theorem~\ref{thm:free-classification}, we obtain a necessary and sufficient condition for deciding directly whether two given self-maps are freely homotopic.
Suppose
\[
[f]_*=(A(m,k),d,c),
\]
\[
[f']_*=(A(m',k'),d',c')
\]
are both pointed homotopy classes of the double-reflection type.
Then
\begin{equation}\label{eq:free-test}
f\simeq f'\quad\Longleftrightarrow\quad\begin{cases}
|d'|=|d|,\\
\operatorname{sgn}(d')m'=\operatorname{sgn}(d)m,\\
k'\equiv k\pmod2,\\
c'=c\ \text{in }\ZZ/|d|\ZZ.
\end{cases}
\end{equation}
We now verify this necessary and sufficient condition directly.
First suppose $f\simeq f'$.
By Lemma~\ref{lem:free-orbit}, there exists
$g=t^ra^\delta\in G,$
such that
$[f']_* = g\cdot[f]_*.$
Put
$\varepsilon=(-1)^\delta.$
By Lemma~\ref{lem:free-reflections},
\[
(A(m',k'),d',c')=(A(\varepsilon m,\varepsilon k+2r),\varepsilon d,c).
\]
Comparing coordinates gives
\[
m'=\varepsilon m,\qquad k'=\varepsilon k+2r,
\]
\[
d'=\varepsilon d,\qquad c'=c.
\]
Since $|\varepsilon|=1$,
$|d'|=|\varepsilon d|=|d|.$
Since $d,d'\ne0$,
$\operatorname{sgn}(d')=\operatorname{sgn}(\varepsilon d)=\varepsilon\operatorname{sgn}(d).$
Therefore
\[
\begin{aligned}
\operatorname{sgn}(d')m'&=\varepsilon\operatorname{sgn}(d)\cdot\varepsilon m\\
&=\varepsilon^2\operatorname{sgn}(d)m\\
&=\operatorname{sgn}(d)m.
\end{aligned}
\]
Moreover, since
$\varepsilon\equiv1\pmod2,$
we have
$k'=\varepsilon k+2r\equiv k\pmod2.$
At the same time the secondary coordinates satisfy
$c'=c.$
Hence all four conditions on the right-hand side of \eqref{eq:free-test} hold.
Conversely, suppose all four conditions hold.
Put
$\varepsilon=\operatorname{sgn}(d')\operatorname{sgn}(d).$
Since $d,d'\ne0$,
$\varepsilon\in\{1,-1\}.$
From $|d'|=|d|$, we obtain
\[
\begin{aligned}
\varepsilon d&=\operatorname{sgn}(d')\operatorname{sgn}(d)d\\
&=\operatorname{sgn}(d')|d|\\
&=\operatorname{sgn}(d')|d'|\\
&=d'.
\end{aligned}
\]
Also, from
$\operatorname{sgn}(d')m'=\operatorname{sgn}(d)m,$
multiplying both sides by $\operatorname{sgn}(d')$ gives
\[
\begin{aligned}
m'&=\operatorname{sgn}(d')\operatorname{sgn}(d)m\\
&=\varepsilon m.
\end{aligned}
\]
Since
$k'\equiv k\pmod2$
and
$\varepsilon k\equiv k\pmod2,$
we have
$k'-\varepsilon k\equiv0\pmod2.$
Thus we can define the integer
\[
r=\frac{k'-\varepsilon k}{2}.
\]
Consequently,
$k'=\varepsilon k+2r.$
Take
\[
g=\begin{cases}t^r,&\varepsilon=1,\\t^ra,&\varepsilon=-1.\end{cases}
\]
By Lemma~\ref{lem:free-reflections},
\[
\begin{aligned}
g\cdot[f]_*&=g\cdot(A(m,k),d,c)\\
&=(A(\varepsilon m,\varepsilon k+2r),\varepsilon d,c)\\
&=(A(m',k'),d',c')\\
&=[f']_*.
\end{aligned}
\]
Thus by Lemma~\ref{lem:free-orbit},
$f\simeq f'.$
This proves \eqref{eq:free-test}.
In particular, when two pointed self-maps of the double-reflection type induce the same homomorphism on the fundamental group and the same homomorphism on the second homotopy group, that is,
\[
m'=m,\qquad k'=k,\qquad d'=d
\]
then \eqref{eq:free-test} reduces to
\begin{equation}\label{eq:free-fixed-coordinates}
f\simeq f'\quad\Longleftrightarrow\quad c=c'.
\end{equation}
Thus all $|d|$ distinct pointed homotopy classes obtained for fixed $(A,d)$ in Theorem~\ref{thm:classification}(2) remain pairwise distinct under free homotopy.

\begin{corollary}[Counting classes for a fixed conjugacy class of fundamental-group endomorphisms]\label{cor:free-count}
The conjugacy class of a double-reflection fundamental-group endomorphism \(\varphi_{A(m,k)}\) under conjugation in the target group is completely determined by
$\ell=|m|,\qquad \kappa=k\bmod 2.$
Fixing this conjugacy class and a positive odd integer \(D=|d|\), the number of corresponding free homotopy classes is
\[
\begin{cases}
D,&\ell=0,\\
2D,&\ell>0.
\end{cases}
\]
\end{corollary}

\begin{proof}
Take any $g=t^ra^\delta\in G$, where $r\in\ZZ$ and $\delta\in\{0,1\}$, and put $\varepsilon=(-1)^\delta$.
By \eqref{eq:free-inner} and the composition formula for fundamental-group endomorphisms established above,
\[
\begin{aligned}
\operatorname{Ad}_g\circ\varphi_{A(m,k)}
 &=\varphi_{A(\varepsilon,2r)A(m,k)}\\
 &=\varphi_{A(\varepsilon m,\varepsilon k+2r)}.
\end{aligned}
\]Here the matrix product is computed as follows.

\[
\begin{aligned}
A(\varepsilon,2r)A(m,k)
&=\begin{pmatrix}\varepsilon&2r\\0&1\end{pmatrix}
\begin{pmatrix}m&k\\0&1\end{pmatrix}\\
&=\begin{pmatrix}\varepsilon m&\varepsilon k+2r\\0&1\end{pmatrix}\\
&=A(\varepsilon m,\varepsilon k+2r).
\end{aligned}
\]

Thus conjugation by an element of the target group changes the parameters $(m,k)$ to $(\varepsilon m,\varepsilon k+2r)$.
Since $\varepsilon\in\{1,-1\}$, we have
\[
|\varepsilon m|=|m|,\qquad \varepsilon k+2r\equiv k\pmod2.
\]
Hence $|m|$ and $k\bmod 2$ are invariant under conjugation.
Conversely, suppose $|m'|=|m|$ and $k'\equiv k\pmod2$.
If $m\ne0$, then $m'=m$ or $m'=-m$, so there is $\varepsilon\in\{1,-1\}$ such that $m'=\varepsilon m$.
If $m=0$, then $m'=0$, and we may take $\varepsilon=1$.
Since $\varepsilon\equiv1\pmod2$, we have
$k'-\varepsilon k\equiv k'-k\equiv0\pmod2.$
Thus
\[
r=\frac{k'-\varepsilon k}{2}\in\ZZ.
\]
Take $g=t^r$ when $\varepsilon=1$, and $g=t^ra$ when $\varepsilon=-1$. Then
$\varepsilon m=m',\qquad \varepsilon k+2r=k'.$
Substituting these into the above conjugation formula gives
$\operatorname{Ad}_g\circ\varphi_{A(m,k)}=\varphi_{A(m',k')}.$
Hence $(\ell,\kappa)$ completely determines the conjugacy class.

Now fix the conjugacy class $(\ell,\kappa)$ and the positive odd integer $D$, and count the corresponding free homotopy classes.
By Theorem~\ref{thm:free-classification}, these free homotopy classes correspond bijectively to the normal forms $(A(\mu,\kappa),D,c)$ satisfying
$|\mu|=\ell,\qquad c\in\ZZ/D\ZZ$
 
Indeed, the first parameter of the normal form satisfies $\mu=\operatorname{sgn}(d)m$, so $|\mu|=|m|=\ell$; the second parameter is the fixed $\kappa$, and the third is the fixed $D$.
Conversely, by the realization part of Theorem~\ref{thm:free-classification}, every normal form satisfying these conditions gives a free homotopy class whose induced fundamental-group endomorphism belongs to the fixed conjugacy class and whose multiplier on the second homotopy group has absolute value $D$.

If $\ell=0$, the only possible value is $\mu=0$. The normal forms are then
\[
(A(0,\kappa),D,c),\qquad c\in\ZZ/D\ZZ.
\]
Since $\ZZ/D\ZZ$ has exactly $D$ elements, uniqueness of the normal forms gives exactly $D$ free homotopy classes.

If $\ell>0$, then $\mu=\ell$ or $\mu=-\ell$. The normal forms are then
\[
(A(\ell,\kappa),D,c)\quad\text{or}\quad(A(-\ell,\kappa),D,c),
\qquad c\in\ZZ/D\ZZ.
\]
For each fixed $\mu$, there are exactly $D$ possible secondary coordinates; since $\ell\ne-\ell$, the two sets of normal forms are disjoint.
By uniqueness of normal forms, these $2D$ classes are pairwise distinct; by realizability of normal forms, all these classes exist.
Therefore the number of corresponding free homotopy classes is $2D$.
\end{proof}

\smallskip\noindent\textbf{Remark.}
In particular, after fixing the induced fundamental-group homomorphism $\varphi_{A(m,k)}$, we can also explain directly why the case $m=0$ introduces no additional identifications of secondary coordinates.
Let $g=t^ra^\delta\in G$ and put $\varepsilon=(-1)^\delta$. By the computation in Lemma~\ref{lem:free-reflections},
$\operatorname{Ad}_g\circ\varphi_{A(m,k)}=\varphi_{A(\varepsilon m,\varepsilon k+2r)}.$
By the one-to-one correspondence between matrices and the corresponding group endomorphisms,
$\operatorname{Ad}_g\circ\varphi_{A(m,k)}=\varphi_{A(m,k)}$
if and only if
$\varepsilon m=m,\qquad \varepsilon k+2r=k.$
When $\varepsilon=1$, the second equality gives $2r=0$, and hence $r=0$. Thus $g=1_G$.
When $\varepsilon=-1$, the first equality gives $-m=m$, that is, $m=0$; the second equality gives $-k+2r=k$, that is, $r=k$. Thus $g=t^ka$.
Combining these cases, the subgroup fixing $\varphi_{A(m,k)}$ is
\[
\{g\in G\mid\operatorname{Ad}_g\circ\varphi_{A(m,k)}=\varphi_{A(m,k)}\}
=\begin{cases}
\{1_G\},&m\ne0,\\
\{1_G,t^ka\},&m=0.
\end{cases}
\]
If $m=0$, then by Lemma~\ref{lem:free-reflections},
\[
\begin{aligned}
t^ka\cdot(A(0,k),d,c)
&=(A(-0,2k-k),-d,c)\\
&=(A(0,k),-d,c).
\end{aligned}
\]
Since $d$ is odd, $d\ne0$, and hence $-d\ne d$. Therefore, although this nontrivial group element fixes the induced fundamental-group homomorphism, it changes the sign of the multiplier on the second homotopy group and leaves the secondary coordinate $c$ unchanged.
Thus, with $\varphi_{A(m,k)}$ and $d$ fixed, even when $m=0$ no new identifications between different secondary coordinates arise under free homotopy.

\begin{corollary}[The group of self-homotopy equivalences up to free homotopy]\label{cor:free-homotopy-equivalences}
Let $[M,M]^\times$ be the unit group of the monoid $[M,M]$ of free homotopy classes, namely the group of free homotopy classes of all self-homotopy equivalences of $M$, with multiplication induced by composition. Then there is a group isomorphism
$[M,M]^\times\cong\ZZ/2\ZZ\times\ZZ/2\ZZ.$
More explicitly, in the normal-form coordinates of \cref{thm:free-classification}, the elements of $[M,M]^\times$ are precisely
\[
Q\bigl((A(\mu,\kappa),1,0)\bigr),\qquad
\mu\in\{-1,1\},\quad\kappa\in\{0,1\}.
\]
Writing such an element simply as $(\mu,\kappa)$, the group operation is
\begin{equation}\label{eq:free-unit-product}
(\mu,\kappa)\star(\nu,\lambda)
=(\mu\nu,\overline{\kappa+\lambda}),
\end{equation}
where $\overline{\kappa+\lambda}$ denotes the unique representative in $\{0,1\}$ of $\kappa+\lambda$ modulo $2$.
\end{corollary}

\begin{proof}

By compatibility of composition with free homotopy, for any two pointed self-maps $f,g:M\to M$ we have
\[
\begin{aligned}
Q([f]_*\star[g]_*)&=Q([f\circ g]_*)\\
&=[f\circ g]\\
&=[f]\star[g]\\
&=Q([f]_*)\star Q([g]_*).
\end{aligned}
\]
Also, since
$Q([\id_M]_*)=[\id_M],$
the forgetful map
$Q:[M,M]_*\longrightarrow[M,M]$
is a monoid homomorphism.
By \cref{lem:free-orbit}, $Q$ is surjective.
We now prove
\begin{equation}\label{eq:free-unit-image}
Q\bigl([M,M]_*^\times\bigr)=[M,M]^\times.
\end{equation}
First, if
$[f]_*\in[M,M]_*^\times,$
then there exists $[g]_*\in[M,M]_*$ such that
$[f]_*\star[g]_*=[\id_M]_*$
and
$[g]_*\star[f]_*=[\id_M]_*.$
Since $Q$ is a monoid homomorphism, we obtain
$Q([f]_*)\star Q([g]_*)=[\id_M]$
and
$Q([g]_*)\star Q([f]_*)=[\id_M].$
Therefore
$Q([f]_*)\in[M,M]^\times.$
Hence
$Q\bigl([M,M]_*^\times\bigr)\subseteq[M,M]^\times.$
Conversely, take any
$[f]\in[M,M]^\times.$
By surjectivity of $Q$, we may choose a pointed self-map
$f_0:(M,x_0)\longrightarrow(M,x_0)$
satisfying
$Q([f_0]_*)=[f].$
Since $[f]$ is a unit of $[M,M]$, $f_0$ is a homotopy equivalence in the free sense.
Thus $f_0$ induces isomorphisms on all homotopy groups. In particular,
\[
(f_0)_{*,1}:\pi_1(M,x_0)\longrightarrow\pi_1(M,x_0),
\]
\[
(f_0)_{*,2}:\pi_2(M,x_0)\longrightarrow\pi_2(M,x_0)
\]
are isomorphisms.
By Corollary~\ref{cor:pointed-homotopy-equivalences}, $f_0$ is a pointed homotopy equivalence, and therefore
$[f_0]_*\in[M,M]_*^\times.$
Thus
$[f]=Q([f_0]_*)\in Q\bigl([M,M]_*^\times\bigr).$
This proves the reverse inclusion.
Combining the two inclusions, we obtain
$Q\bigl([M,M]_*^\times\bigr)=[M,M]^\times.$

By \cref{cor:pointed-homotopy-equivalences}, the unit group of the monoid of pointed homotopy classes is
\[
[M,M]_*^\times=\left\{(A(m,k),d,0)\ \middle|\
m,d\in\{-1,1\},\ k\in\ZZ\right\}.
\]
By \cref{thm:free-classification}, the free homotopy class corresponding to any pointed self-map class in the double-reflection case
$(A(m,k),d,c)$
has a unique normal form
$(A(\mu,\kappa),D,c),$
where
\[
\mu=\operatorname{sgn}(d)m,
\]
\[
\kappa\equiv k\pmod2,\qquad\kappa\in\{0,1\},
\]
and
$D=|d|.$
For
$(A(m,k),d,0)\in[M,M]_*^\times,$
we have $m,d\in\{-1,1\}$, so
$\mu=\operatorname{sgn}(d)m=dm\in\{-1,1\},$
and
$D=|d|=1.$
Therefore, every element of $[M,M]^\times$ can be represented uniquely as
$Q\bigl((A(\mu,\kappa),1,0)\bigr),$
where
$\mu\in\{-1,1\},\qquad\kappa\in\{0,1\}.$
Conversely, for any $\mu,\kappa$ satisfying the above conditions, \cref{thm:classification} gives a pointed self-map class
$(A(\mu,\kappa),1,0).$
Since
$\det A(\mu,\kappa)=\mu=\pm1,$
by \cref{cor:pointed-homotopy-equivalences},
$(A(\mu,\kappa),1,0)\in[M,M]_*^\times.$
and hence
$Q\bigl((A(\mu,\kappa),1,0)\bigr)\in[M,M]^\times.$
Thus uniqueness of the normal forms gives a bijection of sets
$[M,M]^\times\cong\{-1,1\}\times\{0,1\}.$
In particular,
$|[M,M]^\times|=4.$

Take any
$Q\bigl((A(\mu,\kappa),1,0)\bigr)$
and
$Q\bigl((A(\nu,\lambda),1,0)\bigr)$
in $[M,M]^\times$.
Since $Q$ is a monoid homomorphism,
\[
\begin{aligned}
&Q\bigl((A(\mu,\kappa),1,0)\bigr)
\star Q\bigl((A(\nu,\lambda),1,0)\bigr)\\
&=Q\bigl((A(\mu,\kappa),1,0)\star(A(\nu,\lambda),1,0)\bigr).
\end{aligned}
\]
By the composition formula of \cref{cor:double-reflection-composition},
\[
\begin{aligned}
&(A(\mu,\kappa),1,0)\star(A(\nu,\lambda),1,0)\\
&=\bigl(A(\mu,\kappa)A(\nu,\lambda),1,0\bigr).
\end{aligned}
\]
Matrix multiplication gives
\[
\begin{aligned}
A(\mu,\kappa)A(\nu,\lambda)
&=\begin{pmatrix}\mu&\kappa\\0&1\end{pmatrix}
\begin{pmatrix}\nu&\lambda\\0&1\end{pmatrix}\\
&=\begin{pmatrix}\mu\nu&\kappa+\mu\lambda\\0&1\end{pmatrix}\\
&=A(\mu\nu,\kappa+\mu\lambda).
\end{aligned}
\]
Therefore
\[
\begin{aligned}
&Q\bigl((A(\mu,\kappa),1,0)\bigr)
\star Q\bigl((A(\nu,\lambda),1,0)\bigr)\\
&=Q\bigl((A(\mu\nu,\kappa+\mu\lambda),1,0)\bigr).
\end{aligned}
\]
Since the second coordinate $\kappa+\mu\lambda$ need not belong to $\{0,1\}$, we reduce the resulting homotopy class to normal form.
Let $\tau\in\{0,1\}$ be the unique integer satisfying
$\tau\equiv\kappa+\mu\lambda\pmod2.$
By \cref{thm:free-classification},
\[
Q\bigl((A(\mu\nu,\kappa+\mu\lambda),1,0)\bigr)
=Q\bigl((A(\mu\nu,\tau),1,0)\bigr).
\]
Since $\mu\in\{-1,1\}$, we have
$\mu\equiv1\pmod2.$
Consequently,
$\kappa+\mu\lambda\equiv\kappa+\lambda\pmod2.$
By the definition of $\tau$,
$\tau=\overline{\kappa+\lambda}.$
Thus, in normal-form coordinates, the group operation is
$(\mu,\kappa)\star(\nu,\lambda)=(\mu\nu,\overline{\kappa+\lambda}).$

Define the map
\[
\Theta:[M,M]^\times\longrightarrow\ZZ/2\ZZ\times\ZZ/2\ZZ
\]
by
\[
\Theta\bigl(Q((A(\mu,\kappa),1,0))\bigr)
=\left(\left[\frac{1-\mu}{2}\right]_2,[\kappa]_2\right).
\]
Since $\mu\in\{-1,1\}$, we have
\[
\frac{1-\mu}{2}=\begin{cases}0,&\mu=1,\\1,&\mu=-1,\end{cases}
\]
Hence the map
\[
\{-1,1\}\longrightarrow\ZZ/2\ZZ,\qquad
\mu\longmapsto\left[\frac{1-\mu}{2}\right]_2
\]
is a bijection.
Also, $\kappa\mapsto[\kappa]_2$ is a bijection from $\{0,1\}$ to $\ZZ/2\ZZ$. Together with the normal-form classification established above, this implies that $\Theta$ is a bijection of sets.

To prove that $\Theta$ preserves the group operation, note that
\[
\frac{1-\mu}{2}+\frac{1-\nu}{2}-\frac{1-\mu\nu}{2}
=\frac{(1-\mu)(1-\nu)}{2}.
\]
Since $1-\mu$ and $1-\nu$ both belong to $\{0,2\}$, the right-hand side is even. Hence
\[
\left[\frac{1-\mu\nu}{2}\right]_2
=\left[\frac{1-\mu}{2}\right]_2+
\left[\frac{1-\nu}{2}\right]_2.
\]
On the other hand,
$[\overline{\kappa+\lambda}]_2=[\kappa]_2+[\lambda]_2.$
By \eqref{eq:free-unit-product},
\[
\begin{aligned}
&\Theta\bigl(Q((A(\mu,\kappa),1,0))\star Q((A(\nu,\lambda),1,0))\bigr)\\
&=\Theta\bigl(Q((A(\mu\nu,\overline{\kappa+\lambda}),1,0))\bigr)\\
&=\left(\left[\frac{1-\mu\nu}{2}\right]_2,[\overline{\kappa+\lambda}]_2\right)\\
&=\left(\left[\frac{1-\mu}{2}\right]_2,[\kappa]_2\right)
+\left(\left[\frac{1-\nu}{2}\right]_2,[\lambda]_2\right)\\
&=\Theta\bigl(Q((A(\mu,\kappa),1,0))\bigr)
+\Theta\bigl(Q((A(\nu,\lambda),1,0))\bigr).
\end{aligned}
\]
Therefore $\Theta$ is a group homomorphism. Since $\Theta$ is also bijective, it is a group isomorphism, that is,
$[M,M]^\times\cong\ZZ/2\ZZ\times\ZZ/2\ZZ.$
\end{proof}

\begin{corollary}[Homotopy groups do not detect free homotopy classes]\label{cor:free-homotopy-groups-do-not-detect}
For every $\mathbf A\in\mathscr M$ and every odd integer $d$, the restriction of the forgetful map to $[M,M]_*^{\varphi_{\mathbf A},d}$,
\[
Q\big|_{[M,M]_*^{\varphi_{\mathbf A},d}}:
[M,M]_*^{\varphi_{\mathbf A},d}\longrightarrow[M,M]
\]
is injective. In particular,
$\left|Q\left([M,M]_*^{\varphi_{\mathbf A},d}\right)\right|=|d|.$
Furthermore, there exist $|d|$ pairwise non-freely-homotopic pointed self-maps
\[
f_0,f_1,\ldots,f_{|d|-1}:(M,x_0)\longrightarrow(M,x_0),
\]
such that for any
$i,j\in\{0,1,\ldots,|d|-1\}$
and every integer $n\ge1$, we have
\[
(f_i)_{*,n}=(f_j)_{*,n}:\pi_n(M,x_0)\longrightarrow\pi_n(M,x_0).
\]
Thus, when $|d|>1$, the induced homomorphisms on homotopy groups in all dimensions do not completely distinguish the free homotopy classes of self-maps of $M$.
\end{corollary}

\begin{proof}

Fix
$\mathbf A=A(m,k)\in\mathscr M$
and an odd integer $d$.
By Theorem~\ref{thm:classification}, the elements of the set
$[M,M]_*^{\varphi_{\mathbf A},d}$
are uniquely represented as
\[
(A(m,k),d,c),\qquad c\in\ZZ/|d|\ZZ.
\]
Take any
$[f_0]_*,[f_1]_*\in[M,M]_*^{\varphi_{\mathbf A},d},$
and write
\[
[f_0]_*=(A(m,k),d,c_0),
\]
\[
[f_1]_*=(A(m,k),d,c_1),
\]
where
$c_0,c_1\in\ZZ/|d|\ZZ.$
Suppose
$Q([f_0]_*)=Q([f_1]_*).$
By \cref{lem:free-orbit}, there exists $g\in G$ such that
$[f_1]_*=g\cdot[f_0]_*.$
Write $g$ uniquely as
$g=t^ra^\delta,$
where
$r\in\ZZ,\qquad\delta\in\{0,1\}.$
Put
$\varepsilon=(-1)^\delta.$
By \cref{lem:free-reflections}, the left action of the fundamental group $G$ on the double-reflection part is given by
\[
g\cdot(A(m,k),d,c_0)=\bigl(A(\varepsilon m,\varepsilon k+2r),\varepsilon d,c_0\bigr).
\]
Thus, from $[f_1]_*=g\cdot[f_0]_*$, we obtain
\[
\begin{aligned}
(A(m,k),d,c_1)&=g\cdot(A(m,k),d,c_0)\\
&=(A(\varepsilon m,\varepsilon k+2r),\varepsilon d,c_0).
\end{aligned}
\]
By uniqueness of the classification coordinates in \cref{thm:classification}, comparing the second coordinate on both sides gives
$d=\varepsilon d.$
Since $d$ is odd, $d\ne0$, and therefore
$\varepsilon=1.$
Since
$\varepsilon=(-1)^\delta,\qquad\delta\in\{0,1\},$
we obtain
$\delta=0.$
Comparing the matrix coordinates on both sides gives
$A(m,k)=A(m,k+2r).$
By the definition of the matrix $A(m,k)$,
\[
\begin{pmatrix}m&k\\0&1\end{pmatrix}
=
\begin{pmatrix}m&k+2r\\0&1\end{pmatrix}.
\]
Hence
$k=k+2r,$
that is,
$r=0.$
Therefore
$g=t^0a^0=1_G.$
Consequently,
\[
[f_1]_* =g\cdot[f_0]_* =1_G\cdot[f_0]_*=[f_0]_*.
\]
In particular, by uniqueness of the classification coordinates,
$c_0=c_1.$
This proves that
$Q\big|_{[M,M]_*^{\varphi_{\mathbf A},d}}$
is injective.

By \cref{thm:classification}, there is a bijection of sets
$[M,M]_*^{\varphi_{\mathbf A},d}\cong\ZZ/|d|\ZZ.$
Therefore
$\left|[M,M]_*^{\varphi_{\mathbf A},d}\right|=|d|.$
On the other hand, by the injectivity established above,
$Q\big|_{[M,M]_*^{\varphi_{\mathbf A},d}}$
is injective and thus induces a bijection of sets
\[
[M,M]_*^{\varphi_{\mathbf A},d}
\xrightarrow{\ \cong\ }Q\left([M,M]_*^{\varphi_{\mathbf A},d}\right).
\]
Consequently,
$Q\left([M,M]_*^{\varphi_{\mathbf A},d}\right)\cong\ZZ/|d|\ZZ,$
and in particular,
$\left|Q\left([M,M]_*^{\varphi_{\mathbf A},d}\right)\right|=|d|.$

For every integer
$i\in\{0,1,\ldots,|d|-1\},$
write $[i]_{|d|}$ for the residue class of $i$ in the cyclic group
$\ZZ/|d|\ZZ.$
By Theorem~\ref{thm:classification}, for each such $i$ there is a unique pointed homotopy class
$[f_i]_*=(A(m,k),d,[i]_{|d|}).$
For each $[f_i]_*$, choose a pointed representative map
$f_i:(M,x_0)\longrightarrow(M,x_0).$
Take any distinct integers
$i,j\in\{0,1,\ldots,|d|-1\}.$
Since $i\ne j$ and $i,j$ belong to a complete system of residues modulo $|d|$, we have
$[i]_{|d|}\ne[j]_{|d|}.$
By uniqueness of the classification coordinates in \cref{thm:classification},
$[f_i]_*\ne[f_j]_*.$
By the injectivity established above,
$Q([f_i]_*)\ne Q([f_j]_*),$
that is,
$[f_i]\ne[f_j]\in[M,M].$
Thus $f_i$ and $f_j$ are not freely homotopic.
Since $i,j$ were arbitrary distinct indices, the maps
$f_0,f_1,\ldots,f_{|d|-1}$
are pairwise not freely homotopic.

On the other hand, since for every $i$ we have
$[f_i]_*\in[M,M]_*^{\varphi_{\mathbf A},d},$
by \cref{cor:homotopy-groups-do-not-detect}, for any
$i,j\in\{0,1,\ldots,|d|-1\}$
and any integer $n\ge1$, we have
\[
(f_i)_{*,n}=(f_j)_{*,n}:\pi_n(M,x_0)\longrightarrow\pi_n(M,x_0).
\]
In particular, on the first and second homotopy groups respectively, we have
$(f_i)_{*,1}=(f_j)_{*,1}=\varphi_{\mathbf A}$
and (using the previously fixed isomorphism $\pi_2(M,x_0)\cong\ZZ$)
$(f_i)_{*,2}=(f_j)_{*,2}=d\,\id_{\ZZ}.$
In conclusion, there are exactly $|d|$ free homotopy classes realized by the set $[M,M]_*^{\varphi_{\mathbf A},d}$, and they have identical induced homomorphisms on homotopy groups in all dimensions.

When $|d|>1$, these free homotopy classes are distinct, but their pointed representative maps induce the same homomorphisms on all homotopy groups.

Therefore, even knowing the homomorphisms induced by a pointed self-map on every homotopy group does not completely determine its free homotopy class.
\end{proof}

\appendix
\section{Proofs for the geometric model}\label{app:geometry}

\subsection{Proof of \cref{prop:punctured-model}}\label{app:proof-prop:punctured-model}
\begin{proof}Write $\RP^3$ as
\(D^3/(x\sim -x\text{ for }x\in\partial D^3).\)
Remove a small open three-ball centered at the origin in the interior of $D^3$. The remaining shell is homeomorphic to
\(S^2\times[0,1].\)
The inner boundary $S^2\times\{0\}$ is the sphere used in forming the connected sum; the outer boundary $S^2\times\{1\}$ is still identified by
\((x,1)\sim(-x,1)\)
Thus
\[
N:=\RP^3\setminus\operatorname{int}(D^3)
\cong
\frac{S^2\times[0,1]}{(x,1)\sim(-x,1)}.
\]
Contraction in the $[0,1]$ direction deformation retracts $N$ onto $S^2/(x\sim -x)=\RP^2$. Hence
\(N\simeq\RP^2, \qquad \pi_1(N)\cong C_2.\)
Glue two copies $N_a,N_b$ along their spherical boundaries to obtain
\(M=N_a\cup_{S^2}N_b.\)
Reparametrizing the collars gives the useful quotient model
\begin{equation}\label{eq:connected-sum-quotient}
M\cong
\frac{S^2\times[0,1]}
{(x,0)\sim(-x,0),\ (x,1)\sim(-x,1)}.
\end{equation}
\end{proof}

\subsection{Proof of \cref{prop:fundamental-group}}\label{app:proof-prop:fundamental-group}
\begin{proof}Choose open sets $U,V\subset M$ by adjoining small collars to $N_a,N_b$, respectively, so that
\[
U\simeq N_a\simeq\RP^2,
\qquad
V\simeq N_b\simeq\RP^2.
\]
Their intersection deformation retracts onto the gluing sphere:
\(U\cap V\simeq S^2.\)
Choose the basepoint in $U\cap V$. Then
\[
\pi_1(U)\cong C_2,
\qquad
\pi_1(V)\cong C_2,
\qquad
\pi_1(U\cap V)=1.
\]
By the Seifert--van Kampen theorem\cite[Chapter~1, \S1.2]{hatcher},
\[
\pi_1(M)
\cong
\pi_1(U)*_{\pi_1(U\cap V)}\pi_1(V)
\cong C_2*C_2.
\]
Writing $a,b$ for the nonidentity elements of the two $C_2$ factors gives
\begin{equation}\label{eq:G-presentation}
G=\pi_1(M,x_0)
\cong
\langle a,b\mid a^2=b^2=1\rangle.
\end{equation}
This is a free product: there are no additional relations beyond $a^2=b^2=1$.
\end{proof}

\subsection{Proof of \cref{lem:proper-action}}\label{app:proof-lem:proper-action}
\begin{proof}
First suppose $t^n(x,u)=(x,u)$. The second coordinate gives $u-2n=u$, hence $n=0$. If
$t^na(x,u)=(x,u)$, the first coordinate gives $-x=x$, which is impossible on $S^2$. Thus no nonidentity element has a fixed point.

To prove proper discontinuity, take a compact subset $K\subset S^2\times\RR$. The continuous projection
\(\operatorname{pr}_2:S^2\times\RR\longrightarrow\RR\)
has compact image, so there is $R>0$ with
\(K\subset S^2\times[-R,R].\)
Each $g\in G$ is uniquely $t^n$ or $t^na$. If $t^nK\cap K\ne\varnothing$, there are
$u,v\in[-R,R]$ such that $v=u-2n$. Thus
\(2|n|=|u-v|\le 2R.\)
If $t^naK\cap K\ne\varnothing$, there are $u,v\in[-R,R]$ with $v=-u-2n$. Hence
\(2|n|=|u+v|\le 2R.\)
Consequently,
\begin{equation}\label{eq:finite-intersection-set}
F_K:=\{g\in G\mid gK\cap K\ne\varnothing\}
\end{equation}
is finite.

Fix $z\in\widetilde M$. Since $S^2\times\RR$ is locally compact and Hausdorff, choose an open neighborhood
$V$ of $z$ with compact closure $K:=\overline V$, and put $F=F_K$. Freeness implies
$gz\ne z$ for every $g\in F\setminus\{1\}$. The Hausdorff property supplies disjoint open neighborhoods
$A_g$ of $z$ and $B_g$ of $gz$. Set
\(U_g=A_g\cap g^{-1}(B_g).\)
Then $U_g\cap gU_g=\varnothing$. Since $F\setminus\{1\}$ is finite,
\(U:=V\cap\bigcap_{g\in F\setminus\{1\}}U_g\)
is still an open neighborhood of $z$. If $g\in F\setminus\{1\}$, then $U\subset U_g$, so
$U\cap gU=\varnothing$. If $g\notin F$, then $U\subset K$ and $gK\cap K=\varnothing$,
which again implies $U\cap gU=\varnothing$. This proves the assertion.
\end{proof}

\subsection{Proof of \cref{thm:universal-cover}}\label{app:proof-thm:universal-cover}
\begin{proof}
Every $G$-orbit meets $D=S^2\times[0,1]$, and its intersection with
$S^2\times(0,1)$ contains at most one point. At the boundaries $u=0,1$ we have, respectively,
\(a(x,0)=(-x,0), \qquad b(x,1)=(-x,1).\)
Thus
\[
(S^2\times\RR)/G
\cong
\frac{S^2\times[0,1]}
{(x,0)\sim(-x,0),\ (x,1)\sim(-x,1)}
\cong M
\]
where the last homeomorphism is \eqref{eq:connected-sum-quotient}. By \cref{lem:proper-action},
the action is free and properly discontinuous, so the quotient map is a covering map. Since
$S^2\times\RR\simeq S^2$ is simply connected, $p$ is the universal covering and its deck group
is $G$.

A covering induces isomorphisms on homotopy groups in degrees $k\ge2$, giving
\[
\pi_2(M)\cong\pi_2(S^2\times\RR)\cong\pi_2(S^2)\cong\ZZ,
\]
and
\[
\pi_3(M)\cong\pi_3(S^2\times\RR)\cong\pi_3(S^2)\cong\ZZ.
\]
By \eqref{eq:deck-tn}--\eqref{eq:deck-tna}, $t^n$ acts identically on the sphere coordinate,
whereas $t^na$ acts by the antipodal map
\(A:S^2\to S^2, \qquad A(x)=-x.\)
Since $\deg A=-1$, translations act by $+1$ on $\pi_2$ and reflections by $-1$. This is precisely
the action described by $\rho$.
\end{proof}

\section{Proofs on equivariant lifts and evaluation fibrations}\label{app:lifts}

\subsection{Proof of \cref{lem:equivariance-pi1}}\label{app:proof-lem:equivariance-pi1}
\begin{proof}
We first prove the correspondence on maps. Let
\(f:(M,x_0)\longrightarrow(M,x_0)\)
be pointed with $f_*=\varphi$. Since $\widetilde M=S^2\times\RR$ is simply connected,
the covering lifting criterion supplies a lift of $fp:\widetilde M\to M$. Prescribing its value at
$\widetilde x_0$ gives a unique map
\[
\widetilde f:(\widetilde M,\widetilde x_0)
\longrightarrow(\widetilde M,\widetilde x_0)
\]
satisfying
\begin{equation}\label{eq:lift-f}
p\widetilde f=fp.
\end{equation}
We prove step by step that
\begin{equation}\label{eq:phi-equivariance}
\widetilde f(gz)=\varphi(g)\widetilde f(z)
\qquad(g\in G,\ z\in\widetilde M).
\end{equation}
Fix $g\in G$ and a pointed loop $\alpha:I\to M$ representing it. Let
$\widetilde\alpha$ be its unique lift starting at $\widetilde x_0$. By
\eqref{eq:path-end-deck}, $\widetilde\alpha(1)=g\widetilde x_0$. The path
$\widetilde f\circ\widetilde\alpha$ satisfies
\[
\begin{aligned}
p\bigl(\widetilde f(\widetilde\alpha(t))\bigr)
&=(p\widetilde f)(\widetilde\alpha(t))\\
&=(fp)(\widetilde\alpha(t))\\
&=f(\alpha(t)),
\end{aligned}
\]
and
\[
(\widetilde f\circ\widetilde\alpha)(0)
=\widetilde f(\widetilde x_0)=\widetilde x_0.
\]
It is therefore the unique lift of $f\circ\alpha$ starting at $\widetilde x_0$. Since
\([f\circ\alpha]=f_*[\alpha]=\varphi(g),\)
the identification with deck transformations says that its endpoint is
$\varphi(g)\widetilde x_0$. On the other hand, its endpoint is
$\widetilde f(\widetilde\alpha(1))=\widetilde f(g\widetilde x_0)$. Hence
\begin{equation}\label{eq:lift-equivariance-basepoint}
\widetilde f(g\widetilde x_0)=\varphi(g)\widetilde x_0.
\end{equation}
Define
\(F_1(z)=\widetilde f(gz),\qquad F_2(z)=\varphi(g)\widetilde f(z).\)
Since $g$ and $\varphi(g)$ are deck transformations,
\[
\begin{aligned}
pF_1(z)&=p\widetilde f(gz)=f(p(gz))=f(p(z)),\\
pF_2(z)&=p\bigl(\varphi(g)\widetilde f(z)\bigr)
=p\widetilde f(z)=f(p(z)).
\end{aligned}
\]
Thus $F_1,F_2$ both lift $fp$. Moreover, by
\eqref{eq:lift-equivariance-basepoint}, 
\[
F_1(\widetilde x_0)=\varphi(g)\widetilde x_0
=\varphi(g)\widetilde f(\widetilde x_0)=F_2(\widetilde x_0).
\]
Connectedness of $\widetilde M$ and uniqueness of lifts imply $F_1=F_2$, proving
\eqref{eq:phi-equivariance}.

Conversely, suppose a continuous map
\(F:\widetilde M\longrightarrow\widetilde M\)
satisfies
\[F(gz)=\varphi(g)F(z),\qquad F(\widetilde x_0)=\widetilde x_0.\]
Equivariance and the deck transformation property of $\varphi(g)$ imply
\(pF(gz)=p\bigl(\varphi(g)F(z)\bigr)=pF(z).\)
Thus $pF$ is constant on $G$-orbits. By the universal property of the quotient
$p:\widetilde M\to\widetilde M/G\cong M$, there is a unique continuous map
$f:M\to M$ satisfying
\begin{equation}\label{eq:descend-F}
pF=fp.
\end{equation}
Furthermore,
\[
f(x_0)=f(p(\widetilde x_0))=pF(\widetilde x_0)
=p(\widetilde x_0)=x_0.
\]
We check that its induced fundamental-group homomorphism is exactly $\varphi$. For
$g=[\alpha]\in G$, let $\widetilde\alpha$ be the lift of $\alpha$ starting at
$\widetilde x_0$. Then
$\widetilde\alpha(1)=g\widetilde x_0$. By \eqref{eq:descend-F},
\[
p\circ F\circ\widetilde\alpha
=f\circ p\circ\widetilde\alpha=f\circ\alpha,
\qquad
(F\circ\widetilde\alpha)(0)=F(\widetilde x_0)=\widetilde x_0.
\]
Thus $F\circ\widetilde\alpha$ is the unique lift of $f\circ\alpha$ starting at
$\widetilde x_0$, and its endpoint is
\[
\begin{aligned}
(F\circ\widetilde\alpha)(1)
&=F(g\widetilde x_0)\\
&=\varphi(g)F(\widetilde x_0)\\
&=\varphi(g)\widetilde x_0.
\end{aligned}
\]
The identification of the fundamental group with deck transformations gives
\([f\circ\alpha]=\varphi(g).\)
and hence
\(f_*(g)=f_*[\alpha]=[f\circ\alpha]=\varphi(g).\)
Since $g$ was arbitrary, $f_*=\varphi$. The condition
$F(\widetilde x_0)=\widetilde x_0$ is essential here: fixing the lifted basepoint makes the induced homomorphism
exactly $\varphi$, rather than determining it only up to inner conjugacy.

Starting with $f$, uniqueness of descent recovers $f$ after lifting and descending. Starting with $F$,
descending and then taking the pointed lift recovers $F$ by uniqueness, since $F$ itself is such a lift.
Thus the two constructions are mutually inverse on maps.

We now check compatibility with homotopies. Let
\[
f_0,f_1:(M,x_0)\longrightarrow(M,x_0),\qquad
(f_0)_*=(f_1)_*=\varphi,
\]
and let $F_0,F_1:\widetilde M\to\widetilde M$ be the corresponding unique lifts, so that
\[
pF_i=f_ip,\qquad
F_i(\widetilde x_0)=\widetilde x_0,\qquad
F_i(gz)=\varphi(g)F_i(z)\quad(i=0,1).
\]
Suppose $f_0\simeq_*f_1$, and choose a pointed homotopy
\(H:M\times I\longrightarrow M\)
satisfying
\[
H(x,0)=f_0(x),\qquad H(x,1)=f_1(x),\qquad H(x_0,s)=x_0.
\]
Define
\[
\overline H:\widetilde M\times I\longrightarrow M,\qquad
\overline H(z,s)=H(p(z),s).
\]
At $s=0$,
\(\overline H(z,0)=f_0(p(z))=p(F_0(z)).\)
so $F_0$ lifts $\overline H|_{\widetilde M\times\{0\}}$. The homotopy lifting property of a covering
gives a unique continuous map
\(\widetilde H:\widetilde M\times I\longrightarrow\widetilde M\)
satisfying
\begin{equation}\label{eq:lifted-homotopy}
p\widetilde H(z,s)=H(p(z),s),\qquad
\widetilde H(z,0)=F_0(z).
\end{equation}

First we show that this lift fixes $\widetilde x_0$ throughout the homotopy. Consider the path
\(\gamma(s)=\widetilde H(\widetilde x_0,s).\)
By \eqref{eq:lifted-homotopy} and $H(x_0,s)=x_0$,
\(p(\gamma(s))=H(p(\widetilde x_0),s)=H(x_0,s)=x_0.\)
and
\(\gamma(0)=F_0(\widetilde x_0)=\widetilde x_0.\)
Thus $\gamma$ lifts the constant path at $x_0$, starting at $\widetilde x_0$.
Uniqueness of path lifting gives
\begin{equation}\label{eq:lifted-homotopy-basepoint}
\widetilde H(\widetilde x_0,s)=\widetilde x_0\qquad(s\in I).
\end{equation}

Next we show $\widetilde H(-,1)=F_1$. By \eqref{eq:lifted-homotopy},
\(p\widetilde H(z,1)=H(p(z),1)=f_1(p(z)).\)
so $\widetilde H(-,1)$ lifts $f_1p$. By
\eqref{eq:lifted-homotopy-basepoint}, its value at $\widetilde x_0$ is
$\widetilde x_0$. The map $F_1$ also lifts $f_1p$, with
$F_1(\widetilde x_0)=\widetilde x_0$. Uniqueness implies
\(\widetilde H(-,1)=F_1.\)

Finally, we prove equivariance at every time. Fix $g\in G$ and define
\[
A_g(z,s)=\widetilde H(gz,s),\qquad
B_g(z,s)=\varphi(g)\widetilde H(z,s).
\]
Since $g$ and $\varphi(g)$ are deck transformations,
\[
\begin{aligned}
pA_g(z,s)&=H(p(gz),s)=H(p(z),s),\\
pB_g(z,s)&=p\widetilde H(z,s)=H(p(z),s).
\end{aligned}
\]
both $A_g$ and $B_g$ lift $\overline H$. At
$(\widetilde x_0,0)$, equivariance of $F_0$ gives
\[
\begin{aligned}
A_g(\widetilde x_0,0)
&=F_0(g\widetilde x_0)\\
&=\varphi(g)F_0(\widetilde x_0)\\
&=\varphi(g)\widetilde x_0
=B_g(\widetilde x_0,0).
\end{aligned}
\]
Connectedness of $\widetilde M\times I$ and uniqueness of lifts imply $A_g=B_g$. Hence
\[
\widetilde H(gz,s)=\varphi(g)\widetilde H(z,s)
\qquad(g\in G,\ z\in\widetilde M,\ s\in I).
\]
Together with the endpoint and basepoint conditions, this proves that $\widetilde H$ is a
$\varphi$-equivariant homotopy from $F_0$ to $F_1$ fixing $\widetilde x_0$.

Conversely, let $\widetilde H:\widetilde M\times I\to\widetilde M$ be a
$\varphi$-equivariant homotopy between $F_0,F_1$ fixing $\widetilde x_0$. Then
\[
p\widetilde H(gz,s)
=p\bigl(\varphi(g)\widetilde H(z,s)\bigr)
=p\widetilde H(z,s).
\]
Thus $p\widetilde H$ is constant on the orbits of the action on the first factor.
The covering $p$ is an open quotient map, as is $p\times\id_I$. Its universal property supplies
a unique continuous map
\[
H:M\times I\longrightarrow M
\]
satisfying
\[
H(p(z),s)=p\widetilde H(z,s).
\]
At $s=0$,
\(H(p(z),0)=pF_0(z)=f_0(p(z));\)
and at $s=1$,
\(H(p(z),1)=pF_1(z)=f_1(p(z)).\)
Surjectivity of $p$ yields $H(-,0)=f_0$ and $H(-,1)=f_1$. Moreover,
\[
H(x_0,s)=H(p(\widetilde x_0),s)
=p\widetilde H(\widetilde x_0,s)=p(\widetilde x_0)=x_0.
\]
so $H$ is pointed. We have proved
\[
f_0\simeq_*f_1
\quad\Longleftrightarrow\quad
F_0\simeq_{\varphi,\widetilde x_0}F_1.
\]
The mutually inverse constructions preserve these homotopy relations, and hence induce
the asserted natural bijection on homotopy classes.
\end{proof}

\subsection{Proof of \cref{lem:section-evaluation-fibration}}\label{app:proof-lem:section-evaluation-fibration}
\begin{proof}
First verify the $G$-cofibration condition. Choose a CW structure on $M$ with $x_0$ a vertex.
The closed three-manifold $M=\RP^3\#\RP^3$ admits a triangulation having the prescribed point as a vertex,
and hence such a CW structure. Lifting characteristic maps along $p$ produces a CW structure on $X$.
Deck transformations permute the lifts of each cell.
Since the action is free, the lifted cells occur as $G\times D^n$,
attached along $G\times S^{n-1}$. Thus $X$ is a free $G$-CW complex.
The lifts of $x_0$ form exactly $A$, a zero-dimensional $G$-CW subcomplex.

The inclusion of a $G$-CW subcomplex has the equivariant homotopy extension property,
and is therefore a $G$-cofibration\cite[Chapter~I, \S\S2--3, pp.~16--19]{may-equivariant}.
This can also be seen cell by cell: on a representative cell $D^n$, the homotopy is already defined on
\(D^n\times\{0\}\ \cup\ S^{n-1}\times I\)
Extend over $D^n\times I$ using the standard cylinder retraction, then define it on all translates by equivariance.
The definitions agree on cell boundaries; the weak topology of the CW complex ensures continuity.
Thus $A\hookrightarrow X$ has the required equivariant homotopy extension property.

To handle arbitrary parameters simultaneously, put
\(C=(X\times\{0\})\cup(A\times I).\)
The homotopy extension property supplies an equivariant retraction
\[
R:X\times I\longrightarrow C,\qquad R|_C=\id_C.
\]
Explicitly, use target $C$, initial map $x\mapsto(x,0)$,
and homotopy $(a,t)\mapsto(a,t)$ on $A$; an extension is the required $R$.
For any parameter space $Z$ with trivial action, $\id_Z\times R$ is again an equivariant retraction,
giving extensions with parameters.

Consider an arbitrary homotopy lifting problem for $r$:
\[
\begin{tikzcd}[column sep=5em,row sep=3em]
Z\times\{0\} \arrow[r,"u_0"] \arrow[d,hook]
 & \Map_G(X,Y) \arrow[d,"r"]\\
Z\times I \arrow[r,"v"'] \arrow[ur,dashed,"\widetilde v"]
 & \Map_G(A,Y).
\end{tikzcd}
\]
Commutativity means $r(u_0(z))=v(z,0)$. By exponential adjunction, $u_0$ and $v$
correspond to continuous equivariant maps
\[
U_0:Z\times X\to Y,\qquad U_0(z,x)=u_0(z)(x),
\]
\[
K:Z\times A\times I\to Y,\qquad K(z,a,t)=v(z,t)(a).
\]
and commutativity becomes
\(K(z,a,0)=U_0(z,a).\)
We can therefore define on $Z\times C$
\(L(z,x,0)=U_0(z,x),\qquad L(z,a,t)=K(z,a,t).\)
The definitions agree on $Z\times A\times\{0\}$, so the closed gluing lemma gives continuity.
Since $U_0$ and $K$ are equivariant, so is $L$.
Define
\(\widehat L(z,x,t):=L\bigl(z,R(x,t)\bigr).\)
This map is continuous and equivariant, and satisfies
\(\widehat L(z,x,0)=U_0(z,x),\qquad \widehat L(z,a,t)=K(z,a,t),\)
since $R$ restricts to the identity on $C$. The extension is depicted by
\[
\begin{tikzcd}[column sep=6em,row sep=3em]
Z\times C \arrow[r,"L"] \arrow[d,hook]
 & Y\\
Z\times X\times I \arrow[ur,dashed,"\widehat L"']
 & {}
\end{tikzcd}
\]
Applying exponential adjunction again, put
\(\widetilde v(z,t)(x)=\widehat L(z,x,t).\)
Then $\widetilde v:Z\times I\to\Map_G(X,Y)$ is continuous and
$\widetilde v(z,0)=u_0(z)$, $r(\widetilde v(z,t))=v(z,t)$.
This solves every lifting problem in the first diagram, proving that $r$ is a Hurewicz fibration.

We now identify the base. For a discrete group $G$, subgroup $H\le G$, and $G$-space $Y$,
evaluation at the identity coset gives
\(\Map_G(G/H,Y)\cong Y^H;\)
The inverse sends $y\in Y^H$ to $gH\mapsto gy$.
The condition $y\in Y^H$ ensures independence of the coset representative, while equivariance determines the map
by its value at the identity coset\cite[Definition~7.1, p.~249]{davis-luck}.

In our situation, freeness of the deck action gives
\(g\widetilde x_0=\widetilde x_0\ \Longrightarrow\ g=1.\)
Thus $A$ is a free orbit, equivariantly homeomorphic to the discrete $G$-set $G/\{1\}$.
For $y\in\widetilde M$, define
\(h_y(g\widetilde x_0)=\varphi(g)y.\)
If $g\widetilde x_0=g'\widetilde x_0$, freeness implies $g=g'$, so the definition is unambiguous.
For $k\in G$,
\[
h_y(kg\widetilde x_0)=\varphi(kg)y
=\varphi(k)\varphi(g)y=\varphi(k)h_y(g\widetilde x_0),
\]
so $h_y$ is equivariant. It is continuous because $A$ is discrete.
Conversely, every equivariant $h:A\to Y$ satisfies
\(h(g\widetilde x_0)=\varphi(g)h(\widetilde x_0),\)
hence $h=h_{h(\widetilde x_0)}$. Since $h_y(\widetilde x_0)=y$, the maps $e$ and
$y\mapsto h_y$ are mutually inverse as set maps.

To check the topologies, note that $e$ is evaluation at a point and is continuous for the compact-open topology.
Compact subsets of the discrete space $A$ are finite, so the compact-open topology on $\Map(A,Y)$
is the product topology on $\prod_{a\in A}Y$.
Each coordinate of $y\mapsto h_y$, indexed by $g\widetilde x_0$, is the continuous map
$y\mapsto\varphi(g)y$. Hence the map into the product, and into the equivariant subspace, is continuous.
The same holds after compact generation. Thus $e$ is a homeomorphism.

Since $(e\circ r)(F)=F(\widetilde x_0)$, there is a strictly commutative diagram
\[
\begin{tikzcd}[column sep=5em,row sep=3em]
\mathcal E_\varphi \arrow[r,"r"] \arrow[dr,"\operatorname{ev}_{\widetilde x_0}"']
 & \Map_G(A,Y) \arrow[d,"e","\cong"']\\
 & \widetilde M.
\end{tikzcd}
\]
Composing a fibration with a homeomorphism of its base gives a fibration. Its inverse image over $\widetilde x_0$ is
\[
\{F\in\mathcal E_\varphi\mid F(\widetilde x_0)=\widetilde x_0\}
=\mathcal F_\varphi.
\]
This proves \eqref{eq:equivariant-evaluation-fibration}.
\end{proof}

\subsection{Proof of \cref{thm:remove-lift-basepoint}}\label{app:proof-thm:remove-lift-basepoint}
\begin{proof}
Write
\[
e=\operatorname{ev}_{\widetilde x_0}:
\mathcal E_\varphi\longrightarrow\widetilde M.
\]
By \cref{lem:section-evaluation-fibration},
$e$ is a Hurewicz fibration, and
\(\mathcal F_\varphi=e^{-1}(\widetilde x_0).\)

If $\mathcal E_\varphi=\varnothing$, then
$\mathcal F_\varphi=\varnothing$, and the conclusion follows.
Henceforth assume $\mathcal E_\varphi\neq\varnothing$.

We first prove that $j_*$ is surjective.
Let $F\in\mathcal E_\varphi$.
Since $\widetilde M$ is path connected, there exists a path
\[
\sigma:I\longrightarrow\widetilde M,
\qquad
\sigma(0)=e(F),\quad
\sigma(1)=\widetilde x_0.
\]
By the path lifting property, there exists
\[
\widetilde\sigma:I\longrightarrow\mathcal E_\varphi,
\qquad
\widetilde\sigma(0)=F,\quad
e\circ\widetilde\sigma=\sigma.
\]
Consequently,
\(e(\widetilde\sigma(1)) =\sigma(1)=\widetilde x_0,\)
so $\widetilde\sigma(1)\in\mathcal F_\varphi$.
Since $\widetilde\sigma$ joins $F$ to $\widetilde\sigma(1)$,
\[
[F]_{\mathcal E_\varphi}
=j_*[\widetilde\sigma(1)]_{\mathcal F_\varphi}.
\]
Thus $j_*$ is surjective.

We next prove that $j_*$ is injective.
Let $F_0,F_1\in\mathcal F_\varphi$, and suppose that
\(j_*[F_0]_{\mathcal F_\varphi} =j_*[F_1]_{\mathcal F_\varphi}.\)
There is therefore a path
\[
\Gamma:I\longrightarrow\mathcal E_\varphi,
\qquad
\Gamma(0)=F_0,\quad
\Gamma(1)=F_1.
\]
Set
\(\alpha=e\circ\Gamma:I\longrightarrow\widetilde M.\)
Since $F_0,F_1\in e^{-1}(\widetilde x_0)$,
\(\alpha(0)=\alpha(1)=\widetilde x_0.\)
Moreover,
\[
\widetilde M=S^2\times\mathbb R\simeq S^2,
\qquad
\pi_1(\widetilde M,\widetilde x_0)=0.
\]
Hence $\alpha$ admits a nullhomotopy relative to its endpoints:
there exists a continuous map
\(K:I\times I\longrightarrow\widetilde M\)
such that
\[
\begin{aligned}
K(s,0)&=\alpha(s),&
K(s,1)&=\widetilde x_0,\\
K(0,t)&=\widetilde x_0,&
K(1,t)&=\widetilde x_0.
\end{aligned}
\]

Put
\(C=(I\times\{0\})\cup(\{0,1\}\times I),\)
and define
\(L:C\longrightarrow\mathcal E_\varphi\)
by
\(L(s,0)=\Gamma(s),\qquad L(0,t)=F_0,\qquad L(1,t)=F_1.\)
Because $\Gamma(0)=F_0$ and $\Gamma(1)=F_1$,
these definitions agree on the intersections.
The closed pasting lemma therefore implies that $L$ is continuous.
Furthermore,
\(e\circ L=K|_C.\)

Since $(I,\{0,1\})$ is a CW pair,
the relative homotopy lifting property of the fibration $e$
provides a continuous map
\[
\widetilde K:I\times I\longrightarrow\mathcal E_\varphi
\]
satisfying
\(e\circ\widetilde K=K, \qquad \widetilde K|_C=L.\)
Thus the lift preserves both the prescribed path $\Gamma$
on the bottom edge and the constant paths on the two side edges.

Define
\(\Lambda(s)=\widetilde K(s,1).\)
For every $s\in I$,
\(e(\Lambda(s))=K(s,1)=\widetilde x_0.\)
Hence $\Lambda:I\to\mathcal F_\varphi$ is a continuous path.
Its endpoints are
\[
\Lambda(0)=\widetilde K(0,1)=F_0,
\qquad
\Lambda(1)=\widetilde K(1,1)=F_1.
\]
Therefore
\([F_0]_{\mathcal F_\varphi} =[F_1]_{\mathcal F_\varphi},\)
which proves that $j_*$ is injective.

Combining surjectivity and injectivity yields
\[
j_*:\pi_0(\mathcal F_\varphi)
\xrightarrow{\cong}
\pi_0(\mathcal E_\varphi).
\]
Composing with the natural bijection of
\cref{lem:equivariance-pi1}, we obtain
\[
\boxed{
[M,M]_*^\varphi
\cong\pi_0(\mathcal F_\varphi)
\cong\pi_0(\mathcal E_\varphi).
}
\qedhere
\]
\end{proof}

\section{Proof of the path double-coset lemma}\label{sec:double-coset-appendix}

\subsection{Proof of \cref{lem:double-coset-path}}\label{app:proof-lem:double-coset-path}
\begin{proof}
Write $\mathcal H=\mathcal H(A,C;B)$ with basepoint $(a_0,\omega,c_0)$.
The path model gives $\mathcal H=A\times_B^h C$. The homotopy Mayer--Vietoris sequence
\cite[\S1.5.3, p.~79; Exercise~1.7.22, pp.~89--90]{malkiewich}, 
has the relevant portion
\[
\begin{gathered}
\cdots\longrightarrow\pi_2(\mathcal H)
\longrightarrow\pi_2(A)\oplus\pi_2(C)
\xrightarrow{\ i_*-\tau_\omega j_*\ }\pi_2(B)\\
\longrightarrow\pi_1(\mathcal H)
\longrightarrow\pi_1(A)\times\pi_1(C)
\xrightarrow{\delta}\pi_1(B)
\xrightarrow{q}\pi_0(\mathcal H)\\
\longrightarrow\pi_0(A)\times\pi_0(C)=\{*\}.
\end{gathered}
\]
Here $\tau_\omega $ is basepoint transport along the reference path; the final terms are interpreted through pointed sets and group actions.
Explicitly, identify components of the fixed-endpoint path space with
$\pi_1(B,i(a_0))$ by $[\gamma]\mapsto[\gamma*\overline\omega]$. The endpoint loops $(u,v)$ act by
\(g\longmapsto i_*(u)^{-1}\,g\,\tau_\omega j_*(v).\)
Since $A,C$ are path connected, low-dimensional orbit exactness says that $q$ is surjective and its fibers are precisely these orbits.
Set $H=i_*\pi_1(A,a_0)$ and $K=\tau_\omega j_*\pi_1(C,c_0)$. The orbit of $g$ is
$HgK$, so
\[
\pi_0(\mathcal H)\cong H\backslash\pi_1(B,i(a_0))/K.
\]
If $\pi_1(B,i(a_0))$ is abelian, then $HgK=g+(H+K)$, giving
\[
\pi_0(\mathcal H)\cong \pi_1(B,i(a_0))/(H+K).\qedhere
\]
\end{proof}

\end{document}